\documentstyle [11pt,amstex,amscd,epsfig]{article}

\def\al{\alpha}
\def\b{\betta}

\def\vcup{\cup_{V}}

\def\bd{\partial}
\def\L{{\cal L}}

\def\dim{{\rm dim\; }}
\def\de{\delta}

\def\ga{\gamma}
\def\text#1{{\em #1}}

\def\la{\lambda}

\def\ep{\varepsilon}

\def\be{\begin{equation}}
\def\ee{\end{equation}}
\def\bear{\begin{eqnarray}}
\def\eear{\end{eqnarray}}
\def\best{\begin{eqnarray*}}
\def\eest{\end{eqnarray*}}
\def\pf{{\bf Proof}: }

\renewcommand{\theequation}{\arabic{section}.\arabic{equation}}

\newtheorem{theorem}{Theorem}[section]
\newtheorem{prop}[theorem]{Proposition}

\newtheorem{lemma}[theorem]{Lemma}
\newtheorem{cor}[theorem]{Corollary}

\newtheorem{defn}[theorem]{Definition}

\renewcommand{\thetheorem}{\arabic{section}.\arabic{theorem}}

\newtheorem{remark}[theorem]{Remark}
\newenvironment{rem}{\begin{remark}\rm}{\end{remark}}

\newtheorem{example}[theorem]{Example}
\newenvironment{ex}{\begin{example}\rm}{\end{example}}

\def\non{\noindent}
\def\pf{\non {\bf Proof. }}
\def\qed{\nopagebreak \hskip .1in { $\Box$ }\penalty10000 %
\hskip\parfillskip \par  }

\def\non{\noindent}
\def\pf{\non {\bf Proof. }}
\def\qed{\nopagebreak \hskip .1in { $\Box$ }\penalty10000 %
\hskip\parfillskip \par  }

\def\ra{\rightarrow}

\def\rg{\rangle}
\def\lg{\langle}
\def\r#1{\right#1}
\def\l#1{\left#1}

\def\ma#1{\mathop {#1} \limits}

\def\b{\beta}

\def\Si{\Sigma}

\def\ti{\times}

\def\Z{{ \Bbb Z}}
\def\R{{ \Bbb R}}
\def\P{{ \Bbb P}}
\def\Q{{ \Bbb Q}}
\def\cx{{ \Bbb C}}

\def\ev{\mbox{\rm ev}}

\def\wt#1{\widetilde{#1}}

\def\ov#1{\overline{#1}}

\def\w{\omega}

\def\M{{\cal M}}

\def\w{\omega}

\def\M{{\cal M}}

\def\U{{\cal U}}

\def\ker{\mbox{Ker }}
\def\cok{\mbox{Coker }}

\def\Mf{\M^{\small\em flat}}

\title{\bf The Symplectic Sum Formula for Gromov-Witten Invariants \vskip.2in}
\author{ Eleny-Nicoleta Ionel\thanks{both authors partially supported by
the N.S.F.} \\ University of Wisconsin\\
       Madison, WI  53706 \and Thomas H. Parker{\footnotesize *}\\ Michigan
State University\\ East Lansing,
MI   48824}
\date{\empty}

\addtocounter{section}{0}

\begin{document}

\maketitle

\vskip.15in

\vskip.4in

\begin{abstract}
In the symplectic category there is a `connect sum' operation that
glues symplectic manifolds
by identifying neighborhoods of  embedded codimension two
submanifolds. This paper establishes a
formula for the Gromov-Witten invariants of a symplectic sum $Z=X\#
Y$ in terms of the relative
GW invariants of
$X$ and $Y$.  Several applications to enumerative geometry are given.
\end{abstract}

\vskip.4in

    Gromov-Witten invariants are counts of holomorphic maps into
symplectic manifolds.  To define them on a
symplectic manifold $(X,\w)$ one introduces an almost complex structure
$J$ compatible with the symplectic
form $\w$ and forms the moduli space of $J$-holomorphic maps from
complex curves into $X$ and its
compactification, called the space of stable maps.  One then imposes
constraints on the stable maps,
      requiring the domain to have a certain form and the image to pass
through fixed homology
cycles in $X$. When the right number of constraints are imposed there
are only finitely many maps
satisfying the constraints; the  (oriented) count of these is the
corresponding GW invariant. For
complex algebraic manifolds these symplectic invariants can also be
defined by algebraic geometry, and
in important cases the invariants are the same as the curve counts
that are the subject of classical
enumerative algebraic geometry.

\smallskip

In the past decade  the foundations for this theory were  laid
and the invariants  were used to
solve several long-outstanding problems.   The focus now is on finding
effective ways of computing the
invariants.  One useful technique is the method  of `splitting the domain', in
which  one localizes the invariant to the set of
maps whose domain curves have two irreducible components  with the
constraints distributed
between them.   This produces recursion relations relating the
desired GW invariant to invariants with
lower degree or genus.   This paper establishes a general formula
describing the behavior of GW
invariants under the analogous operation of `splitting the target'.
Because we work in the context
of symplectic manifolds the natural splitting of the target is
the one associated with the symplectic
cut operation and its inverse, the symplectic sum.

\smallskip

The symplectic sum is defined by gluing along codimension two
submanifolds.  Specifically, let $X$ be a symplectic $2n$-manifold
with a symplectic $(2n-2)$-submanifold $V$.  Given a similar pair
$(Y,V)$ with a symplectic identification between the 2 copies of $V$
and a complex anti-linear isomorphism between the normal bundles
$N_X$ and $N_Y$ of
$V$ in $X$ and in $Y$ we can form the symplectic sum $Z=X\#_{V}Y$.
Our main theorem is a `Symplectic Sum Formula' which expresses the GW
invariants of the sum $Z$ in
terms of relative GW invariants of $(X,V)$ and $(Y,V)$ introduced in
\cite{ip4}.

\smallskip

The symplectic sum is perhaps more naturally seen  not as a single
manifold but as a family
depending on a `squeezing parameter'.  In Section 2 we construct a
family $Z\to D$ over
the disk whose fibers $Z_\la$ are smooth and symplectic for $\la\neq
0$ and whose central fiber
$Z_0$ is the singular manifold $X\cup_V Y$.     In a neighborhood of
$V$,  the total space $Z$ is
$N_X\oplus N_Y$, regarded as a subset of $X\ti Y$ by the symplectic
neighborhood theorem,
and the fiber $Z_\la$ is defined by the equation
$x y=\la$ where $x$ and $y$ are coordinates in the normal bundles
$N_X$ and $N_Y\cong
N_X^*$.   The fibration $Z\to D$ extends  away from $V$ as the
disjoint union of $X\times D$ and
$Y\ti D$.  The smooth fibers $Z_\la$, depicted in Figure 1, are
symplectically isotopic to one
another; each is a model of the symplectic sum.

\smallskip

The overall strategy for proving the symplectic sum formula is  to
relate the holomorphic
maps into $Z_0$  (which are simply  maps into $X$ and $Y$ which match
along $V$) with the
holomorphic maps into $Z_\la$ for $\la$ close to zero.   This
strategy involves two
parts: limits and gluing.  For the limiting process we consider
sequences of stable maps into the family $Z_\la$ of symplectic sums as
the `neck size' $\la\to 0$.  In Section 3 we show that these limit
to maps into the singular manifold $Z_0$ obtained by identifying $X$
and $Y$ along $V$.  Along the way several things become apparent.

\smallskip

First, the limit maps are holomorphic only if the almost complex
structures on $X$ and $Y$ match along $V$.  To ensure that we impose
the ``$V$-compatibility'' condition (1.10) on the almost complex
structure.  There is a price to pay for that.  In the symplectic
theory of Gromov-Witten invariants we are free to perturb $(J,\nu)$
without changing the invariant; that freedom can be used to ensure
that intersections are transverse.  After imposing the
$V$-compatibility condition we can no longer perturb $(J,\nu)$ along
$V$ at will, and hence we cannot assume that the limit curves are
transverse to $V$.  In fact, the components of the limit maps meet $V$
at points with   multiplicities and, worse, some components
may lie entirely in $V$.

    To count such maps into $Z_0$ we look first on the $X$ side, ignore
the maps with components in $V$, separate the moduli space of stable
maps into components $\M_s(X)$ labeled by the multiplicities
$s=(s_1,\dots, s_\ell)$ of their intersection points with $V$.  We
showed in \cite{ip4} how these spaces $\M_s(X)$ can be compactified
and used to define relative Gromov-Witten invariants $GW^V_X$.  The
definitions are briefly reviewed in Section 1.

\best
\begin{minipage}{2.5in} {\vskip-.9in Figure  1:\ \  Limiting curves
\\ in $Z_\la=X\#_\la Y$
as $\la\to 0$.}\end{minipage}
&\hskip.3in{\psfig{figure=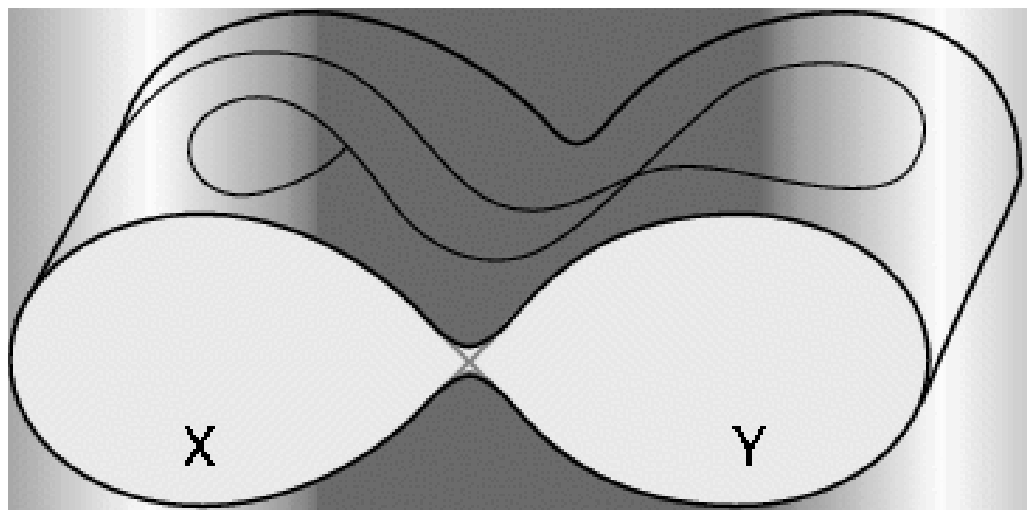,width=2.0in}}
\eest
\medskip

Second, as Figure 1 illustrates, connected curves in $Z_\la$ can limit
to curves whose restrictions to $X$ and $Y$ are not connected.  For
that reason the GW invariant, which counts stable curves from a
connected domain, is not the appropriate invariant for expressing a
sum formula.  Instead one should work with the `Taubes-Witten'
invariant $TW$, which counts stable maps from domains that need not be
connected.  Thus we seek a formula of the general form
\bear
TW_X^{V}\, *\, TW_Y^{V}\, = \, TW_Z
\label{0.1}
\eear
where $*$ is some operation that adds up the ways curves on the $X$
and $Y$ sides match and are identified with curves in $Z_\la$.  That
necessarily involves keeping track of the multiplicities $s$ and the homology
classes.  It also involves accounting for the limit maps with
non-trivial components in $V$; such curves are not counted by the
relative invariant and hence do not contribute to the left side of
(\ref{0.1}).  We postpone this issue by first analyzing limits of
curves which are {\em $\delta$-flat} in the sense of Definition
\ref{defnFlat}.

A more precise analysis reveals a third complication: the squeezing
process is not injective.  In Section 5 we again consider a sequence
of stable maps $f_n$ into $Z_\la$ as $\la\to 0$, this time focusing on
their behavior near $V$, where the $f_n$ do not uniformly converge.
We form renormalized maps $\hat{f}_n$ and prove that both the domains
and the images of the renormalized maps converge.  The images converge
nicely according to the leading order term of their Taylor expansions,
but the domains converge only after fixing certain roots of
unity.

These roots of unity are apparent as soon as one writes down
formulas.  Each stable map $f:C\to
Z_0$ decomposes into a pair of maps $f_1:C_1\ra X$ and $f_2:  C_2\ra
Y$ which agree at the nodes
of $C=C_1\cup C_2$.  For a specific example, suppose that $f$ is such 
a map that intersects
$V$ at a single point $p$ with
multiplicity three.  Then we can choose local coordinates $z$ on 
$C_1$ and $w$ on
$C_2$ centered at the node,
and coordinates $x$ on $X$ and
$y$ on $Y$ so that $f_1$ and $f_2$ have  expansions
$x(z)=az^{3}+\cdots $ and
$y(w)=bw^{3}+\cdots$.  To find  maps into $Z_\la$  near
$f$, we smooth the domain $C$ to the curve $C_\mu$ given locally near 
the node by $zw=\mu$ and
require that the
image of the smoothed map lie in $Z_\la$, which is locally  the locus of
$xy=\la$.  In fact,  the leading terms in the formulas for $f_1$ and
$f_2$ define a map
$F:C_\mu\to Z_\la$ whenever
$$
\la \, =\,  xy \, =\,  az^3 \cdot bw^3 \, =\,  ab\,(zw)^3  \, =\,    ab\, \mu^3
$$
and conversely any family of smooth maps with limit to $f$ satisfy
this equation in the limit (c.f.
Lemma \ref{lambdazlemma}).  Thus
$\la$ determines the domain $C_\mu$ up to a cube root of unity.  That
means that this particular $f$  is, at
least {\em a priori},   close to {\em three} smooth maps into $Z_\la$
--- a `cluster' of order
three.

\smallskip

Other maps $f$ into $Z_0$ have larger associated clusters  (the order
of the cluster is the product of the
multiplicities with which $f$ intersects $V$).  Within a cluster, the
maps have the same
leading order formula but have different smoothings of the domain.
As $\la\to 0$ the
maps within the cluster coalesce, limiting to the  single map $f$.

\medskip

    This clustering
phenomenon greatly  complicates the analysis.  To distinguish
the curves within each cluster  and make the analysis uniform in $\la$
as $\la\to 0$ it is necessary to use `rescaled' norms and
distances which magnify distances as the clusters form.  With  the
right choice of norms,
the distances between the maps within a cluster are bounded away
from zero as $\la\to 0$ and become the fiber of a covering of the
space of limit maps.   Sections
4-- 6 introduce  the required  norms, first on the space of
curves, then on the space of
maps.

\smallskip

    For maps we use
a  Sobolev norm  weighted in the directions perpendicular to $V$; the weights
are chosen so the norm dominates the $C^0$ distance between the
renormalized maps $\hat{f}$.
       On the space of curves
we require a stronger metric than the usual complete metrics on
$\ov{\M}_{g,n}$.  In section 4 we
define  a complete metric on $\ov{\M}_{g,n}\setminus {\cal
N}$ where ${\cal N}$ is the set of all nodal curves.  In this metric
the distance
between two sequences that approach ${\cal N}$ from different
directions (corresponding to the
roots of unity mentioned above) is bounded away from zero; thus  this
metric  separates the domain
curves of maps within a cluster.  The metric leads to a
compactification of $\ov{\M}_{g,n}\setminus {\cal N}$  in which the
stratum ${\cal N}_\ell$  of $\ell$-nodal curves is
replaced by a bundle over ${\cal N}_\ell$ whose fiber is the real torus
$T^\ell$.

\bigskip\bigskip

The limit process is reversed by constructing a space  of
approximately holomorphic
maps and showing it is diffeomorphic to the space of stable maps
into $Z_\la$.  The space of
approximate maps is described in Section 6, first intrinsically, then
as a subset ${\cal AM}_s$ of the space of maps.
For each $s$ and $\la$  it  is a covering of the
space $\M_s(Z_0)$ of the $\delta$-flat maps into $Z_0$ that meet $V$ at
       points with multiplicities $s$.  The fibers of this covering are
the clusters
-- they are distinct maps  into $Z_\la$ which converge to the same
limit as $\la\to 0$.

\smallskip

      From there the analysis follows  the   standard
technique that goes back to Taubes and Donaldson: correct the
approximate maps to true
holomorphic maps by  constructing a partial right
inverse to the linearization $D$ and applying  a
fixed point theorem.   That involves (a) showing that    the operator
$D^*D$ is  uniformly invertible as $\la\to 0$, and (b) proving
{\em a priori} that every solution
is close to an approximate solution, close enough to be in the domain
of the fixed
point theorem.  Proposition \ref{nearapprox} shows that (b) follows from the
renormalization analysis of Section 4.  But the  eigenvalue estimate
(a) proves to be
surprisingly  delicate and seems to succeed only with a very specific
choice of  norms.

\smallskip

The difficulty, of course, is that $Z_\la$ becomes singular along $V$
as $\la\to 0$.  However,
for small $\la$ the bisectional curvature in the neck region is
negative; a Bochner
formula
then shows that eigenfunctions with small eigenvalue cannot be
concentrating in the neck.  One
can then reason that since the cokernel of $D$ vanishes on
$Z_0$ (for generic $J$) it should also vanish on $Z_\la$ for small
$\la$.  We  make that
reasoning rigorous by  introducing exponential weight functions into
the norms, thereby making
the linearizations $D_\la$  a continuous family of Fredholm maps.
That in turn necessities
further work on the  Bochner formula, bounding the  additional term
   that arises  from the derivative of the weight functions.
These estimates are carried
out in Section 8.

\medskip

The upshot of the analysis is a diffeomorphism between the
approximate moduli space  and the
true moduli spaces
\best
{\cal AM}_s(Z_\la) \ \overset{\cong}\longrightarrow\ \M_s(Z_\la)
\eest
which intertwines with the attaching map of the domains and the
evaluation map into the target
(Theorem \ref{gluing thm}).  We then pass to homology, comparing and
keeping track of the
homology classes of the maps, the domains, and the constraints.  This
involves several
difficulties, all ultimately due to the fact that $H_*(Z_\la)$ is
different from both
$H_*(Z_0)$ and $H_*(X)\oplus H_*(Y)$.  This is sorted out in Section
10, where we define the
convolution operation and prove a first
Symplectic Sum Theorem:  formula (\ref{0.1}) holds when all stable
maps are  $\delta$-flat.

\smallskip

In Sections 11 and 12 we remove the flatness assumption by partitioning
the neck into a large number of segments   and using the
pigeon-hole principle as in Wieczorek \cite{w}. For that we construct spaces
$Z_\la^N(\mu_1, \dots,\mu_{2N+1})$, each symplectically isotopic to
$Z_\la$.   As  $(\mu_1, \dots,\mu_{2N+1})\to 0$
\best
&{\psfig{figure=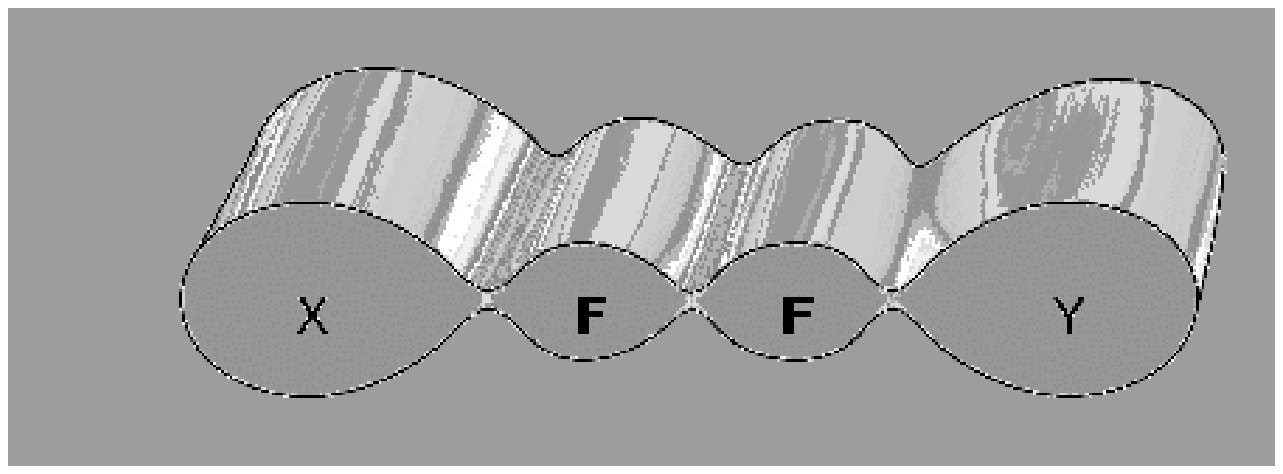,width=3.5in}}
\\
&\mbox{Figure  2:\ \  $Z_\la(\mu,\mu,\mu)$ for $|\mu|<<|\la|$}
\eest

\medskip

    \noindent   these degenerate to  the singular
space obtained by connecting
$X$ to $Y$ through a series of $2N$ copies of the  rational ruled
manifold ${\Bbb F}_V$ obtained
by adding an infinity section to the normal bundle
  to $V$.   An energy bound shows that for large $N$ each
map into $Z_\la(\mu_1,
\dots,\mu_{2N+1})$ must be flat in   most necks. Squeezing some or
all of the flat necks decomposes
the curves in $Z_\la$ into curves in $X$ joined to curves in $Y$ by a
chain of curves in
intermediate   spaces ${\Bbb F}_V$.   The limit maps are then
$\delta$-flat, so
formula (\ref{0.1}) applies
to each.  This process counts each stable map many times (there are
many choices of where to
squeeze) and in fact gives an  open cover of the moduli space.
Working through the
combinatorics and inverting a power series, we show that the total
contribution of the entire
neck  region between $X$ and
$Y$ is given by a certain $TW$ invariant of ${\Bbb F}_V$ --- the
$S$-matrix of Definition
(\ref{defSmatrix}).

\smallskip

The  $S$-matrix  keeps track of how the genus, homology class, and
intersection points with $V$
   change as the images of stable maps pass through the middle region
of Figure 2.
Observing  this back in the model of Figure 1, one  sees these
quantities changing abruptly
as the map passes through the neck --- the maps are ``scattered'' by
the neck.  The scattering
occurs when some of the stable maps contributing to the TW invariant
of $Z_\la$ have components
that lie entirely in $V$ in the limit as $\la\to 0$.   Those maps are
not $V$-regular, so are not
counted in the relative invariants of $X$ or $Y$.  But by moving to the spaces
of Figure 2 this complication can be analysized and related to the
relative invariants of the
ruled manifold   ${\Bbb F}_V$.

\medskip

The  $S$-matrix is the final subtlety.  With it in hand, we can at
last state our main result.

\bigskip\medskip

\noindent{\bf Symplectic Sum Theorem}{\em \quad   Let $Z$ be the symplectic
sum of $(X,V)$ and  $(Y,V)$ and suppose that
$\alpha\in{\Bbb T}(H_*(Z))$ splits as $(\al_X,\;\al_Y)$ as in
Definition \ref{def.al.separates}.
Then the TW invariant of $Z$ is given
in terms of the relative invariants of $X$ and $Y$  by
\bear
TW_Z(\al)\ =\   TW_X^{V}(\al_X) \, * \,
S_V\, * \, TW_Y^{V}(\al_Y)
\label{0.3}
\eear
where $*$ is the convolution operation (\ref{wteddiagonal1}) and
$S_V$ is the $S$-matrix (\ref{defSmatrix}).}

\bigskip

\noindent  A  detailed statement of this theorem is given in Section
12 and its extension to general constraints $\al$  is discussed in
Section 13.  We actually
state and prove (\ref{0.3})  as a formula for the {\em relative}
invariants of $Z$
in terms of the relative invariants of $X$ and $Y$  (Theorem
\ref{mainthm}).  In that form the
formula can be iterated.

\bigskip

Of course, (\ref{0.3}) is of limited use unless we can compute the
relative invariants of $X$
and $Y$ and the associated $S$-matrix.  That turns out to be
perfectly feasible, at least for
simple spaces.
   In Section 14 we build a collection of  two and four dimensional
spaces whose relative GW
invariants we can compute.  We also prove that the  $S$-matrix is the
identity in several cases
of particular interest.

\medskip

The last  section presents  applications.    The examples of section
14 are used as building
blocks to give short proofs of three recent results in
enumerative geometry: (a) the Caporaso-Harris formula for the number of
nodal curves in $\P^2$ \cite{ch}, (b)  the formula for the Hurwitz
numbers counting branched
covers of $\P^1$ (\cite{gjv} \cite{lzz}),  and (c) the ``quasimodular
form'' expression for  the rational enumerative
invariants of the rational elliptic surface  (\cite{bl}).  In
hindsight, our proofs of (a) and
(b) are essentially the same as those in the literature; using the symplectic
sum formula makes the proof considerably shorter and more
transparent, but the key ideas are the
same.   Our proof of (c), however, is completely different from that
of Bryan and Leung in
\cite{bl}. It is worth outlining here.

\smallskip

The rational elliptic surface $E$ fibers over $\P^1$ with a section
$s$ and fiber $f$.    For
each $d\geq 0$ consider the invariant
   $GW_{d}$ which  counts the number of connected rational stable maps in
the class $s+df$.  Bryan and Leung showed that the generating series
$F_0(t)=\sum GW_{d}\  t^{d}$ is
\bear
F_0(t)\ =\ \left(\ma\prod_{d}{1 \over 1-t^d}\right)^{12}.
\label{0.modular}
\eear
This formula is related to the work of Yau-Zaslow \cite{yz} and
is one of the simplest instances of some general conjectures
concerning counts of nodal curves in complex surfaces --- see
\cite{go}.

\smallskip

   While the intriguing form (\ref{0.modular}) appears in (\cite{bl}) for purely
combinatorial reasons,  it arises in our proof because of a
connection with  elliptic curves.  In
fact, our proof begins by relating $F_0$ to a similar series $H$
which  counts   {\em elliptic}
curves in $E$.  We then  regard $E$ as the fiber sum  $E\#(T^2\times
S^2)$ and apply the
symplectic sum formula.  The relevant relative invariant on the
$T^2\times S^2$ side is easily
seen to the generating function $G(t)$ for the number of degree $d$
coverings of the
torus $T^2$ by the torus.  The symplectic sum formula  reduces to  a
differential equation
relating $F_0(t)$ with $G(t)$,  and integration yields the quasimodular
form (\ref{0.modular}).  The details, given in section 15.3, are
rather formal;  the
needed geometric input is mostly contained in the symplectic sum  formula.

\medskip

All three of the applications in section 15 use the idea of `splitting the
target' mentioned at the beginning of this introduction.  Moreover,
all three follow  from rather
simple cases of the Symplectic Sum Theorem --- cases where the
$S$-matrix is the identity and where at least one of the relative
invariants in (\ref{0.3}) is
readily computed using elementary methods.  The full strength of the
symplectic sum theorem has
not yet  been used.

\bigskip\bigskip

This paper is a sequel to \cite{ip4}; together with \cite{ip4} it
gives a complete detailed exposition of the results announced in
\cite{ipa}.  Further applications have already appeared in \cite{ip2}
and \cite{i}.  Li and Ruan also have a sum formula  \cite{lr}.
Eliashberg, Givental, and Hofer are developing a general theory for
invariants of symplectic manifolds glued along contact boundaries
\cite{egh}.

\vskip.4in

\noindent{\bf\Large Contents}

\vskip.2in

\begin{minipage}{2.5in}
{\small
1.\ GW and TW Invariants \\
2.\ Symplectic Sums \\
3.\ Degenerations of symplectic sums  \\
4.\ The Space of Curves\\
5.\ Renormalization at the Nodes\\
6.\ The Space of Approximate Maps\\
7.\ Linearizations \\
8.\ The Eigenvalue Estimate
}
\end{minipage}
\hskip.1in
\begin{minipage}{3.4in}
{\small
9.\ The Gluing Diffeomorphism\\
10.\ Convolutions and the Sum Formula for Flat Maps\\
11.\ The space ${\Bbb F}_V$ and the S-matrix\\
12.\ The General Sum Formula\\
13.\ Constraints Passing Through the Neck\\
14.\ Relative GW Invariants in Simple Cases\\
15.\ Applications of the Sum Formula \\
     Appendix:\ Expansions of Relative TW Invariants
}
\end{minipage}

\vskip.1in


\vskip.4in


\setcounter{equation}{0}
\section{GW and TW Invariants}
\bigskip

For stable maps and their associated invariants we will use the
definitions and notation of \cite{ip4}; those are based on the
Gromov-Witten invariants as defined by Ruan-Tian \cite{rt1} and
Li-Tian \cite{lt}.  In summary, the definition goes as follows. A {\em
bubble domain} $B$ is a finite connected union of smooth oriented
2-manifolds $B_i$ joined at nodes together with $n$ marked points,
none of which are nodes. Collapsing the unstable components to points
gives a connected domain $st(B)$.  Let $\ov{\cal U}_{g,n}\to
\ov\M_{g,n}$ be the universal curve over the Deligne-Mumford space of
genus $g$ curves with $n$ marked points.  We can put a complex
structure $j$ on $B$ by specifying an orientation-preserving map
$\phi_0: st(B) \to \ov{\cal U}_{g,n}$ which is a diffeomorphism onto a
fiber of $\ov{\cal U}_{g,n}$.  We will often write $C$ for the
curve $(B,j)$ and use the notation $(f,C)$ or $(f,j)$ instead of
$(f,\phi)$.

A $(J,\nu)$-{\em holomorphic map} from $B$ is then a map $(f,\phi):B
\to X\times\ov{\cal U}_{g,n}$ where $\phi=\phi_0\circ st$ and which
satisfies $\bar{\partial}_Jf\ =\ \phi^*\nu$ on each component $B_i$ of
$B$.  A {\em stable map} is a $(J,\nu)$-holomorphic map for which the
energy
\bear
E(f,\phi)\ =\ \frac12 \int|d\phi|^2+|df|^2\
\label{2.defstable}
\eear
is positive on each component $B^i$. This means that  each component $B_i$
is either a stable curve or the restriction of $f$ to $B_i$ is non-trivial
in homology.

For generic $(J,\nu)$  the moduli space ${\M}_{g,n}(X,A)$ of stable
$(J,\nu)$-holomorphic maps
representing a class $A\in H_2(X)$ is a smooth orbifold of
(real) dimension
\bear
-2K_{X}[A]-\frac12(\dim X-6)\chi+2n
\label{dimM}
\eear
Its compactification carries a (virtual) fundamental class whose pushforward
     under the map
$$
\ov{\M}_{g,n}(X,A) \ \ \overset{{\rm {st}} \times {\rm{ev}}}\longrightarrow
\ \ \ov{\M}_{g,n}\ti X^n
$$
defined by stabilization and evaluation at the marked points is the
Gromov-Witten invariant $GW_{X,A,
g,n}\in H_*(\ov{\M}_{g,n}\ti X^n)$.  These can be  assembled   into a
single invariant  by setting
$\overline{\M} = \bigcup_{g,n}\overline{\M}_{g,n}$,    and introducing 
variables  $\lambda$ to keep track of the euler class and $t_A$ satisfying
$t_At_{B}=t_{A+B}$ to   keep
track of $A$.  The total GW invariant  of $(X,\w)$ is then
the formal series
\begin{equation}
GW_{X}\ =\ \sum_{A,g,n} \frac{1}{n!}\,GW_{X,A,g,n}\ t_A\ \la^{2g-2}.
\label{defn.7}
\end{equation}
whose coefficients lie in
$H^*(\ov\M)\otimes  {\Bbb T}(X)$ where ${\Bbb T}(X)$ denotes the
total tensor algebra ${\Bbb T}(H^*(X))$. This in turn defines the
``Taubes-Witten'' invariant
$$
TW_X\ =\ e^{GW_X}
$$
whose coefficients count holomorphic curves whose domains need not be
connected (as occur in \cite{t}).

\bigskip

The dimension (\ref{dimM}) is the index of the  linearization  the
$(J,\nu)$-holomorphic equation,
which is obtained as follows. A variation
of a map $f$ is specified by a $\xi\in \Gamma(f^*TX)$, thought of as a
vector field along the image, and a variation  in the complex structure
of  $C=(B,j,x_1,\dots,x_{n})$ is  specified by
\bear
k\in T_C\M_{g,n}\ \cong\ H^{0,1}_j\l(TB\otimes{\cal O}\l(-\sum{x_i}\r)\r)
\label{def.TM}
\eear
(tensoring with ${\cal O}(- x)$ accounts for the variation in the marked
point $x$).  Calculating the
variation in the path
\bear
(f_t, j_t)\ =\ \l(\exp_f (t\xi),\ j +tk\r)
\label{defbyexp}
\eear
one finds that the linearization at $(f,j)$ is the operator
\bear
D_{f,j}: \Gamma(f^*TX) \oplus T_C\M_{g,n} \to \Lambda^{0,1}(f^*TX)
\label{2.firstlinearization}
\eear
given by $D_{f,j}(\xi, k)\ =\ L(\xi) + Jf_*k$ with
\bear
     L(\xi)(w)\ = \ \frac12\left[ \nabla_w \xi+J\nabla_{jw}\xi +
\frac12(\nabla_\xi J)(f_*jw +Jf_*w-2J\nu(w))\right]- (\nabla_\xi \nu)(w)
\label{def.L}
\eear
where $w$ is a vector tangent to the domain and $\nabla$ is the  pullback
connection on $f^*TX$.  Writing $L$ as the sum of its $J$-linear component
${1\over 2}(L+JLJ)=\ov\partial_{f}+S$ and its $J$-antilinear component
$T$,
we have
\bear
L(\xi)(w) =  \ov\partial_{f,j}\xi
(w)+S(\xi,f_*w,f_*jw,w)+T(\xi,f_*w,f_*jw,w).
\label{2.L=d+T}
\eear
Here $\ov\partial=\sigma_J\circ \nabla$ with $\sigma_J$ the $J$-linear
part
of the symbol of $L$, $T$ is the tensor on $X\ti {\cal U}$ with
$J\,T(\xi, X,Y,w)$  given by
\best
\frac12\l[(\nabla_X J) +J(\nabla_{Y} J)\r]\xi
+{1\over 4}\l[ (\nabla_{J\xi} J)-J(\nabla_{\xi} J)\r](Y+JX-2J\nu(w))+
(\nabla_{J\xi}\nu- J\nabla_\xi\nu)(w)
\eest
and $S$ is a similarly looking tensor. Note that since the first two
terms of $L$ are complex linear we have, for  complex valued functions
$\phi$,
\bear
L(\phi\xi)=\ov\partial\phi\cdot \xi + \phi L (\xi)+(\ov \phi-\phi) T(\xi).
\label{L.alm.cx}
\eear

\bigskip

The invariant $GW_X$ was generalized in \cite{ip4} to an invariant of
$(X,\w)$ relative to a
codimension 2 symplectic submanifold $V$.  To  define it, we fix a pair
$(J,\nu)$  which is `$V$-compatible' in the sense of Definition 3.2 in
\cite{ip4}, that is, so that along $V$ the normal components  of
$\nu$ and of the tensor $T$ in
(\ref{2.L=d+T}) satisfy
\bear
\label{defcompatibleJnu}
&(a)& \mbox{$V$ is $J$-invariant and $\nu^N=0$, and}\\
&(b)&  \mbox{$T^N(\xi,X,JX-\nu,w)=0$  for
all $\xi\in N_V$, $X\in TV$ and $w\in TC$.} \nonumber
\eear
A stable map into $X$ is  called {\em $V$-regular} if no component of the
domain is mapped entirely into
$V$ and no marked point or node is mapped into $V$.  Any such  map has
only finitely many points $x_1, \dots, x_\ell$ in $f^{-1}(V)$.  After
numbering these, their degrees of contact with $V$ define a
multiplicity vector $s=(s_1, \dots , s_{\ell})$ and three associated integers:
\bear
\ell(s)=\ell, \qquad \deg s=\sum s_i,\qquad |s|=\prod s_i.
\label{def.ls}
\eear
The space of all $V$-regular maps is the union of components
$$
\M_{\chi,n,s}^V(X,A)\ \subset\ \M_{\chi,n+\ell}(X,A)
$$
labeled by vectors $s$ of length $\ell(s)$. This has a compactification
that comes with `evaluation'
maps
\bear
\ep_V:   \ov{\M}_{\chi,n,s}^V(X,A) \to \wt{\M}_{\chi,n}  \ti X^n \ti {\cal
H}_{X,A,s}^V.
\label{MHVSV}
\eear
Here $\wt{\M}_{\chi,n}$ is the  space of stable curves
with finitely many components, Euler class $\chi$ and $n$ marked points,
and ${\cal H}_{X,A,s}^V$ is   the
`intersection-homology' space  described in section 5 of \cite{ip4}.  There
is a covering map
$\ep:{\cal H}_{X,A,s}^V\to H_2(X)\ti  V_s$ whose first component records
the class $A$
and whose component in the space $V_s\cong V^{\ell(s)}$  records the image
of the last $\ell(s)$ marked
points. This covering is a necessary complication to the definition of
relative GW invariants.

\medskip

The complication occurs because of ``rim tori''.  A rim torus is an element of
\bear
{\cal R}\ =\ \mbox{ker }(\iota_*:H_2(X\setminus V)\to H_2(X))
\label{def.Rim}
\eear
where $\iota$ is the inclusion.  Each such element can be represented
as $\pi^{-1}(\gamma)$ where $\pi$ is the projection $S_V\to V$ from
the boundary of a tubular neighborhood of $V$ (the ``rim of $V$'') and
$\gamma:S^1\to V$ is a loop in $V$.  The group ${\cal R}$ is the group
of deck transformation of the covering
\begin{equation}\begin{array}{cccl}
{\cal R} & \longrightarrow &{\cal H}_{X}^V & \\
& & \Big\downarrow \ep&  \\
& &   H_2(X)\ti \ma \bigsqcup_s V_s. &
\end{array}
\label{HVcover}
\end{equation}
When there
are no rim tori (as is the case if $V$ is simply connected) ${\cal
H}_{X,s}^V$ reduces to $H_2(X)\ti V_s$ and the evaluation map (\ref{MHVSV})
is more easily described.

\medskip

The tangent space to $\M_{\chi,n,s}^V(X,A)$ is modeled on $\mbox{ker }D_s$
where $D_s$ is the
restriction  of (\ref{2.firstlinearization}) to the subspace where
$\xi^N$ has a zero of order
$s_i$ at the marked points $x_i$, $i=1,\dots,\ell$. It follows that
\bear
\mbox{dim }\M_{\chi,n,s}^V(X,A)\ =\ -2K_X[A]-\frac{\chi}2\,(\dim X-6) +2n-
2(\deg s-\ell(s))
\label{2.maindimformula}
\eear

\bigskip

With this understood, the definition of the relative GW invariant parallels
the above definition
     of $GW_X$: the image moduli space under (\ref{MHVSV}) carries a homology
class which,
after summing on $\chi, n$ and $s$, can be thought of as a map
\bear
GW_{X,A}^V:{\Bbb T}\l(H^*(X)\r) \ \longrightarrow\
     H_*(\ov \M\ti{\cal H}_X^V;\Q[\la]).
\label{relativeGW}
\eear
This gives the expansion
\bear
GW_X^V =\ \sum_{A,\; g}\;
     \sum_{s\mbox{ \tiny ordered seq}\atop \mbox{\tiny deg } s=A\cdot V}
{1\over \ell(s)!}\;
GW^V_{X,A, g,s}\ t_A\ \la^{2g-2}
\label{defnRelInvt2}
\eear
whose coefficients are (multi)-linear maps
${\Bbb T}\l(H^*(X)\r)\ra H_*\left(\ov \M\ti {\cal H}^V_{X,A,s}\right)$
(dividing by  $\ell(s)!$ eliminates
the redundancy associated with renumbering the last $\ell$ marked points).
The corresponding relative
Taubes-Witten invariant is again given by
\bear
TW_X^V= \mbox{exp}\,(GW^V_X).
\label{2.relTW}
\eear
After imposing constraints one can expand $TW_X^V$ in power series.  That
is done in the appendix
under the assumption that there are no rim tori.

\vskip.4in


\setcounter{equation}{0}
\section{Symplectic Sums}
\label{section1}
\bigskip

Assume $X$ and $Y$ are $2n$-dimensional symplectic manifolds each
containing symplectomorphic copies of a codimension two symplectic
submanifold $(V,\w_V)$. Then the normal bundles are oriented, and we
assume they have opposite Euler classes:
\bear
e(N_X V)+e(N_Y V)=0.
\label{1.1}
\eear
We can then fix a symplectic bundle isomorphism $\psi:N_X ^*V \to N_Y V$.

This data determines a family of symplectic  sums
$Z_\la=X\#_{V,\la} Y$ parameterized by $\la$ near 0 in $\cx$; these
have been described in \cite{gf} and \cite{mw}.  In fact, this family
fits together to form a smooth $2n+2$-dimensional symplectic manifold
$Z$ that fibers over a disk. In this section we will construct $Z$ and
describe its properties.
\medskip

\begin{theorem}
Given the above data, there exists a  $2n+2$-dimensional
symplectic manifold $(Z,\w)$ and a fibration
$\la:Z\to D$ over a disk $D\subset \cx$.  The center fiber $Z_0$ is the
singular symplectic
manifold $X\cup_VY$, while for $\la\neq 0$, the fibers $Z_\la$ are
smooth compact symplectic submanifolds --- the symplectic connect sums.
\label{thm1.1}
\end{theorem}

This displays the $Z_\la$ as deformations, in the symplectic category,
of the singular space $X\cup_VY$.  For $\la\neq 0$ these are
symplectically isotopic to one another and to the sums described in
\cite{gf} and \cite{mw}.

\bigskip

The proof of Theorem \ref{thm1.1} involves the following construction.
Given a complex line bundle $\pi: L\ra V$ over $V$, fix a hermitian
metric on $L$, set $\rho(x)=\frac12|x|^2$ for $v\in L$, and choose a
compatible connection on $L$.  The connection defines a real-valued
1-form $\alpha$ on $L\setminus\{\mbox{zero section}\}$ with
$\alpha(\partial/\partial\theta)=1$ (identify the principal bundle
with the unit circle bundle and pull back the connection form by the
radial projection). The curvature $F$ of $\alpha$ pulls back to
$\pi^*F=d\alpha$.  Then the 2-form
\begin{equation}
\w\ =\ \pi^*(\w_V)+\rho \pi^*(F)+ d\rho\wedge\alpha
\label{1.defw}
\end{equation}
is $S^1$-invariant, closed, and non-degenerate for small $\rho$.  The
moment map for the circle action $v\mapsto e^{i\theta}$ is the
function $-\rho$ because $i_{\frac{\partial}{\partial\theta}}\w\ =\
i_{\frac{\partial}{\partial\theta}}(d\rho\wedge\alpha ) \ =\ -d\rho$.

We can extend $\w$ to a compatible triple $(\w,J,g)$ as follows.  Fix
a metric $g^V$ and an almost complex structure $J_V$ on $V$ compatible
with $\w_V$ in the sense that
$$
g_V(X,Y)=\w_V(X,J_VY)
$$
for all tangent vectors $X$ and $Y$.  At each $x\in
L\setminus\{\mbox{zero section}\}$, there is a splitting $T_xL=V\oplus
H$ into a vertical subspace $V=\mbox{ker }\pi_*$ and a horizontal
subspace $H=\mbox{ker }d\rho \cap \mbox{ker }\alpha$.  We can
therefore identify $V=L_x$ and $H=T_{\pi(x)}V$ and define an almost
complex structure on the total space of $L$ by $J=J_L\oplus J_V$.
Writing $r(x)=|x|$ and $F_J(X,Y)=F(X,JY)$, one can then check that the
metric
\begin{equation}
g=\pi^*(g^V+\rho F_J) +(dr)^2 +r^2\,\alpha \otimes \alpha
\label{defg}
\end{equation}
is compatible with $J$ and $\w$.

The dual bundle $L^*$ has a dual metric $\rho^*(v^*)=\frac12|v^*|^2$
and connection $\alpha^*$ with curvature $-F$.  This gives a
symplectic form similar to (\ref{1.defw}) on $L^*$ and hence one on
$\pi:L\oplus L^*\to V$, namely
\begin{equation}
\w\ =\ \pi^*[\w_V+(\rho-\rho^*) F]+ d\rho\wedge\alpha - d\rho^*\wedge\alpha^*.
\label{1.defwZ}
\end{equation}
Below, we will denote points in $L\oplus L^*$ by triples $(v,x,y)$
where $v\in V$ and $(v,x, y)$ is a point in the fiber of $L\oplus L^*$
at $v$.  This space has
\bear
\mbox{(a)  a circle action $(x,y)\mapsto (e^{i\theta}x,e^{-i\theta}y)$ with
Hamiltonian $t(v,x,y)=\rho^*-\rho$}
\label{3.momentmap}
\eear

\ \hskip.23in (b) a natural $S^1$ invariant map $L\oplus L^*\to \cx$ by
$\la(z,x,y)=xy\in \cx.$

\medskip

\noindent Repeating the above construction of $J$ and $g$ gives an $S^1$
invariant compatible structure
$(\w,J,g)$ on $L\oplus L^*$.

\medskip

\noindent{\bf Proof of Theorem \ref{thm1.1}.}\ Let $L$ by the complex
line bundle with the same Euler class as $N_XV$ and give $L$ the above
structure $(\w,J,g)$.  Using $\psi$ and the Symplectic Neighborhood
Theorem, we symplectically identify a neighborhood of $V$ in $X$ with
the disk bundle of radius $\ep$ in $L$ and a neighborhood of $V$ in
$Y$ with the $\ep$-disk bundle in $L^*$.  Let $D$ denote the disk of
radius $\ep$ in $\cx$.

The space $Z$ is constructed from three open pieces: an $X$ end
$E_X=(X\setminus V) \times D$, a $Y$ end $E_Y=(Y\setminus V) \times
D$, and a ``neck'' modeled on the open set
\bear\label{u}
N=\{ \; (v,x,y)\in L\oplus L^*\; |\  |x|\le \ep,\; |y|\le \ep \}
\eear
These are glued together by the diffeomorphisms
\best
\psi_X:N\to (N_XV\setminus V) \times D  \qquad & \mbox{by     } &
(v,x,y)\mapsto (v,x,\la(x,y))
\\
\psi_Y:N\to (N_YV\setminus V) \times D \qquad & \mbox{by     } &
(v,x,y)\mapsto (v,y,\la(x,y))
\eest
This defines $Z$ as a smooth manifold.  The function $\la$ extends
over the ends as the coordinate on the $D$ factor, giving a projection
$\la:Z\to D$ whose fibers are smooth submanifolds $Z_\la$ for small
$\lambda\ne 0$.


\makebox(-100,130)[l]{Figure  3:\ \  Construction of $Z_\la$}
\centerline{\begin{picture}(110,150)(0,-10)
\put(-6,-6){\makebox(0,0){$V$}}
\put(0,0){\vector(1,0){100}}
\put(-9,100){\makebox(0,0){$L^*$}}
\put(100,-9){\makebox(0,0){$L$}}
\put(0,0){\vector(0,1){100}}
\put(0,0){\circle*{3}}
\put(3,40){\line(-1,0){6}}
\put(3,80){\line(-1,0){6}}
\put(-10,80){\makebox(0,0){$\ep$}}
\put(80,-10){\makebox(0,0){$\ep$}}
\put(60,-10){\makebox(0,0){$\la^{1/4}$}}
\put(40,3){\line(0,-1){6}}
\put(80,3){\line(0,-1){6}}
\put(55,70){\makebox(0,0){$Z_\la$}}
\put(50,65){\vector(-2,-1){25}}
\put(53,62){\vector(-1,-2){15}}
\put(15,50){\vector(0,1){45}}
\put(22,90){\makebox(0,0){$t$}}
\put(50,15){\line(1,0){40}}
\put(85,0){\dashbox(40,35)}
\put(80,0){\dashbox(50,35)}
\put(105,15){\makebox(0,0){$X\times D$}}
\put(75,0){\dashbox(60,35)}
\put(0,75){\dashbox(35,60)}
\put(0,80){\dashbox(35,50){$Y\times D$}}
\put(0,85){\dashbox(35,40)}
\put(50,50){\oval(70,70)[bl]}
\put(21,21){\circle*{3}}
\put(60,0){\oval(18,45)[t]}
\put(60,8){\makebox(0,0){$A$}}
\end{picture}  }
\bigskip 

In the region on the $X$ side near $|x|=\la^{1/4}$ (region $A$ in
Figure 3), we can merge the form $\psi_X^*\w$ into the symplectic
form (\ref{1.defwZ}) on $N$ by replacing $\alpha$ by
$\eta\alpha+(1-\eta)d\theta$ where $\eta(t)$ is a cutoff function with
$\eta=1$ for $|x|\leq\la^{1/4}$ and $\eta=0$ for $|x|\geq 2\la^{1/4}$.
The form (\ref{1.defwZ}) then extends over the $X$ end of $Z$. Doing
the same on the $Y$ side, we obtain a well-defined global symplectic
form $\w$ on $Z$. The restriction of $\w$ to a level set $Z_\la\cap
U_X$ is the original symplectic form $\w_X$ on $X$; similarly, its
restriction to $Z_\la\cap U_Y$ is $\w_Y$.  Finally, along $Z_\la\cap
U$ we have $\alpha^*=-\alpha$, so $\w$ restricts to
$$
\w_\la\ =\ \pi^*(\w_V-tF)- dt\wedge\alpha.
$$
This is non-degenerate for small $\la$.  Thus after possibly making $\ep$
smaller, we have a fibration $\la:Z\to D$ with  symplectic fibers.
\qed

\bigskip

This construction shows that the neck region $U$ of $Z$ has a  symplectic
$S^1$ action with Hamiltonian $t$.
This action preserves $\la$,  so  restricts to a Hamiltonian action on each
$Z_\la$.  In fact,  $t$ gives a
parameter along the neck, splitting each $Z_\la$ into manifolds with boundary
$$
Z_\la\ =\ Z_\la^- \ \cup\ Z_\la^+
$$
where $Z_\la^-$ is $Z_\la\cup U_X$ together with the part of $Z_\la\cup U$
with $t\leq 0$.  From this decomposition we can
recover  the symplectic manifolds $X$ and $Y$ in two ways:

\begin{enumerate}
\item  as $\la\to 0$,  $Z_\la^-$ (resp.  $Z_\la^+$) converges
to $X$ (resp. $Y$) as symplectic manifolds, or
\item  $X$ (resp. $Y$) is the symplectic cut  of $Z_\la^-$ (resp.
$Z_\la^+$) at $t=0$ (cf. \cite{l}).
\end{enumerate}
Thus we have collapsing maps
\begin{equation}
\begin{array}{ccc}
X\sqcup Y\ \ \  && \quad\ \  Z_\la
\\
{\scriptsize \mbox{$\pi_0$}}\searrow & & \swarrow {\scriptsize
\mbox{$\pi_\la$}}
\\
&Z_0 &
\end{array}
\label{2.collapsingmaps}
\end{equation}
and $\pi_\la$ is a deformation equivalence on the set where $t\neq 0$.

\medskip

The proof of Theorem \ref{thm1.1} constructs a structure $(\w, J, g)$ on
$Z$ whose restriction to $Z_\la$
on the $X$  end agrees with the given structure $(\w_X, J_X, g_X)$ on $X$
(and similarly on the $Y$ end).
More generally,  given $V$-compatible pairs $(J_X, \nu_X)$ and
$(J_Y, \nu_Y)$  which agree
along $V$ under the the map $\psi$ of (\ref{1.1}) and with the normal
components of $\nu_X$ and
$\nu_Y$ vanishing along $V$, then we can extend them to $(J,\nu)$ on
the entire fibration $Z$.

\medskip

We finish this section will a useful lemma comparing the canonical class of
the symplectic sum with the
canonical classes $K_X$ and $K_Y$ of $X$ and $Y$.

\begin{lemma} If $A\in H_2(Z_{\la};\Z)$, $\la\neq 0$,  is homologous in
$Z$ to the union $C_1\cup C_2\subset X\cup_V Y$ of cycles
$C_1$ in $X$ and $C_2$ in $Y$, then
$$
K_{Z_{\la}}[A]\ =\ K_Z[A]\ =\ K_X[C_1]+K_Y[C_2]+2\beta
$$
where $\beta$ is the intersection number $V \cdot [C_1] = V
\cdot[C_2]$. In particular, $K_{Z_{\la}}[R]=0$ for
any rimmed torus $R$ in (\ref{def.Rim}).
\label{addcanonicalclass}
\end{lemma}
\pf  For $\la\neq 0$, the normal bundle to $Z_{\la}$ has a
nowhere-vanishing section
$\partial/\partial\la$.  Thus the canonical bundle of $Z_{\la}$ is the
restriction of the canonical
bundle of $Z$, giving
$$
K_{Z_{\la}}[A]\ =\ K_Z[A]\ =\ K_Z[C_1]+K_Z[C_2].
$$
Outside the neck region of $X$, the tangent bundle to $Z$ decomposes
as $TX\oplus \cx$.  Inside
the neck region we have
$$
TZ\ =\ TX\oplus \pi^*\psi^*N_YV\ \cong\ TX\oplus \pi^*(N_XV)^{-1}
$$
where $\pi$ is the projection $N_XV\to V$. But the Poincar\'{e}
dual of $V$ in $X$, regarded as an element of $H^2_{cpt}(X)$, is the
chern class
$c_1(\pi^*N_XV)$. Since the
canonical class is minus the first chern class of the tangent bundle we
conclude that
$$
K_Z[C_1]\ =\  K_{X}[C_1] + V\cdot [C_1]
$$
and similarly on the $Y$ side.
\qed

\vskip.4in


\setcounter{equation}{0}
\section{Degenerations of symplectic sums}\label{S3}
\bigskip

The Gromov-Witten invariants of the symplectic sum $Z_\la$ are  defined in
terms of stable
pseudo-holomorphic maps from complex curves into the $Z_\la$.  The basic
idea of our connect sum formula
is to approximate the maps in $Z_\la$ by certain maps into the singular
space $Z_0$.  The first step is
to understand exactly which maps into $Z_0$ are limits of stable maps into
the $Z_\la$ as $\la\to 0$.  This section
gives a description of the limits of flat stable maps.  This `flat'
condition, defined below,
ensures that the limit has no components mapped into $V$.

\medskip

Fix a small $\de>0$. Given a  map $f$ into $Z_\la$, we can restrict
attention to that part of the image that lies in the `$\de$-neck'
\begin{equation}
     Z_\la(\de)\ =\ \{z=(v,x,y)\in Z_\la\ |\ ||x|^2-|y|^2|\leq \de\}.
\label{defndeltaneck}
\end{equation}
This is a narrow region symmetric about the middle of the neck in
Figure 3.  The energy of $f$
in this region is
\be
E_{\de}(f)\ =\ \frac12 \int|d\phi|^2+|df|^2
\label{2.defE(f)}
\end{equation}
where the integral is over $f^{-1} \l(Z_\la(\de)\r)$.

By  Lemma 1.5 of \cite{ip4} there is a constant $\alpha_V<1$,
depending only on
$(J_V,\nu_V)$ such that every component of every stable
$(J_V,\nu_V)$-holomorphic  map $f$ into
$V$ has energy
\bear
E(f)\ge \alpha_V.
\label{introalpha0}
\eear

\begin{defn}{\bf (Flat Maps)}
A  stable $(J,\nu)$-holomorphic  map $f$ into $ Z$ is {\em
flat} (more precisely $\de$-flat) if the energy in the
$\de$-neck is at most half $\alpha_V$, that is
\bear
      E_{\de}(f)\ \leq \ \alpha_V/2.
\label{defnFlateq}
\eear
\label{defnFlat}
\end{defn}
             For each small $\la$, let
$$
\Mf_{\chi,n}(Z_{\la},A)
$$
denote the set of flat maps in $\M_{\chi,n}(Z_{\la},A)$.  These are a
family of subsets of the space of stable maps and we write
\begin{equation}
\lim_{\la\to 0}\ \Mf_{\chi,n}(Z_{\la},A)
\label{2.flatintersection}
\end{equation}
for the set of limits of sequences of flat maps into $Z_\la$ as
$\la\to 0$. Because (\ref{defnFlateq}) is a closed condition this
limit set is a {\em closed} subspace of $\ov\M(Z)$. The remainder of
this section is devoted to a precise description of the space
(\ref{2.flatintersection}).

\medskip

\begin{lemma}
Each element of (\ref{2.flatintersection}) is a stable  map $f$ to
$Z_0=X\vcup Y$ with no irreducible components of the domain mapped entirely
into $V$.
\label{flatlemma1}
\end{lemma}
\pf  Each sequence  in (\ref{2.flatintersection}) has a subsequence $f_k$
converging in the
space of stable maps $\ov\M_{\chi,n}(Z,A)$ to a limit
$f:C\to  Z$.  In particular, the images converge pointwise, so lie in $Z_0$.

Suppose that the image of some component $C_i$ of $C$ lies in $V$.  Then
the restriction $f_i$ of $f$ to
that component satisfies $E(f_i)\leq E_\delta(f)$. Furthermore, by Theorem
1.6 of \cite{ip4} the sequence
$f_k$ (after precomposing with diffeomorphisms) converges in $C^0$ and in
$L^{1,2}$, so
$E_{\delta}(f)=\lim E_{\delta}(f_k) \leq  \alpha_V/2$.  This contradicts
(\ref{introalpha0}).
\qed

\bigskip

We can be very specific about how the images of the maps in
(\ref{2.flatintersection}) hit $V$.  By  Lemma
\ref{flatlemma1} and Lemma 3.4 of \cite{ip4},  at each point $p\in f^{-1}(V)$
the normal component of $f$ has a local expansion
$a_0z^d+\dots $.  This defines a local `degree of contact' with $V$
\begin{equation}
d\ =\ \mbox{deg}(f,p)\  \geq \ 1
\label{localdegree}
\end{equation}
and implies that $f^{-1}(V)$ is a finite set of
points.  Restricting $f$ to one component $C_i$
of $C$ and removing the points  $f^{-1}(V)$ gives a map from a connected
domain to the
disjoint union of $X\setminus V$ and $Y\setminus V$.  Thus the components
of $C$ are of two types:
those  components $C^X_i$ whose image lies in $X$, and those components
$C^Y_i$ whose image lies in $Y$. We  can therefore split   $f$ into two parts:
     the union of the components whose image lies in $X$ defines a map
$f_1:C_1\to X$,
from  a (possibly disconnected, prestable) curve  $C_1$, and the remaining
components define
a similar map $f_2:C_2\to Y$.

\begin{lemma}
     $f^{-1}(V)$ consists of nodes of $C$.  For each node $x=y\in f^{-1}(V)$
$$
     \mbox{\em deg}(f_1,x)\ =\  \mbox{\em deg}(f_2,y).
$$
\label{sumofdegrees}
\end{lemma}
\vspace{-.3in}
\pf The  local  degree (\ref{localdegree}) is a linking number.
Specifically, let $N_X(V)$ be a tubular
neighborhood of  $V$ in $X$ and let $\mu_X$ be the generator  of
$H_1(N_X(V)\setminus V)=\Z$ oriented as the boundary of a holomorphic disk
normal to $V$.  If $\mu_Y$ is the corresponding generator on the $Y$
side, then $\mu_X=-\mu_Y$ in $H_1$ of the neck $Z_\la(\delta)$. For each point
$x$ in $f_1^{-1}(V)$ and each small circle $S_\ep$ around $x$, the local
degree $d$ satisfies
$$
d\cdot \mu= [f_1(S_\ep)].
$$
If $x$ is not a node of $C$ then  by Theorem 1.6 of \cite{ip4}
$f_k$ converges to $f_1$ in
$C^1$ in a disk $D$ around $x$.  But then for large
$k$ $d\cdot \mu=[f(S_\ep)]=[f_k(S_\ep)]=[f_k(\partial D)]=0$, contradicting
(\ref{localdegree}).

Next consider  a node  $x=y$ of $C$  which is mapped into $V$.  Choose
holomorphic  disks $D_1=D(x,\ep)$ and  $D_2=D(y,\ep)$ that contain
no other points of $f^{-1}(V)$ and let $S_i=\partial D_i$. Then $S_1\cup
S_2$ bounds in $C$, so
$[f_k(S_1)]+[f_k(S_2)]=0$ in $H_1$ of the neck $Z_\la(\delta)$. Again,
$f_k\to f$ in $C^0$, so $0=[f(S_1)]+[f(S_2)]=d_1\mu_1+d_2\mu_2$ where
$\mu_i$ is either $\mu_X$ or
$\mu_Y$, depending on which side
$f(S_i)$ lies. Since $d_i>0$ the only
possibility is that $x=y$ is a node between a component in
$X$ and one in $Y$ and  $d_1=d_2$. \qed

\bigskip

Lemmas \ref{flatlemma1} and \ref{sumofdegrees} show that each map $f$ in
the limiting set
(\ref{2.flatintersection}) splits into $(J, \nu)$-holomorphic maps
$f_1:C_1\to X$ and $f_2:C_2\to Y$.
Numbering  the nodes in $f^{-1}(V)$ gives   extra marked points
$x_1,\dots,x_\ell$ on $C_1$ and matched $y_1,\dots,y_\ell$ on $C_2$ with
$s_i=\mbox{deg } x_i=\mbox{deg } y_i$.
Furthermore, the Euler characteristics $\chi_1$ of $C_1$ and  $\chi_2$ of
$C_2$ satisfy
\bear
\chi_1+\chi_2-2\ell=\chi.
\label{2.sumofchis}
\eear
\best
&{\psfig{figure=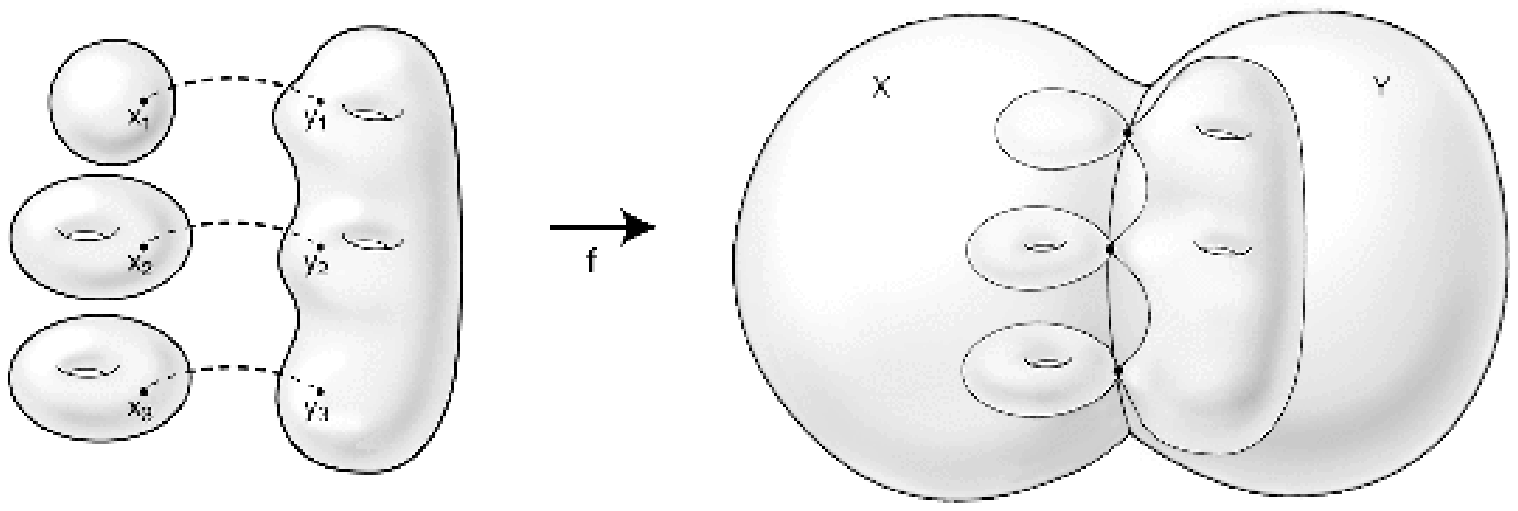,width=3.0in}}
\\
&\mbox{Figure  4:\ \ The map $f_0=(f_1,f_2)$ into $Z_0=X\cup_V Y$ }
\eest
\medskip

\begin{rem} To simplify the exposition we will assume for the rest of
the paper that
(a) all the components are stable, and
(b) $C_1$ and $C_2$ have no non-trivial automorphisms.
In fact, one can deal with unstable components by first stabilizing as
in \cite{lt}, and deal with automorphisms by lifting to a cover of the
universal curve as in \cite{rt2}, section 2. One can easily check
that the analytic arguments below, which are local on the moduli space
of holomorphic maps, carry through the stabilization and lifting
procedures.  Under these assumptions   the moduli
spaces are generically smooth and their virtual fundamental class is equal
to the actual fundamental class after taking  quotients as in \cite{rt2}.
\label{R.exp}
\end{rem}
\bigskip

We can now give a global description of how the limit maps $f$ in
(\ref{2.flatintersection}) are assembled from
their components $f_1$ and $f_2$.  First, consider how the domain curves
fit together in accordance with
(\ref{2.sumofchis}).  Given stable curves $C_1$ and
$C_2$ (not necessarily connected) with Euler
characteristics $\chi_i$ and $n_i+\ell$ marked points, we can construct a
new curve by identifying the last
$\ell$ marked points of
$C_1$ with the last $\ell$ marked points of $C_2$, and then
forgetting the marking of these new
nodes.  This  defines an attaching  map
\bear
\xi_\ell:\wt{\M}_{\chi_1,n_1+\ell}\ti \wt{\M}_{\chi_2,n_2+\ell}\;
\longrightarrow \; \wt{\M}_{\chi_1+\chi_2-2\ell,n_1+n_2}
\label{def.xi}
\eear
whose image is a subvariety of complex codimension $\ell$.  Taking the union
over all
$\chi_1, \chi_2, n_1$ and $n_2$ gives an attaching  map
$\xi_\ell:\wt{\M}\ti\wt{\M}\to\wt{\M}$ for each $\ell$.

Second, consider how the  maps fit together along $V$.  The evaluation map
$$
\ev_s: \M^V_{\chi,n,s}(X) \ti \M_{\chi,n,s}^V(Y)  \ \overset{\ep_V \ti
\ep_V}\longrightarrow \ {\cal H}_X^V
\ti {\cal H}_Y^V
\ \overset{\ep_2\ti \ep_2}\longrightarrow\  V_s\ti V_s.
$$
records the intersection points with $V$ and the pair $(f_1, f_2)$ lies in
the space
\bear
     \M^V(X) \ma\ti_{ev_s} \M^V(Y)\ \overset{\mbox{\small def}}=\
ev_s^{-1}(\Delta_s).
\label{2.inverseofdiagonal}
\eear
where $\Delta_s$ is the diagonal
$$
{\bold \Delta}_s\ \subset \  V_s\ti V_s.
$$
Denote by ${\cal H}_X^V\ti_{\ep} {\cal H}_Y^V=
(\ep_2\ti\ep_2)^{-1}({\bold \Delta})$ the fiber sum of ${\cal H}_X^V$ and
${\cal H}_Y^V$ along the evaluation map $\ep_2$, where
${\bold \Delta}=\ma\sqcup_s {\bold \Delta}_s$. Then we have a well
defined map
\bear
g:{\cal H}_X^V\ma\ti_{\ep} {\cal H}_Y^V \ra H_2(Z)
\label{def.g}
\eear
which describes how the homology-intersection data of $f_1$ and $f_2$
determine the homology class of $f$.

\begin{lemma} For generic $(J,\nu)$ the space (\ref{2.inverseofdiagonal}) is a
smooth orbifold of the same dimension as
$\M_{\chi,n}^{flat}(Z_\la,A)$ given by (\ref{dimM}).
\label{lemma2.4}
\end{lemma}
\pf The dimensions of $\M_{\chi_1,s}^V(X,A_1)$ and
$\M_{\chi_2,s}^V(Y,A_2)$ are given by (\ref{2.maindimformula}).  A
small modification of the proof of Lemma 8.6 of \cite{ip4} shows that
the evaluation map at the last $\ell=\ell(s)$ marked points (i.e. the
intersection points with $V$) is transversal to the diagonal
$\Delta\subset V^{\ell}\ti V^{\ell}$, imposing $\ell\,\dim V=\ell(\dim
X-2)$ conditions. Thus (\ref{2.inverseofdiagonal}) is a smooth
manifold of dimension
$$
-2K_X[A_1]-2K_Y[A_2]-4\deg s-\frac12(\dim X-6)(\chi_1+\chi_2-2\ell)+2n.
$$
The lemma follows by comparing this with (\ref{dimM}) using
(\ref{2.sumofchis}), Lemma \ref{addcanonicalclass}, and the fact that
$\deg s =A_1\cdot V=A_2\cdot V$.
\qed

\vskip.2in

Finally, note that renumbering the pairs $(x_i, y_i)$ of marked points
defines a free action of the symmetric group $S_\ell$ on
(\ref{2.inverseofdiagonal}) and the limit maps in
(\ref{2.flatintersection}) correspond to elements in the
quotient. Moreover, after ordering the double points along $V$
the limit set (\ref{2.flatintersection}) is
a closed subset
\bear
{\cal K}_\de\;\subset\M^V(X) \ma\ti_{ev} \M^V(Y)
\label{def.K}
\eear
which is the disjoint union of open sets labeled by $s$ and which has
compact closure as in \cite{ip4}. Since the maps in
$\M^{flat}(Z_\la)$  are $C^0$ close to  flat maps into $Z_0$ for small
$\la$   there is a decomposition
$$
\M^{flat}(Z_\la)\ = \ \bigsqcup_s \l(\l.\M_s^{flat}(Z_\la)
\r)\r/{S_{\ell(s)}}
$$
as a union of components labeled by ordered sequences $s=(s_1, s_2 
\dots)$.  As in the proof
of Lemma \ref{sumofdegrees}, these $s_i$ are
local winding numbers of the $\ell(s)$ vanishing cycles $S_\ep$.   In 
that form the labeling
extends to all continuous maps $C^0$ close to  flat maps into $Z_0$. 
Thus for small $\la$
$$
\M^{flat}_s(Z_\la)\subset \mbox{Map}_s(Z_\la)
$$
where $\mbox{Map}_s(Z_\la)$, the ``space of labeled maps'', is the 
set of labeled
continuous maps into $Z_\la$ which are $C^0$ close to  flat maps into $Z_0$.

Thus with this notation, the statements of Lemmas \ref{flatlemma1} and
\ref{sumofdegrees} translate into the commutative diagram
\bear
\begin{array}{ccc}
\ma\bigsqcup_s\  \M^V(X) \ma\ti_{ev_s} \M^V(Y)
& \longleftarrow &
\ma\lim_\la\ \l(\ma\bigsqcup_s \M_s^{flat}(Z_\la)\r)\\
\Big\downarrow &  & \Big\downarrow \\
(\wt{\M}\ti \wt{\M})\ti \l( {\cal H}_X\ma\ti_\ep {\cal H}_Y \r)&
\overset{\xi_{\ell(s)}\ti g}\longrightarrow &  \wt{\M}\ti H_2(Z).
\end{array}
\label{2.bigsquare}
\eear
The top arrow shows how  the maps that arise as limits of flat maps
decompose into pairs $(f_1,f_2)$ of $V$-regular maps into $X$ and $Y$,
while the bottom arrow keeps track of
the domains and homology classes
(the vertical maps arise from (\ref{MHVSV}) in the obvious way).

One then expects the top arrow in
(\ref{2.bigsquare}) to be a diffeomorphism for each $s$ and both sides
to be a model for the stable maps into $Z_\la$ for that $s$. The
analysis of the next six sections will show that this is true after
passing to a finite cover.

The necessity of passing to covers is dictated by  clustering 
phenomenon  mentioned in the
introduction: when
$s>1$ each curve in $Z_0$ is close  (in the stable map
topology)  to a cluster of curves in
$Z_\la$ for small $\la$, and these coalesce as $\la\to 0$.   To distinguish
the curves within a cluster and indeed to even verify this statement 
about clustering, it is
  necessary to use stronger norms and
distances --- strong enough that
the distances between the maps within a cluster are bounded away
from zero as $\la\to 0$.  The maps in a cluster can then be thought 
of as the fiber of a
covering of the space of limit maps.    The next three sections 
introduce the required  norms and
construct a first version of the covering.  The first step is to 
define an appropriate distance
function on the space of stable curves.

\vskip.4in


\setcounter{equation}{0}
\section{The Space of Curves}
\label{spaceofcurvessection}
\bigskip

Given two holomorphic maps, one can measure the
distance between their domain curves using some metric on the
Deligne-Mumford space $\ov\M_{g,n}$.  However, it is often more
convenient to fix a diffeomorphism of the domains, regarding the two
curves as two complex structures $j$ and $j'$ on a single 2-manifold,
and measuring the distance between $j$ and $j'$ using a Sobolev norm.
In this section we will define diffeomorphisms between nearby curves
in the universal family, fix  a Sobolev metric, and describe the
corresponding distance function on $\ov\M_{g,n}$.

Our distance function is designed so that a neighborhood of the image
of the attaching map
(\ref{def.xi}) is obtained by  gluing cylindrical ends of the
spaces $\M_{g,n}$.  It is  a complete metric on $\ov{\M}_{g,n}\setminus {\cal
N}$ where ${\cal N}$ is the set of all nodal curves; in particular it
is stronger than the
Weil-Petersson  metric.

\medskip

The construction starts by fixing a Riemannian metric $g_{\cal U}$ on
the universal curve $\ov\U_{g,n}\overset{\pi}\to\M_{g,n}$ compatible with
the complex
structure.  In the fibers of $\ov{\cal U}_{g,n}$ the `special points'
(marked points and nodes) are distinct and hence, by compactness, are
separated by a minimum distance.  After conformally changing the
metric we can assume that the separation distance is at least 4 and
that every fiber is flat in the disk of radius 3 around each of its
nodes.  We also fix a smooth function $\wt\rho$ on $\ov{\cal U}_{g,n}$
equal to the distance to the node in these disks of radius 3.
Finally, we replace $g_{\cal U}$ by a conformal metric that is
singular along the nodal points locus namely
\bear
g\ =\  \wt{\rho}^{-2} g_{\cal U}
\label{4.change1}
\eear

\medskip

To understand the geometry of this metric we focus attention to a
small ball $U$ in the set ${\cal N}_\ell$ of $\ell$ nodal curves and
construct a local model.  Each $C_0=C_0(u)=\pi^{-1}(u)$ with $u\in U$
is the union of not-necessarily-connected curves $C_1$ and $C_2$
intersecting at the nodes where points $x_k\in C_1$ are identified
with $y_k\in C_2$ for $k=1,\dots, \ell$.  For each $k$ we fix local
coordinates $\{z_k\}$ on $C_1$ and $\{w_k\}$ on $C_2$ centered at the
nodes. We can use the construction of Section 1 to form a family of
symplectic sums; for details see \cite{Masur}.  The result is a
holomorphic fibration
$$
{\cal F}\to U\times D^\ell
$$
     whose fibers  $C_\mu(u)$ are given
by $z_k w_k=\mu_k$ in disjoint balls $B_k$ centered on the nodes  and which
has a fixed a
trivialization outside a neighborhood of the nodes (here $D^\ell$ is the
unit disk in $\cx^\ell$).

\medskip

\begin{rem} The gluing parameters
$\{\mu_k\}$  are intrinsically  elements of the bundle
\bear
\ma\bigoplus_{k=1}^\ell\; ({\cal L}_{k}\otimes {\cal L}'_{k})^*
\label{model.N}
\eear
where ${\cal L}_k$ and ${\cal L}_k'$ the
relative cotangent bundle to $C_1$ at $x_k$ and respectively to $C_2$
at $y_k$.  Thus (\ref{model.N}) models the tubular neighborhood of ${\cal
N}_\ell$ in
$\ov\M_{g,n}$.
\end{rem}

\medskip

Fix a metric on ${\cal F}$ which is euclidean in the coordinates $(z_k,
w_k)$ on each $B_k$ and
$B_k$ has radius at least 4.  The induced metric on $C_\mu\cap B_k$ is
\bear g_\mu \ =\
\left. dz\,d\bar z+dw\,d\bar w\phantom{\int}\hspace{-.5em}\r|_{zw=\mu} \ =\
\l(1+{|\mu|^2\over r^4}\r) \l(dr^2+r^2\,d\theta^2\r)
\label{gmuvsg}
\eear
where $r=|z|$ and the distance to the node in $B_k$ is
$\rho^2=|z|^2+|w|^2=r^2+|\mu|^2/ r^2$.
       Switching to the conformal metric $ g\ =\ \rho^{-2}\,g_\mu$ as in
(\ref{4.change1}),  each
nodal curve in ${\cal F}$ has a cylindrical neck in each ball $B_k$.  In
fact, when $\mu_k\ne 0$
we can identify $C_\mu\cap B_k(1)$ with $[-T,T]\times S^1$ by writing
$r=\sqrt{|\mu_k|}\;e^t$ with $T=|\log\sqrt{|\mu_k|}|$.  In these
cylindrical coordinates $\rho^2=2|\mu|\cosh (2t)$ and
\bear
g\ =\ \rho^{-2}\,g_\mu\ =\ \l(r^{-1}dr\r)^2+d\theta^2 \ =\ dt^2+d\theta^2.
\label{conformalmetrics}
\eear
Similarly the curves $C_\mu$ have necks in $B_k$ which become longer and
longer as $\mu_k\to 0$.

\medskip

As in \cite{Masur} there is a biholomorphic map of fibrations from
${\cal F}$ to a neighborhood of $C_0$ in $\ov{\cal U}_{g,n}$ (if $C_0$
has automorphisms we lift to finite covers as in \cite{rt2}).  Because
this map is holomorphic with bounded differential its restriction to
each fiber is conformal and the conformal factor is bounded.
Consequently, the PDE results of the next several sections, all of
which involve only local considerations in the space $\ov{\cal
U}_{g,n}$, can be done in the model space ${\cal F}$ using the metric
(\ref{conformalmetrics}) and the results will apply uniformly on
$\ov{\cal U}_{g,n}$. We will henceforth consistently use this metric
(\ref{conformalmetrics}) on the domains of holomorphic curves.  Note
that the flatness condition (\ref{defnFlateq}) continues to hold
(after a uniform change of constants) because the energy density is
conformally invariant.

\medskip

We next describe how to lift vectors in $D^\ell$ to vectors on the family
${\cal F}=\{C_\mu\}\to
D^\ell$ around an $\ell$-nodal curve $C_0$. Fix smooth curves $C_\mu,
C_{\mu'}\in {\cal F}$. These
are identified outside a neighborhood of the nodes.  That identification
extends to a diffeomorphism $\phi=\phi_{\mu\mu'}: C_\mu\to C_{\mu'}$
as follows.  Using cylindrical coordinates
$z=\sqrt{|\mu|}e^{t+i\theta}$ around each node of $C_\mu$ and
$z'=\sqrt{|\mu'|}e^{t'+i\theta'}$ around the nodes of $C_\mu'$, set
$$
\phi(t,\theta)\ =\ \l(t+(2\al(t)-1)(T'-T),
\ \theta+\al(t)\arg({\mu\over \mu'})\r)
$$
where $T=\frac12|\log |\mu||$ and where $\al(t)$ is a cutoff function
equal to 1 for $t\geq 1$ and 0 for $t\leq -1$. Note that when $|z|\geq
1$, $\phi(t,\theta)= (t+T'-T,\, \theta)=(t',\,\theta')$, so
$\phi(z)=z'$.  The relation $zw=\mu$ similarly implies that
$\phi(w)=w'$ whenever $|w|\ge 1$.  Thus $\phi$ extends as claimed.

The corresponding infinitesimal diffeomorphism defines the lifts:  each
$v=(v_1,\dots,v_\ell)\in T_\mu D$ defines a family $\mu_s=\mu+sv$ and a
vector field
$$
\tilde{v}\ =\ \frac{d\,}{ds}\,\phi_{\mu\mu_s}|_{s=0}
\ =\ \sum_k  \l((\al-\frac12)\,\mbox{Re } \l(\frac{v_k}{\mu_k}\r),\
\al\cdot \mbox{Im } \l(\frac{v_k}{\mu_k}\r)\r)
$$
along $C_\mu$.  Going the other way, given any path $\mu_s$ in the
complement of the nodal set ${\cal N}$ we can lift the vectors
$\dot{\mu}$ as above and
integrate the lifted vector fields to get diffeomorphisms
$\phi_s:C_{\mu_0}\to C_{\mu_s}$.  For each $s$ the variation in the
complex structure is $h_s=\frac{d\,}{ds}\,\phi_s^*j$.  Define a second
distance between the complex structures by
\bear
\mbox{Dist }\l(C_{\mu_0},C_{\mu_1}\r)\ =\ \inf\ \int_0^1 \|h_s\|\ ds
\label{secondDistance}
\eear
where the infimum is over all paths from $\mu_0$ to $\mu_1$ in the
complement of ${\cal N}$ and
where
\bear
\|h\|^2\ =\ \int_{C_\mu}\ |\nabla^2 h|^2+ |\nabla h|^2+ |h|^2.
\label{wtedhnorm}
\eear

Note that in each family $\{C_\mu\}$ the nodal curves correspond to
$\{\mu\,|\mbox{some $\mu_k$ is zero}\}$.

\bigskip

\begin{lemma}
On the complement of the nodal set ${\cal N}=\{\mu\,|\mbox{some $\mu_k$
is zero}\}$ the Riemannian metric (\ref{wtedhnorm}) is uniformly
equivalent to the metric
\bear
     \sum_k\  \frac{1}{|\mu_k|^2}\,\mbox{Re } (d\mu_k)^2.
\label{easymetric}
\eear
\label{hdistancessamelemma}
\end{lemma}
\pf  Calculating $h= \frac{d\,}{ds}\,(d\phi_s^{-1}\cdot j\cdot d\phi_s)$ at
$s=0$, one finds that
$$
h\ =\ jd\tilde{v}-d\tilde{v} j\ =\ \begin{pmatrix}B & A\\ A & -B\end{pmatrix}
\qquad\mbox{where}\qquad
\l\{
\begin{array}{l}
A\ =\ \alpha' \, \mbox{Re } \l(\frac{v}{\mu}\r)\\
B\ =\ \alpha' \,\mbox{Im }\l( \frac{v}{\mu}\r).
\end{array}
\r.
$$
Noting that the  integrals of $|d\alpha|$, $|\nabla d\alpha|$,  and
$|\nabla^2 d\alpha|$ are independent of $\mu$ we then have
$$
\| h\|^2 \ = \
\sum_k  \frac{2|v_k|^2}{|\mu_k|^2} \int_{-1}^1\int_0^{2\pi}
|\nabla^2\alpha|^2+
|\nabla\alpha|^2+|\alpha|^2
\ = \ c\,\sum_k \frac{|v_k|^2}{|\mu_k|^2}. \qquad \Box
$$

\medskip

The metric (\ref{easymetric}) is cylindrical in each coordinate: using
polar coordinates
$\mu=r\,e^{i\theta}$ and
$r=e^t$ we have
$$
\frac{1}{|\mu|^2}\,\mbox{Re } (d\mu_k)^2
\ =\ \frac{1}{r^2}\l( dr^2 + r^2 d\theta^2\r)
\ =\ dt^2+d\theta^2.
$$
The corresponding distance function is that of the cylinder in each
coordinate, so for
$\mu=r\,e^{i\theta}=e^{t+i\theta}$ and $\mu'=r'\,e^{i\theta'}=e^{s+i\theta'}$
\bear
\mbox{dist}^2\ (\mu,\mu')\ =\
\sum_k\,|t_k-t'_k|^2+|\theta_k-\theta_k'|^2
\ =\ {\sum_k \l|\log\l(\frac{\mu_k'}{\mu_k}\r)\r|^2}.
\label{firsthdistance}
\eear
Thus the metric (\ref{easymetric}), defined in a neighborhood of the
nodal set ${\cal N}$, extends to a global metric on $\M_{g,n}=\ov
\M_{g,n}\setminus {\cal N}$ which is not complete ---near the stratum
${\cal N}_\ell$ of curves with $\ell$ nodes it is asymptotic to a
cylinder $W_\ell \ti \R_+^\ell$ where $W_\ell$ is a bundle over ${\cal
N}_\ell$ whose fiber is the real torus $T^\ell$ corresponding to the
bundle (\ref{model.N}). We can compactify this by identifying  the 
end $W_\ell \ti \R_+^\ell$
with  $W_\ell \ti(0,1)^\ell$ and  compactifying  to  $W_\ell 
\ti(0,1]^\ell$. This ``cylindrical end compactification'' projects
down to the Deligne-Mumford compactification $\ov \M_{g,n}$
so that the fiber along the nodal stratum ${\cal N}_\ell$ is a copy of
$W_\ell$.

\vskip.4in


\setcounter{equation}{0}
\section{Renormalization at the Nodes}
\label{S5}
\bigskip

In this section we will consider a sequence of flat
$(J,\nu)$-holomorphic maps
\bear
\label{4.sequence}
f_n:C_{\mu_n} \to Z_{\lambda_n}
\qquad\mbox{with}\qquad  \lambda_n\to 0.
\eear
By the Compactness Theorem for holomorphic maps $Z$ these
converge to a limit map $f_0$ from a nodal curve $C_0$ to $Z_0$ as
described in Section \ref{S3}; the convergence is in
$L^{1,2}\cap C^0$ and in $C^\infty$ on compact sets in the complement
of the nodes. We will refine this by constructing renormalized maps $\hat{f}_n$
around each node and
proving convergence results for the renormalized maps.  This gives
detailed information
about how the original  maps $f_n$ are converging in a
neighborhood of the nodes.

\bigskip

As in Section \ref{S3}, $C_0$ is the union of (not necessarily
connected) curves $C_1$ and $C_2$ which intersect at nodes, and $f_0$
decomposes into maps $f_1:C_1\to X$ and $f_2:C_2\to Y$.  Their images
meet along $V$ with contact vector $s=(s_1,\dots s_\ell)$; that is
there are points $x_k\in C_1$ and $y_k\in C_2$ so that $f_1$ and $f_2$
contact $V$ of order $s_k$ at $f_1(x_k)=f_2(y_k)\in V$.  For short, we
simply write
$$
f_n\ \to\ f_0=(f_1,f_2)\in {\cal K}\subset \M_s\ti_{ev}\M_s.
$$
where ${\cal K}$ is the compact set in (\ref{def.K}). Note that in
particular all
the estimates in the next sections will be uniform on ${\cal K}$.
\medskip

Around each node $x_k$ we can use the coordinates $(z_k, w_k)$ on the
domain described before (\ref{gmuvsg}), and coordinates $(v,x,y)$
centered on the image of the node to write $f_n=(v_n,f^x_n,f^y_n)$ (here
$v$ is a coordinate on $V$ and $(x,y)$ are coordinates in the sum
$N_X\oplus N_Y$ of the normal bundles to $V$).  These can be chosen so
$Z_\la$ is locally the graph of $xy=\la$.  Using the expansions of
$f_0$ provided by Lemma 3.4 of \cite{ip4} and Lemma
(\ref{sumofdegrees}) we can write around each node
\bear
f_0\ =\ (p_k+h^v,\ a_kz^{s_k}+h^x,\ b_kw^{s_k}+h^y)
\qquad\mbox{where}\qquad |h^v| \ \leq\ c\,\rho \mbox{ and }
|h^x|,|h^y|\le c\rho^{s_k+1}.
\label{f0expansion}
\eear

\medskip

\begin{rem}  The coefficient $a_k$  is the $s_k$-jet of the
function $f^x$ at $x_k$ modulo higher order terms, so
$$
a_k\in (T^*_{x_k}C)^{s_k}\otimes N_X
$$
where $N_X$ is the pullback of the normal bundle to $V$ in $X$.  The
evaluation map
$$
\ov\M^V_{g,n,s}\subset \ov\M_{g,n+\ell}\to \ov\M_{g,n}\times X^n\times
V^{\ell}
$$
determines complex line bundles ${\cal N}_X$ whose fiber at a map
$f$ is the normal bundle $N_X$ to $V$ at $f(x_k)$, and relative
cotangent bundles
$\L_k$ as in (\ref{model.N}) for
the last $\ell$ points of the domain. The leading coefficients are thus
sections
\bear
a_k\in\Gamma(\L_k^{s_k}\otimes {\cal N}_X)
\qquad{\mbox{and}}\qquad
b_k\in\Gamma((\L_k')^{s_k}\otimes {\cal N}_Y).
\label{akbkintrinsically}
\eear
\end{rem}
\bigskip

Note that the $f_n$ are nearly holomorphic with respect to the complex
structure $J_0$ defined by the coordinates $(v,x,y)$.  In fact, from
the $(J,\nu)$-holomorphic map equation, we have $\overline{\partial}
f_n=\overline{\partial}_{Jf_n}-(J-J_0)\,df_n \,j\ =\ \nu_n-
(J-J_0)\,df_n \,j$. Because the $f_n$ are converging in $C^0$ to the
continuous function $f_0$ with $|f_0|=0$ at $t=0$, we have the
pointwise bound
\bear
\label{growth0}
|\overline{\partial} f_n|\ \leq\ c |f_n|\,|df_n|\ \leq\ \ep\,|df_n|
\eear
where $\ep$ can be made arbitrarily  small by restricting the domain to a
small annular region in
the neck of $C_{\mu_n}$.  Similarly the metric on the target can be made
arbitrarily close to the
euclidean metric.  In the next lemma we consider such an  annular region in
cylindrical
coordinates
$(t,\theta)$ and estimate the energy
$$
E(f_n,T)\ =\ \frac12\int_{-T}^T\int_0^{2\pi}|df_n|^2\ dt\,d\theta
$$

\begin{lemma}
\label{energygrowthrate}
For small $r$ and large $n$, the energy  $E(t)=E(t,f)$ of $f=f_n$ on
the cylinder
$A(t)=[-t,t]\times S^1$ satisfies
\bear
E(t)\leq\ E(T)\  \rho^{\frac23}.
\label{energygrowthrate0}
\eear
Consequently,  there is a pointwise bound for $|df|$ of the form
\bear
|df|^2\ \leq\  c_1 \ E(t+1)\ \leq\ c_2 E\  \rho^{\frac23}.
\label{energygrowthrate1}
\eear
\end{lemma}
\pf By writing $f=u+iv$ one finds that
$$
4|\overline{\partial} f|^2\ dt\,d\theta\ =\  |df|^2\ dt\,d\theta -2\,d(u\,dv).
$$
Integrating over $A=A(t)$ and using Stokes' theorem  gives
$$
\label{growth1}
\frac12\int_{A} |df|^2\ =\ 2\int_{A} |\overline{\partial} f|^2\ +
\int_{\partial A}
u\ v_\theta\ d\theta.
$$
The boundary term is  an integral over two circles. On each, we can
replace $u$ by
$\tilde{u}=u-\frac{1}{2\pi}\int_0^{2\pi} u
\ d\theta$  and applying the H\"{o}lder and  Poincar\'{e}
inequalities on the circle
\bear
\int u\,v_\theta\ d\theta \ =\ \int \tilde{u}\,v_\theta\ d\theta\ \leq
\|\tilde{u}\|\,\|v_\theta \|\ \leq\
\|\tilde{u}_\theta\|\,\|v_\theta\|\ \leq\
\int|f_\theta|^2.
\label{holderpoincare}
\eear
Furthermore, from the definition $2\overline{\partial} f=f_t +if_\theta$
and the inequality
$(a-2b)^2\leq 2a^2+8b^2$ we obtain
\best
\label{growth3}
3|f_\theta|^2\ =\ 2 |f_\theta|^2 +|f_t-2\overline{\partial} f |^2\ \leq\  2
|df|^2
+8|\overline{\partial} f|^2.
\eest
Combining the previous three displayed equations and using
(\ref{growth0}) shows that
$$
\left(1-4c\ep^2\right)\,\int_A |df|^2\ \leq\
\frac43\left(1-4c\ep^2\right) \
\int_{\partial A} |df|^2.
$$
Taking $r$ small enough that $\ep=\ep(r)$ satisfies $4c\ep^2<1/44$, we obtain
$$
\frac23 E(t)\ \leq\ \ E'(t).
$$
Integrating this differential inequality from $t$ to $T$ yields
(\ref{energygrowthrate0}).

On the cylinder $[-T,T]\times S^1$, each point lies in a unit disk
with euclidean metric, and $f$ satisfies the equation
$\overline{\partial} f=\nu$.  Standard elliptic estimates then bound
$|df|$ at the center point in terms of the energy in that unit disk
(c.f. [PW] Theorem 2.3).  Thus (\ref{energygrowthrate0}) implies
(\ref{energygrowthrate1}).  \qed

\bigskip

In the next several sections we will repeatedly use the fact that in the
cylindrical metric
\bear
\int \rho^\alpha\ dt\ \sim\ c_\alpha\, \rho^\alpha \qquad \mbox{for}\
\alpha\neq 0.
\label{integralofrho}
\eear
Thus, for example, (\ref{growth0}) and (\ref{energygrowthrate1}) give
\bear
\|f^{-1}\overline{\partial}f\|_{p,A(r)}  \ \leq\
\|df\|_{p,A(r)}  \ \leq\ c_p\,\rho^{1/3}.
\label{Lpboundondf}
\eear
We will also use bump functions defined as follows. Fix a smooth
function $\beta:\R\to [0,1]$ supported on $[0,2]$ with $\beta\equiv 1$
on $[0,1]$.  The function $\beta_{\ep}(z,w)=\beta(\rho/\ep)$ has
support where $\rho^2=|z|^2+|w|^2\leq 4\ep^2$.  When restricted to
$C_\mu$, $\beta_{\ep}\equiv 1$ on the `neck' region $A_\mu(\ep)$
where $\rho\leq \ep$, and $d\beta_{\ep}$ is supported on two annular
regions where $\ep\leq \rho\leq 2\ep$.  We can choose $\beta$ so that using
the cylindrical metric (\ref{conformalmetrics})
\bear
|d\beta_\ep|\leq 2.
\label{bumpfunction}
\eear

\medskip

As before, we write $f_n=(v_n,x_n,y_n)$ in coordinates centered on the
image of each node.
\begin{defn}
In the region $\rho\leq 1$ around each node  define {\em renormalized maps}
$\hat{f}_n$ by
$$
\hat{f}_n\ =\ (\hat{v}_n, \hat{x}_n,\hat{y}_n)\ =\ (v^1_n-\bar{v}^1_n,\dots
v^k_n-\bar{v}^k_n,\frac{x_n}{az^s},\frac{y_n}{bw^s})
$$
where $\bar{v}^i_n$ is the average value of $v^i_n$ on the center circle
$\gamma_\mu=\{\rho=\sqrt{\mu}\}$  of $C_\mu$.
\end{defn}

Whenever  $\lambda_n=x_ny_n$ is non-zero $x_n$ has
no zeros and has (local) winding number $s$.  Hence each $\hat{x}_n$ has
winding number zero, so the functions $\log \hat{x}_n$, and similarly $\log
\hat{y}_n$, are 
well-defined.  The convergence (\ref{f0expansion}) shows that  on each set
$\rho\geq
r$ we have $\hat{x}_n\to f_0/az^s = 1+O(r)$ in $C^1$; hence there is a
constant $c$ so that
\bear
\sup_{r\le \rho\le 1}\ |\log \hat{x}_n| +|\log \hat{y}_n|\ \leq\ c\,r
\hskip.4in \forall n\geq N=N(\rho).
\label{c1boundonG}
\eear

\medskip

\begin{lemma} For each sequence (\ref{4.sequence}) we have
$\displaystyle \lim_{n\to\infty}  \frac{\lambda_n}{\mu_n^s}  \ =\  ab.$
\label{lambdazlemma}
\end{lemma}
\pf For $G_n=\log \hat{x}_n$ the integral
$$
\bar{G}_n(\rho)\ =\ \frac{1}{2\pi}\,\int_0^{2\pi} G_n\ d\theta
$$
over the circles with fixed $\rho$ satisfies
$$
\frac{d\,}{dt}\,\bar{G}_n\ =\ \int \partial_t G_n\ =\
\int (\partial_t +i\partial_\theta)G_n\ =\
2\int x_n^{-1}\,\ov\partial x_n\ d\theta.
$$
To bound $x_n^{-1}\,\ov\partial x_n$ we regard $f$ locally as a map into
$V\times \cx$ with its
product almost complex structure $J_V$ and a product metric $g_V$.  Then
$\nu^N$ and $J-J_V$
are both $O(R)$ where $R^2=|x|^2+|y|^2$, so a slight modification of
(\ref{growth0}) gives
$|\ov\partial x_n| \leq  cR\,|df|_{g_V}$.  Also noting that
$g=\l(1+|\la|^2/|x|^4\r)g_V$ as in (\ref{gmuvsg}), we obtain
\bear
|x_n^{-1}\,\ov\partial x_n|^2\ \leq \
c\left( \frac{|x|^2+|y|^2}{|x|^2}\right)\, |df|^2_{g_V}
\ =\ c \left(1+\frac{\la^2}{|x|^4}\right)\, |df|_{g_V}^2
\ =\ c\,|df|^2.
\label{boundepartialx}
\eear
These equations and Lemma \ref{energygrowthrate}  give
$|\frac{d\, }{dt}{\bar{G}}_n|\ \leq\ c_1\,\rho^{1/3}$.  Hence for
$\rho\leq r$ and $n>N(r)$
\bear
|\bar{G}_n(\rho)|\ \leq\ |\bar{G}_n(r)|\ +\ c_1\int^r_\rho
\rho^{1/3}\ dt\ \leq\
|\bar{G}_n(r)|\ +\  c_2\,r^{1/3}
\ \leq\ c_3\,r^{1/3}
\label{avgGbound}
\eear
(the last inequality uses (\ref{c1boundonG})).  This implies that
$$
\left|\frac{1}{2\pi}\,\int_{\ga_{\mu_n}} \log \hat{x}_n\ \right|\ =\
\l|\bar{G}_n\l(\sqrt{|\mu_n|}\r)\r|\ \to\ 0
$$
as $n\to \infty$.  The same limit statement holds with
$G$ replaced by $\log \hat{y}$. For each $n$ we can then integrate the constant
$$
\log\left(\frac{\lambda_n}{ab\mu_n^s}\right)\ =\
\log\left(\frac{x_n y_n}{az^s\,bw^s}\right)\ =\
\log\left(\hat{x}_n\,\hat{y}_n\right)
$$
over $\gamma_{\mu_n}$  to see that
$$
2\pi \,\log\left(\frac{\lambda_n}{ab\mu_n^s}\right)\ =\
\int_{\ga_{\mu_n}} \log\left(\frac{\lambda_n}{ab\mu_n^s}\right)\ =\
\int_{\ga_{\mu_n}} \log \hat{x}_n\ + \log \hat{y}_n\ \ \to \ 0.
$$
The lemma follows.
\qed

\vskip.3in

\begin{lemma}
Set $\xi^V_n=v_n-\ov{v}_n$ and $\xi^N_n=(x_n, y_n)$.  Then for each $p\geq
2$ there are
constants $C$ and $N=N(r)$ so that whenever $\delta\leq \frac13$, $r\le 1$
and $n\geq N$
\bear
\int_{\rho<r} \rho^{-p\delta/2}\,\l( |\nabla\xi_n^V|^p + |\xi_n^V|^p +
\rho^{(1-s)p} |\nabla\xi_n^N|^p
+\rho^{(1-s)p} |\xi_n^N|^p \r) \ \leq\ C_p\,r^{p/6}.
\label{firstnormxi}
\eear
\label{1pconvergence}
\end{lemma}
\pf  Write the annular regions  $A_\mu(r)=\{ \rho<r \}$ in the neck as the
union of annuli $A_k=\{k\leq
t\leq k+1\}$ of unit size and let $\rho_k$ be the value of $\rho$ at one
end of $A_k$.     Since
$|d\xi^V_n|=|dv_n|\leq c\rho^{1/3}$ by (\ref{energygrowthrate1}) we  have
\bear
\sup_{A(r)}\ |\xi^V_n|
\ \leq\  \sum_k |\mbox{osc}_{A_k}\ \hat{v}|
\leq\ C\,\sum_k \ \|d v\|_{4,A_k}
\leq\ C\,\sum_k \ \rho_k^{1/3}
\leq\ C\, r^{1/3}
\label{sup<dv}
\eear
where the last inequality comes from the  Riemann sum for $\int
\rho^{1/3} dt$.  Thus $|\nabla\xi_n^V|^p + |\xi_n^V|^p\leq c\rho^{p/3}$
pointwise.  Integrating via
(\ref{integralofrho}) then gives  the first half of (\ref{firstnormxi}).

Next,  the Calderon-Zygmund inequality of \cite{is} shows that $G\ =\ \log
\hat{x}_n$ satisfies
$$
\|dG\|_{p,A(r)}\ \leq\
C\,\|\overline{\partial}(\beta_r G)\|_{p,A(2r)}\ \leq\
C\,\left(\|d\beta_r\cdot G\|_{p,A(2r)\setminus A(r)}\ +\
\|x_n^{-1}\overline{\partial}x_n\|_{p,A(2r)}\r)
$$
We can integrate (\ref{boundepartialx}) as in (\ref{Lpboundondf}), and use
(\ref{bumpfunction})
and  the bound (\ref{c1boundonG})  in the region $r\leq\rho\leq 2r$ where
$d\beta_r\neq 0$. These imply that the $L^p$ norm of
$dG$ is bounded by $c\,r^{1/3}$.  But then for each annulus $A\subset A(r)$
with unit
diameter we can use (\ref{avgGbound}) and a Sobolev inequality to obtain
$$
\sup_A\ |G|\ \leq\ |\mbox{avg}_{\bd A}\ G|\ +\ \mbox{osc}_A\ G
\leq\ c\,r^{1/3} \ +\ C\,\|dG\|_{4,A(r)}
\ \leq\ c\,r^{1/3} \qquad \mbox{for all } n\geq N(r).
$$
Exponentiating this bound on $G$ shows that $|\hat{x}_n-1|\leq
cr^{1/3}$ in $A(r)$, and that in turn gives $|d\hat{x}_n| =
|\hat{x}_n\,dG|\leq c\,|dG|$. Consequently $\xi_n^x=  \hat{x}_n(az^s)$
satisfies
\bear
|\xi^x_n|\leq c\rho^{s+1/3}
\qquad\mbox{and}\qquad
|d\xi^x_n|\leq c\rho^s\l(1+|dG|\r)
\label{pointwiseboundsonxi1}
\eear
(after noting that $|dz/z|$ is bounded).  Of course the same bounds hold
for the $y$ components, so
integration, combined with the $L^p$ bound
on $dG$,  gives the second half of (\ref{firstnormxi}).
\qed

\vskip.4in


\setcounter{equation}{0}
\section{The Space of Approximate Maps}
\bigskip

The limit argument of section 3 shows that as $\la\to 0$ holomorphic
maps $f_\la$ into $Z_\la$ converge to maps into $X\cup Y$ with
matching conditions along $V$, i.e. to a maps in
$\M_s^V(X)\ti_{ev}\M_s^V(Y)$.  Lemma \ref{1pconvergence} gives further
information about the convergence near the matching points; it shows
that for small $\la$ the maps $f_\la$ are closely approximated by maps
$g(z,w)= (\ov{v}, az^s, bw^s)$ in local coordinates.  Over the next
four sections we will reverse this process, showing how one can use
$\M_s\ti_{ev}\M_s$ to construct a model ${\cal AM}_s(\la)$ for the space of
stable maps into $Z_\la$.  The final result is stated as Theorem
\ref{gluing thm}.

The construction has two main steps.  In the first, maps $f$ in a
compact set ${\cal K}\subset \M_s^V(X)\ti_{ev}\M_s^V(Y)$ are smoothed in a
canonical way to construct maps $F$ into $Z_\la$ which are nearly
holomorphic.  The second step corrects those approximate maps $F$ to
make them truly holomorphic.  This section describes the canonical
smoothing and the resulting space of approximate maps and introduces
norms on the space of maps which capture the convergence of the
renormalized maps.  Those norms lead to a precise statement that the
approximate maps are nearly $(J,\nu)$-holomorphic.

\medskip

The maps alone cannot be canonically smoothed --- more data are
needed.  This harks back to the comment at the end of Section 3 that
each $f$ will generally be the limit of many maps into $Z_\la$.
Recall that $f\in \M_s^V(X)\ti_{ev}\M_s^V(Y)$ is a map from an
$\ell(s)$-nodal curve $C_0$ whose nodes $x_k=y_k$ are mapped into $V$
with contact of order $s_k$.  As in section 3 $C_0$ has an $\ell$
dimensional family of smoothings $C_\mu$, $\mu=(\mu_1,\dots,
\mu_\ell)$.  Lemma \ref{lambdazlemma} shows that $C_0$ is the limit of
maps into $Z_\la$ only if $\mu$ satisfies $a_kb_k\,\mu_k^{s_k}=\la$.
That leaves $|s|= s_1 s_2\cdots s_\ell$ possibilities for $\mu$
corresponding to the different choices of root for each $\mu_k$. Thus
the maps into $Z_\la$ near $f$ are specified by pairs $(f,\mu)$, with
the $\mu$ specifying the deformation of the domain.

\medskip

Globally, we have $\la\in N_X\otimes N_Y\cong \cx$ (via a fixed
trivialization), so
(\ref{akbkintrinsically}) implies that at each node  the coefficients
$a_k,b_k$  determine a section
$$
\frac{\la}{a_k b_k}\in\Gamma\left(\L_k^{*}\otimes (\L_k')^*\right)
$$
over $ \M_s^V(X)\ti_{ev}\M_s^V(Y)$.  The ${s_k}^{th}$ root of this section is a
multisection of $\L_k^{*}\otimes (\L'_k)^*$; considering all $k$ at
once defines a multisection of the direct sum of the $\L_k^{*}\otimes
(\L_k')^*$. This gives an intrinsic model for our space ${\cal AM}_s(\la)$ of
approximate
maps:
\begin{defn}
\label{defofAsla}
For each $s$ and $\la\neq 0$, the {\em model space} ${\cal AM}_s(\la)$ is
the multisection of
\best
\ma \bigoplus_{k=1}^{\ell}\; [\L_k^{*}\otimes (\L_k')^*] \to
\M_s^V(X)\ti_{ev}\M_s^V(Y)
\eest
consisting at $f_0$ of those $\mu=(\mu_1,\dots,\mu_\ell)$ which satisfy
\bear
\mu_k^{s_k}\ =\ \frac{\la}{a_kb_k}
\qquad{\mbox{for each $k$.}}
\label{defofmuk}
\eear
\end{defn}

This model space is an $|s|$-fold cover of
$\M_s^V(X)\ti_{ev}\M_s^V(Y)$, and hence is a manifold for generic $(J,\nu)$.
Elements of the model space are pairs $(f,\mu)$ where $f:C_0\to Z_0$
and $\mu$ satisfies (\ref{defofmuk}).  Each such
element  gives rise to an approximate holomorphic
map as follows.

\medskip

\begin{defn}
For each $(f,\mu)\in {\cal AM}_s(\la)$, $\la\neq 0$,  define an {\em
approximate holomorphic
map} $F=F_{f,\mu}:C_{\mu}\to Z_{\la}$ by
\bear
F\ =\ f- \sum \,\beta_k (f-p_k)
\label{defapproxmaps}
\eear
where $\beta_k$ is bump function (\ref{bumpfunction}) with
$\ep=|\la|^{1/4{s_k}}$ in coordinates $(z_k,w_k)$ around the
$k^{\mbox{\scriptsize th}}$ node, $p_k$ is the image of the node in those
coordinates, and $f$ is the restriction of $f(z,w)=(v(z,w),x(z),y(w))$
to $C_{\mu}$.
\label{defApproxMaps}
\end{defn}

Altogether, the association $(f,\mu)\mapsto (F_{f,\mu}, C_{\mu})$ defines
a `gluing map'
\bear
\Gamma_\la:  {\cal AM}_s(\la) \ \to\  \mbox{Maps}_s(C,Z_\la\ti {\cal U})
\label{firstGamma}
\eear
This map  is injective:  if $\Gamma_\la(f,\mu)= \Gamma_\la(f',\mu')$ then $f$
and $f'$  are  $(J,\nu)$-holomorphic maps which agree on the set where
$\rho>1$  and therefore, by the
unique continuation property of elliptic equations, agree everywhere.

In section \ref{GluingDiffeomorphismsection} we will show that
$\Gamma_\la$ is a diffeomorphism onto a submanifold.  Here, as a
preliminary, we introduce norms which make the space of maps in
(\ref{firstGamma}) into a Banach manifold.

\bigskip\bigskip

\non{\bf Norms}. We will use weighted Sobolev norms tailored for our
problem.  On the domain we continue to use the cylindrical metric
(\ref{conformalmetrics}) and to use (\ref{secondDistance}) to measure
distance between curves.  In the target we identify a neighborhood of
$V$ in $Z$ with the disk bundle of the bundle $N_X\oplus N_Y$ over
$V$; in that neighborhood we can then decompose vector fields $\xi$
into components $(\xi^V,\xi^x,\xi^y)$ where $\xi^V$ is horizontal with
respect to the connection on $N_X\oplus N_Y$ and the $\xi^x$ and
$\xi^x$ are tangent to the fibers of $N_X$ and $N_Y$.

In a neighborhood of each node $p_k$ let
$\gamma_k$ be the circle $\rho=\sqrt{|\mu|}$ in $C_\mu$ and let
$\beta_k$ be as in
(\ref{defapproxmaps}).  Given $\xi$,
subtract the
(extended) average value of the $V$ components, defining
\bear
\zeta\ =\ \xi^V-\b_k \ov{\xi}
\qquad \mbox{where}\qquad
\ov \xi \ =\ \frac{1}{2\pi}\int_{\gamma_k}\xi^V\in T_{p_k}V.
\eear
These  averaged vectors $\ov\xi_k$ at the different nodes can be assembled
into a single vector
$\ov\xi\in T_pV^\ell$ where $p=(p_1,\dots, p_\ell)$. Similarly, the
$(\zeta_k, \xi_k^x,\xi_k^y)$ extend to a global vector field on $C_\mu$
\bear
\zeta=\xi-\sum \beta_k\bar{\xi}_k
\label{def.zeta}
\eear
where $\beta_k$ is the bump function (\ref{bumpfunction}) with $\ep=1$
centered on the node
$p_k$. Fix $\delta>0$ and set
\bear
\|\zeta\|_{1,p,s}^p \ =\ \int_{C_\mu} \rho^{-\delta p/2}\, |\nabla (W\zeta)|^p
+ |W\zeta|^p
\label{norms}
\eear
where $\nabla$ is the covariant derivative of the cylindrical metric on
the domain and the metric induced on $Z_\la$ from $Z$ while the
endomorphism $W:\zeta\to W\zeta$ weights the normal components around
each node:
\bear
W\zeta\ =\ (1-\sum \b_k)\zeta\ +\ \sum
\beta_k\l(\zeta^V,\frac{\zeta_k^x}{{z}^{s_k-1}},\,
\frac{\zeta_k^y}{{w}^{s_k-1}}
\r).
\label{def.W}
\eear
Here a complex valued function  $\phi=u+iv$  acts on (real) vector
field $\zeta$ by  $\phi\cdot \zeta=u\cdot \zeta+v\cdot J\zeta$, with
     $|\phi\zeta|=|\phi|\cdot |\zeta|$.

\begin{defn}  Given a tangent vector $(\xi, h)$ to the space
$\mbox{Map}_s(C,Z_\la\ti {\cal U})$, we
form the triple $(\zeta,\bar{\xi},h)$ as in (\ref{def.zeta}) and define the
weighted
$L^1_s$ norm
\bear
\|(\xi, h)\|_1\ =\  \|\zeta\|_{1,2,s} +\|\zeta\|_{1,4,s} \ +\  |\bar{\xi}|
\ + \|h\|
\label{defofnorms}
\eear
where $\|h\|$ is given by (\ref{wtedhnorm}). For  1-forms
$\eta\in \Omega^{0,1}(f_\mu^* TZ_\la)$ we do the same without averaging:
\bear
\|\eta\|_1\ =\  \|\eta\|_{1,2,s} +\|\eta\|_{1,4,s}.
\label{defofnorms2}
\eear
The weighted $L^0_s$ norm $\|\cdot \|_0$ and the weighted
$L^2_s$ norm $\|\cdot \|_2$ are defined similarly.
\label{def.norms}
\end{defn}

\medskip

The norm $\|(\xi, h)\|_1$ dominates the $C^0$ norm  (since
$L^{1,4}\hookrightarrow C^0$). Hence we can use it to complete the
space of $C^\infty$ maps, making $\mbox{Map}_s(C,Z_\la\ti {\cal U})$
a Banach manifold with   neighborhoods modeled by
$(\xi,h)$.

\medskip

\begin{rem} Note that the norms defined above make sense also at $\la=0$, where
the average value $\ov\xi$ is equal to the value of $\xi$ at the
double points $x_k, y_k$.
\end{rem}

\bigskip

We conclude this section by showing that the approximate maps are nearly
holomorphic.  The specific
statement is that the quantity
$\overline{\partial}F -\nu_{F}$, which measures the failure of the
approximate map to be $(J,\nu)$-holomorphic, is small in the norms just
introduced.

\bigskip

\begin{lemma}
\label{Boundondg-nulemma}
For $\delta\leq \frac13$ and $\la$ sufficiently small, each $F=F_{f,\mu}$
satisfies
$\|\overline{\partial}F -\nu\|_0 \ \leq\  c \,|\la|^{1/6|s|}.$
\end{lemma}
\pf Let $N_k(\mu)$ be the region around the node $p_k$ where $\rho\leq
2\sqrt[4]{|\mu|}$.  Outside
$\cup_k\,N_k(\mu)$  $F\equiv f$ is $(J,\nu)$-holomorphic and therefore
$\Phi=\overline{\partial}F
-\nu_{F}$ vanishes. When $\la\sim\mu^s$ is sufficiently small the image of
each $N_k(\mu)$ lies in a neighborhood of $V$ where we can separate
components tangent and normal
to $V$. Taking $F=f-\beta(f-p)$,
\best
\Phi(z)\ =\ (1-\beta)\overline\partial_{J_F}f +d\beta\l[ (f-p)-J_F(f-p)j\r]
- \nu(z,F(z))
\eest
where $J_F$ means $J$ at the point $F(z)$. Since
$f$ is $(J_f,\nu_f)$-holomorphic and $|d\beta|$ is bounded
\best
|W\Phi| &\leq &  |W((J_F-J_f)df)| +c|W(f-p)| +|W(\nu_F-\nu_f)|.
\eest
Now the local expansions of $x(z)$ and $y(z)$ show that $|W(f-p)|\approx
|(v-v_0,az, bw)|\leq c\rho$ and similarly
$|W((J_F-J_f)df)|\leq c\rho$. Because
$\nu^N$ vanishes along $V$, the normal component of
$\nu_F-\nu_f$ is bounded by $c|F^N|\le c\rho^{s}$, while
$|(\nu_F-\nu_f)^V|\leq c|F-f|\le
c\rho$. Thus $W\Phi\leq c\rho$.   Integrating  over $N_k(\mu)$ using
(\ref{integralofrho}) then gives
\best
\|W\Phi\|_{0,p,s; N_k(\mu)}^p
\ \leq\ c \,|\mu_k|^{\frac{p}{6}} \ \leq\  c\,|\la|^{p/(6s_k)}.
\eest
The lemma follows by summing on $k$ and on $p=2,4$.
\qed

\vskip.4in


\setcounter{equation}{0}
\section{Linearizations}
\bigskip

This section  describes the linearization of the 
$(J,\nu)$-holomorphic map equation as an
operator on the Sobolev spaces of Definition \ref{def.norms}.  We do 
this first for the space
$\M_s^V(X)$ which defines the relative invariants, then for the space 
$\M_s(Z_0)=\M_s^V(X)\ti_{ev}\M_s^V(Y)$ of maps into the singular space 
$Z_0$.  That serves as
background for our main purpose: describing the linearization 
operator ${\bf D}_\mu$  at an
approximate map into $Z_\la$ and its adjoint ${\bf D}^*_\mu$.  These 
are the operators that will
be used in Section 9 to correct the approximate maps into holomorphic maps.

\bigskip

To begin with, let $(f,j)\in \M^V_s(X)$, and consider the
linearization (\ref{2.firstlinearization}) with the norms
(\ref{defofnorms}) and
(\ref{defofnorms2}). It is convenient to decompose $\xi$ into $\zeta$ and
$\ov{\xi}$  as in (\ref{def.zeta}) and to consider  this linearization
as the operator
\bear
{\bf D}_{(f,j)}:L^{1}_s(\Lambda^0(f^*TX))\oplus T_pV^\ell \oplus
T_C\M_{g,n+\ell}\ra L_s(\Lambda^1(f^* TX))
\label{def.DD}
\eear
defined in terms of the operator (\ref{2.firstlinearization}) by
\bear
{\bf D}_{(f,j)}(\zeta,\ov{\xi},h)=D_{(f,j)}\l(\zeta+\sum\b_k \ov{\xi}_k, h\r).
\label{7.2}
\eear
Note that if $(\xi,h)$ satisfies $D_f(\xi,h)=0$ then the
condition that $\xi^N$ has a zero of order $s$ at $x$ is equivalent to
$\|\xi\|_1 <\infty$. One of the implications is obvious, while the
other follows from $\|\rho^{-2\delta}W\zeta\|_{L^4}\leq
\|\xi\|_1<\infty$ by basic elliptic estimates. This means that with
the norms above the domain of ${\bf D}_f$ includes already the
linearization of the contact conditions. Moreover, for generic
$0<\de<1$, ${\bf D}_f$ is a bounded linear Fredholm operator with
respect to these norms and models the space $\M_s^V(X)$. In
particular, for generic $V$-compatible $(J,\nu)$ $\M_s^V(X)$ is an  orbifold of
dimension $2\,{\rm ind}_\cx\; {\bf D}_f$ and we have the following two facts.

\medskip

\begin{lemma} (a)\ \ ${st\ti ev}:\M_s^V(X) \to \M_{g, n+\ell}\ti  V^\ell$  is
a smooth map of Banach
manifolds.\\

(b)\ \ The `leading coefficient map'  (\ref{akbkintrinsically}) defines a
smooth section of the bundle
$\ma\oplus_{k=1}^\ell {\cal L}_k\otimes {\cal N}_X$ over $\M_s^V(X)$.

\label{L.Ms.Banach}
\end{lemma}
\pf (a)\ It suffices to show that the linearization is a smooth map
everywhere. Let $(\xi,h)$ be a tangent vector to $\M_s^V(X)$ and decompose
$\xi=\zeta+\beta \ov\xi$ with $\ov\xi=\xi(x)\in T_pV$. The
linearization of the map $st\ti ev$ is $(\xi,h)\ra (h,\ov\xi)$ which
is obviously smooth with our norms.

\medskip

(b)\   Choose a path $(f_t,j_t)$ in $\M_s^V$ and let
$(\xi,k)$ be its tangent vector at $t=0$. Assume for simplicity that
$\ell=1$ and let  $f_t^N=a_tz^s+O(|z|^{s+1})$ be the expansion near  the
single point $x$ with $f(x)\in V$.  Writing $\xi=\zeta+\beta \ov\xi$
with $\ov\xi=\xi(x)\in T_pV$ and differentiating, we see that
  the tangential component $\dot{\zeta}\in T V$  vanishes and the 
normal component is
$\dot\zeta^N=\dot a_t z^s+O(|z|^{s+1})$. Since  $D_f(\xi,h)=0$, 
elliptic bootstraping
gives $|\dot a_t|\le
|\dot\zeta^N|_{C^s}\leq c\|(\xi,h)\|_1$. \qed

\bigskip

For maps $f_0$ into $Z_0$ the linearization of the 
$(J,\nu)$-holomorphic map equation  has a form
similar to (\ref{def.DD}), as follows.  Thinking of   $f_0$ as a pair of maps
  $(f_1,f_2)\in \M_s^V(X)\ti_{\ev}\M_s^V(Y)$, a variation
$\xi$ of $f_0$ consists of continuous sections on each component of the domain
which have the same value on both sides of each node. This means that
the domain of $D_0$ consists of sections
$\xi=(\zeta_1,\ov{\xi}_1,h_1;\zeta_2,\ov{\xi}_2,h_2)$ with the
matching condition $\ov{\xi}_1=\ov{\xi}_2$ at each node in
$V$. The corresponding operator ${\bf D}_0$ whose kernel models
$T_{f_0}\M_s(Z_0)=T_{(f_1,f_2)} \M_s^V(X)\ti_{\ev}\M_s^V(Y)$ is
\bear
{\bf D}_0:L^{1,p}_s(\Lambda^0(f_0^*TZ_0))\oplus
T_{p} V^\ell\oplus T_{C_1}\wt \M\oplus T_{C_2}\wt \M\ra
L^p(\Lambda^1(f^*_0 TZ_0))
\label{def.DD0}
\eear 
where $\Lambda^i(f_0^*TZ_0))$ means $\Lambda^i(f_1^*TX))\oplus
\Lambda^i(f_2^*TY))$. Again, one can verify that the evaluation map
$ev:\M_s^V\ti \M_s^V \to V\ti V$ is smooth and its image is
tranverse to the diagonal $\Delta$ for generic $(J,\nu)$. Thus
generically $\cok {\bf D}_0=0$ and the space $\M_s^V\ti_{ev} \M_s^V
=ev^{-1}(\Delta)$ is a smooth orbifold as in Lemma \ref{lemma2.4}.

\bigskip

The space of stable maps is defined as the set of
$(J,\nu)$-holomorphic maps {\em modulo diffeomorphisms}.  The compute
a linearization, we choose a path in the moduli space, lift to a local
slice to the action of the diffeomorphism group, and differentiate.
In fact the constructions of Section 4 provide such a local slice at
the approximate maps of Definition \ref{defApproxMaps}.  We will
describe the slice, then use it to compute the linearization operator.

Recall that a $(J,\nu)$-holomorphic map is a pair $(f,\phi):\Sigma\to
Z_\la\ti \ov{\cal U}_{g,n}$ where $\Sigma$ is a smooth 2-manifold and
$\phi$ identifies $\Sigma$ with a fiber  of the universal curve
$\ov{\cal U}_{g,n}$ and where $f$ satifies the $(J,\nu)$-holomorphic map
equation with respect to that complex structure. Given an approximate 
map $(F, \phi):\Si \to Z_\la$ we can construct
1-parameter families of deformations $(F_t,\phi_t)$ as follows.
Fix a section  $\xi$ of $F^*TZ_\la$ over $\Si$ and a vector $v\in
T_C\ov{\M}_{g,n}$  tangent to the space of stable curves at
$C_\mu=\phi(\Si)$. As in section \ref{spaceofcurvessection},
the path $\mu_t=\exp(tv)$ lifts to a path of diffeomorphisms
$\exp(t\tilde{v}):C_\mu\to C_{\mu_t}$.  This gives the family of maps
\bear
(F_t,\phi_t)\ =\ (\exp_{F}\xi,\ \phi\circ \exp^{-1}(t\tilde{v}))
\label{defomationsbyexp}
\eear
and a path $j_t=[\exp(t\tilde{v})]^*j$ of complex structures with
initial tangent $h_v\in \Omega^{01}(TC)$.  By Lemma 
\ref{hdistancessamelemma}  the norm
(\ref{wtedhnorm}) of this $h$ is uniformly equivalent to $|\dot{\mu}|$.

\medskip

That understood,  the linearization at the approximate map $F_\mu=F_{f,\mu}$ is
then given by (\ref{2.firstlinearization}) with $h=h_v$ as above.  As
before, the decomposition (\ref{def.zeta}) of $\xi$ into $\zeta$ and
$\ov{\xi}$ allows us to consider  the linearization as the operator
\bear
{\bf D}_\mu:L^{1}_s(\Lambda^0(F_\mu^*TZ_\la))\oplus T_pV^\ell \oplus
T_C\M_{g,n}\ra L_s(\Lambda^1(F^*_\mu TZ_\la))
\label{def.DD.mu}
\eear
defined by (\ref{7.2}) with $(f,j)$ replaced by the approximate map $F_\mu$.

\begin{lemma} For generic $0<\delta<1$, (\ref{def.DD.mu}) is a
bounded Fredholm operator for each approximate map $F_\mu$   with
\best
\mbox{\rm index}_{\cx }\; {\bf D}_\mu = \frac12\,\dim 
\M_s^V(X)\ti_{ev}\M_s^V(Y)
\eest
\label{Dis.bounded}
\end{lemma}
\vskip-.2in
\pf  Combining (\ref{def.DD.mu}) and (\ref{def.L}) and noting that
$|df|\le c\rho$ and  gives the pointwise bound
\best
|{\bf D}(\zeta,a,h)|\ =\ |L(\zeta)+aL(\beta) +Jdf h|
\ \le\  c\l(|\nabla \zeta|+|\zeta|+ |a|+\rho |h|\r).
\eest
It also shows that  the normal component near each node is
\best
[{\bf D}(\zeta,a,h)]^N=\ov{\partial} \zeta^N+(\nabla_{\zeta^N} J)^N\circ
df\circ j+
(\nabla_{\zeta^V+a} J)^N\circ df^N\circ j+J\circ df^N(h)\ +\ O(\rho)
\eest
(the missing term $(\nabla_{\zeta^V+a} J)^N\circ
df^V$ vanishes at $x_i$ because $a$ is tangent to $V$ and $V$ is
$J$-holomorphic).
Because $\nabla J$ is bounded and $df^N=s(z^{s-1}dz,w^{s-1}dw)+O(\rho^{s})$
this gives
\best
|W{\bf D}(\zeta,a,h)^N|\le c\,(\,|W\ov{\partial}\zeta^N |
+|W\zeta|+ |a|+\rho |h|).
\eest
Differentiating $W^{-1}\zeta^N=(z^{s-1}\zeta^x,w^{s-1}\zeta^y)$ and
noting that $|d z|, |dw|\leq \rho$ shows that
$|W\ov{\partial} \zeta|\le c(|\nabla W\zeta|+|W \zeta|)$.  Hence we
have  the pointwise bound
\best
|W{\bf D}(\zeta,a,h)|\le c\l(|\nabla (W\zeta)|+|W \zeta|+ |a|+\rho |h|\r).
\eest
   Integrating and using the Sobolev embedding on the $h$
shows that ${\bf D}$ is bounded as stated. The fact that ${\bf D}$ is
Fredholm follows from \cite{lockhard}. \qed

\bigskip\bigskip

The adjoint ${\bf D}^*$ of (\ref{def.DD.mu}) with respect to
the weighted $L^2$ norms is determined by the relation
$$
\lg(\zeta,a,h), {\bf D}^*\eta\rg\ = \
\lg{\bf D}(\zeta,a,h),\,\eta\rg
$$
Fixing the map $F$ and putting in ${\bf D}(\zeta,a,h)=L(\zeta+\b a)+J
F_*h$ as in (\ref{def.L}) one finds that
\bear
{\bf D}^{*}\eta\ =\  D^{*}\eta+A\eta+ B\eta
\qquad\mbox{where}\qquad
\left\{
\begin{array}{l}
D^*\eta \ =\ \l(\rho^{-\delta}W\ov W\r)^{-1}L^{*}(\rho^{-\delta}W\ov W\eta)
\\
A\eta\ =\ \int_{C_\mu} \rho^{-\de}\; \left[ (\ov\partial \b)\eta^V+ \b
\lg \nabla J\circ df\circ j,\eta^V \rg \right] \\
B\eta\ =\ -F_*^tJ\eta.
\end{array}\right.
\label{def.D*}
\eear
where $W$ is given by (\ref{def.W}) (and $\ov W$ is the
corresponding weighting by $\ov z$ and $\ov w$), $F_*^t$ is the
transpose of the differential of $F$ and $L^*$ is the  $L^2$ adjoint of
operator $L$ of (\ref{2.L=d+T})
\bear
L^*\eta\ =\ \partial_{f}^*\eta +S^*\eta+T^*\eta
\label{def.L*}
\eear
where $S^*$ and $T^*$ are the
adjoints of $S$ and $T$, and $\partial_{f}^*=- \sigma_J^*\circ\nabla $ where
$\sigma_J^*$ is the  adjoint of the symbol of $\ov\partial_f$.

\vskip.4in


\setcounter{equation}{0}
\section{The Eigenvalue Estimate}
\bigskip

We now come to the key analysis step: obtaining estimates on the
linearization $D$ of the $(J,\nu)$-holmorphic map equation along the
space of approximate maps.  We establish a lower bound for the
eigenvalues of $DD^*$ and construct a right inverse $P$ for $D^*$.
This operator $P$ will be used in the next section to correct
approximate maps to true holomorphic maps.

\bigskip

To get uniform estimates we fix $(J,\nu)$ generic in the sense of
Lemma \ref{lemma2.4}.  We continue to work with $\de$-flat maps, which
we will call  $\de_0$-flat in this section to  avoid confusion with
the exponential weight $\de$ of the norm (\ref{norms}), which will also
appear.  As in (\ref{def.K})  this $\de_0$ defines a compact set ${\cal 
K}_{\de_0}\subset
\M_s^V(X)\ti_{ev}\M_s^V(Y)$ of (\ref{def.K}) and corresponding subsets
\bear
{\cal AM}_\la^{\de_0}\subset {\cal AM}_\la
\qquad\mbox{and}\qquad
{\cal A}_\la^{\de_0}\subset {\cal A}_\la
\label{8.delta}
\eear
of the model space and the space of approximate maps. Thus
${\cal AM}_\la^{\de_0}$ is the inverse image of ${\cal K}_{\de_0}$ under the
covering map of Definition \ref{defofAsla} and ${\cal A}_\la^{\de_0}$ is
the image of ${\cal AM}_\la^{\de_0}$ under the gluing map
(\ref{firstGamma}).  For the maps $(f_1, f_2)$ in ${\cal K}_{\de_0}$ the
leading coefficients $|a_k|,\; |b_k|$ at the nodes are uniformly
bounded away from 0 and $\infty$, and therefore $|\la|$ is uniformly
equivalent to $|\mu_k|^{s_k}$ for each $k$.  That understood, the aim
of this section is to prove the following analytic result.

\begin{prop} There are constants $E,c>0$ independent of $\la$ and of
$f\in {\cal K}_{\de_0}\subset \M_s^V\ti_{ev}\M_s^V$ such that the linearization
$D_\mu$ at an approximate map $F=F_{f,\mu}$  has a
partial right inverse
\best
P_\mu :L_s^0(\Lambda^{0,1}(F^*TZ_\la))\ra
L^{1}_s(\Lambda^0(F^*TZ_\la))\oplus T_{C_\mu}\M
\eest
such that
\bear
c\|\eta\|_0\le \|P_\mu \eta\|_{1} \le E^{-1}\|\eta\|_0
\label{Pmu.bounded}
\eear
\end{prop}
\pf By the spectral theorem for elliptic operators, the domain of
$D^*$ decomposes as the direct sum of finite-dimensional eigenspaces
of $D^*D$ and the target similarly decomposes into eigenspaces of
$DD^*$.  The eigenvalues are non-negative and the eigenfunctions are
smooth.  Lemma \ref{lowevaluebound} below
shows that there is a uniform lower bound $E$ on the first eigenvalue
of $D_\mu D_\mu^*$ for approximate maps $F_\mu$. Using
that, Lemma \ref{bootstraplemma} shows that $D_\mu D_\mu^*$ is uniformly
invertible. Therefore $P_\mu=D_\mu^*(D_\mu D_\mu^*)^{-1}$ is a
partial right inverse for $D_\mu$ that satisfies the required
estimate. \qed

\bigskip

Let $N_k$ be the neck  region  defined by $\rho\leq 1/k$. We start by
proving the following essential estimate:
\begin{prop} For $\delta>0$ small there  are constants $k_0$ and $c$
such that for all $\la$ sufficiently small,  all approximate maps
$F\in {\cal A}_\la^{\de_0}$ and each neck $N_k$ with
$k\geq k_0$, each $\eta\in\Omega^{0,1}(F_\mu^*TZ_\la)$ satisfies
\bear
\int_{N_k} \rho^{\delta}\,\l(|\nabla\eta|^2\ + |\eta|^2\ \r)
\ \leq\ c \int_{N_k} \rho^{\delta}\, |{L^{*}\eta}|^2
+c \int_{\bd N_k}\rho^{\de}\l(|\nabla\eta|^2\ + |\eta|^2\ \r).
\label{Fourierbound1}
\eear
\label{L.Fourier0}
\end{prop}
\pf For $\delta>0$ write $\rho^{\delta}$ as the derivative of
$\psi(t)=\int_0^t \rho^{\delta}(\tau)\, d\tau$  and  integrate by parts:
$$
\int_N \rho^{\delta}\,|\eta|^2
\ =\ \int_N \psi'\,|\eta|^2\ dt\,d\theta
\  \leq \ \int_N|\psi|\cdot
2\langle\eta,\nabla \eta\rangle+ \int_{\bd N} |\psi|\cdot|\eta|^2.
$$
Because $\rho^2=2|\mu|\cosh(2t)$ satisfies $\rho^2\leq 2|\mu| e^{2t}\leq
2\rho^2$, we have $|\psi|\leq
c\rho^{\delta}/|\delta|$, so  the first integrand on the right is bounded
by $\frac12
|\eta|^2 + c_\delta \,|\nabla \eta|^2$.
Rearranging gives
\bear
\int_N \rho^{\delta}\,|\eta|^2
\  \leq\  c \int_N \rho^{\delta}\,|\nabla\eta|^2\ +\  c \int_{\bd N}
\rho^{\delta}\,|\eta|^2.
\label{Fourierboundlemma1}
\eear

Now on the cylinder $N_k$ every $(0,1)$ form $\eta$ can be written
$\eta=\eta_1\,dt-(J\eta_1)\,d\theta$ where $\eta_1$ is a
section of the pullback tangent bundle. Denote by $U=F_*\partial_t$
and $V=F_*\partial_\theta$. Note that in the usual coordinates both
$|\nabla J|$ and $|\nabla \nu |$ are bounded, so when translating into
cylindrical coordinates on the domain we get $|dF|\leq c\rho$ and thus
$$
L^*\eta\ =\  -\nabla_U\eta_1+J\nabla_V \eta_1+O(\rho|\eta|)
$$
Therefore
\best
c(|\rho\eta|^2+\rho\,|\eta| \,|\nabla\eta|)+
|L^*\eta|^2 &\ge&  |\nabla_{U}\eta_1|^2\,+\,
|J\nabla_{V}\eta_1|^2
\,-\,2\langle\nabla_{U}\eta_1,\,J\nabla_{V}\eta_1\rangle \\
&\geq &  \frac12 |\nabla\eta|^2
\,-\,2\langle\nabla_{U}\eta_1,\,J\nabla_{V}\eta_1\rangle
\eest

Next, differentiating  the  1-form
$\w=\langle \eta_1, J\nabla_{U}\eta_1\rg\,dt+
\lg\eta_1, J\nabla_V\eta_1\,\rg d\theta$
and moving $J$ past $\nabla$
\best
d\w & = & \l( 2 \langle \nabla_{U}\eta_1, J\nabla_{V}\eta_1\rangle
\,+\, \langle \eta_1,
\nabla_{U}(J\nabla_{V}\eta_1)-\nabla_{V}(J\nabla_{U}\eta_1)\rangle \r)\
\ dt d\theta \\
     & \geq  &  2 \langle \nabla_{U}\eta_1, J\nabla_{V}\eta_1\rangle
\,+\, \langle \eta_1, J{\cal R}(U, V)\eta_1\rangle -c\,|\nabla
J|\,|dF|\,|\eta|\,|\nabla \eta|
\\
&\ge &2 \langle \nabla_{U}\eta_1, J\nabla_{V}\eta_1\rangle
\,+\, \langle \eta_1, J{\cal R}(U, V)\eta_1\rangle
-c\,\rho\,|\eta|\,|\nabla \eta|
\eest
where ${\cal R}$ is the curvature of $\nabla$. Combining the last two
displayed equations, multiplying by $\rho^{\delta}$, integrating
by parts and using the bound
$2\rho|\eta|\,|\nabla\eta|\le \rho|\nabla\eta|^2+\rho|\eta|^2$ then gives
\bear
\frac12\int_{N_k} \rho^{\delta}|\nabla\eta|^2
&\leq& \int_{N_k}  \rho^{\delta} \l[|L^*\eta|^2
+ \langle {\cal R}(U, V)\eta_1, J\eta_1\rangle\r]
- d\l(\rho^{\delta}\r) \wedge \w
+\ma\int_{\partial {N_k}} \rho^{\delta} \w\;
\label{Fourierboundlemma2}
\\ \nonumber
&&+c\int_{N_k}\rho^{\de+1}(|\eta|^2+|\nabla\eta|^2)
\eear
Because the domain metric is flat, ${\cal R}$ is the curvature of $Z_\la$.
By the Gauss equations
$$
\langle {\cal R}(U, V)\eta_1, J\eta_1\rangle
\ =\ \langle R^Z(U, V)\eta_1, J\eta_1\rangle
-\langle h(\eta_1,V), h(J\eta_1,U)\rangle
+\langle h(J\eta_1,V), h(\eta_1,U)\rangle
$$
where $R^Z$ is the curvature of $Z$ and $h$ is the second fundamental
form of $Z_\la\subset Z$, which satisfies $|h(F_*v,\cdot)|\leq c|v|$
for any $v$.  Since $R^Z$ is bounded then the term containing it is
dominated by $c\rho^2|\eta|^2$.  Also, as in Lemma
\ref{Boundondg-nulemma} $|V-JU|=|\ov\partial F|\leq c\rho$. Hence we
can replace $V$ by $JU$ with small error: 
\bear 
\langle {\cal R}(U, V)\eta_1, J\eta_1\rangle \ \leq\ -\langle h(\eta_1,JU),
h(J\eta_1,U)\rangle +\langle h(J\eta_1,JU), h(\eta_1,U)\rangle
+c\rho|\eta|^2.
\label{Rheq}
\eear
Observe that if we had $\nabla J=0$ along $Z_\la$ then $h$ would be linear
in $J$ and the two $h$ terms above would reduce to
$-2|h(\eta_1,U)|^2\leq 0$. In our case the $\nabla J$ term is of
order $\rho|\eta|$ therefore
\bear
\langle {\cal R}(U, V)\eta_1, J\eta_1\rangle
\ \leq\ c\rho|\eta|^2.
\label{finalRbound}
\eear
which can be absorbed in the last term of (\ref{Fourierboundlemma2}).

It remains to bound the $\w$ term in (\ref{Fourierboundlemma2}). As in
(\ref{gmuvsg}) we can introduce cylindrical coordinates
$\tau=\log|x/\la|$ and $\Theta$ on $N_X$ and normal (Fermi)
coordinates in the $V$ direction. Then the metric on $Z_\la$ is
$R^2(d\tau^2+d\Theta^2)+g^V$ where $R^2=|x|^2+|y|^2
=2|\la|\cosh(2\tau)$ and $g$ is the metric of $V$. The formula for $F$
shows that in this basis
$F_* \partial_\theta= s\partial_\alpha+O(\rho)\partial_i$ and a
computation shows that the
Christoffel symbols are all bounded and those in the $\Theta$ direction are
\best
\Gamma^\Theta_{\Theta \Theta}=\Gamma^\Theta_{\tau \tau}= 0\quad
\Gamma^\Theta_{\Theta\tau}= - \Gamma^\tau_{\Theta \Theta} = \tanh(2\tau).
\eest
Thus $\nabla_\theta = \partial_\theta +\tanh(2\tau) J + A\rho$ where $A$
is bounded.  Recalling the definition of $\rho^2$ from
(\ref{conformalmetrics}), we have
\bear
-d(\rho^\de)\wedge \w\ =\  -\partial_t\rho^\de\,\langle \eta_1,
J\nabla_V\eta_1\rangle\,dt
d\theta=  \delta\rho^\delta\tanh(2t) \langle J\eta_1,
\nabla_V\eta_1\rangle\,dtd\theta
\label{om.term}
\eear
Because $g_\la$ is independent of $\theta$ in these coordinates, using
the same methods as in (\ref{holderpoincare}) combined with the fact
that $|\tanh(2\tau)|\le 1$ we get the bound
\best
-\tanh(2\tau)\int_{S^1} \lg J\eta_1, \partial_\theta \eta_1\rg d\theta\;\le\;
\int_{S^1}|\partial_\theta\eta_1|^2\,d\theta
\eest
Moving all the terms on the same side we get
\best
0\le \int_{S^1}\lg \nabla_\theta\eta_1,\partial_\theta\eta_1\rg
+c\rho\int_{S^1} (|\eta|^2+|\nabla\eta|^2)
\eest
which the implies
\best
\tanh(2\tau)\int_{S^1} \lg J\eta_1, \nabla_\theta \eta_1\rg \,d\theta\;\le\;
\int_{S^1}|\nabla_\theta\eta_1|^2d\theta+
c\rho\int_{S^1} (|\eta|^2+|\nabla\eta|^2)
\eest
But $\tanh(2\tau)=\tanh(2st)+O(\rho)$ and $0\le\tanh(2t)/\tanh(2st)\le
1$ so combining the last displayed equation with (\ref{om.term}) gives
\best
-\int_{N_k} d(\rho^\de)\wedge \w \le \de \int_{N_k} \rho^\de |\nabla\eta|^2 +
c\int_{N_k} \rho^{1+\de}( |\nabla\eta|^2+|\eta|^2)
\eest

Inserting this and (\ref{finalRbound}) into (\ref{Fourierboundlemma2})
including  (\ref{Fourierboundlemma1}) gives (\ref{Fourierbound1})
for small $\delta$ and large $k$. \qed

\bigskip

Write $\nabla^\delta\eta=\rho^{\delta}\,\nabla(\rho^{-\delta}\eta)$
where $\nabla$ is as usual the covariant derivative of  the
cylindrical metric on the domain and the metric induced on $Z_\la$
from $Z$. Note that when $\delta>0$ is small, the $L^{1,2}$ weighted
norm defined using $\nabla^\delta$ is uniformly (in $\la$) equivalent
to the one using $\nabla$. Then Proposition \ref{L.Fourier0} implies:

\begin{cor} For $\delta>0$ small there  are constants $k_0$ and $c$
such that for  all $\la$ sufficiently small and  all approximate maps
$F\in {\cal A}_\la^{\de_0}$  and each neck $N_k$ with
$k\geq k_0$, each $\eta\in\Omega^{0,1}(F^*TZ_\la)$ satisfies
\bear
\|\eta\|_{1,2,N_k}\le c \|D^*\eta\|_{2,N_k}+c\|\eta\|_{1,2,\bd N_k}
\label{Fourierbound2'}
\eear
i.e.
\bear
\int_{N_k} \rho^{-\delta}\,\l(|\nabla^{\delta}W\eta|^2\ + |W\eta|^2\
\r)
\ \leq\ c \int_{N_k} \rho^{\delta}\, |WD^{*}\eta|^2
+c \int_{\bd N_k}\rho^{-\de}\l(|\nabla^{\delta}W\eta|^2\ + |W\eta|^2\ \r).
\label{Fourierbound2}
\eear
\label{Fourierboundlemma}
\end{cor}
\pf Since on each coordinate $\partial\ov W=0$ then relation (\ref{L.alm.cx})
combined with condition (\ref{defcompatibleJnu}) implies that
\best
(\ov W)^{-1}L^* (\ov W \eta)= L^* \eta + O(\rho |\eta|).
\eest
So (\ref{Fourierbound2}) follows from (\ref{Fourierbound1}) after
replacing $\eta$ by $\rho^{-\de}W\eta$ and using (\ref{def.D*}) .\qed

\bigskip

   From now on we will fix $\delta>0$ small and generic.
The following lemma can be compared to Lemma 6.6 in \cite{rt1} and 3.10 in
\cite{lt}.

\medskip

\begin{lemma} There is a constant $E>0$ such that for all $\la$
sufficiently small and all approximate maps $F\in{\cal A}_\la^{\de_0}$,
the first eigenvalue of ${\bf D}{\bf D}^{*}$ is
bounded below by $E$.
\label{lowevaluebound}
\end{lemma}
\pf Suppose the claim is false.  Then there are sequences $\la_n,\mu_n\to
0$, maps $F_n:C_{\mu_n}\to Z_{\la_n}$ in some ${\cal K}_{{\de_0}}$ and
$(0,1)$ forms $\eta_n$ along $F_n$ with
${\bf D}_n{\bf D}_n^{*}\eta_n=\ep_n \eta_n$ with $\ep_n\ra 0$.
In particular,
\bear
\ep_n \int \rho^{-\de}|W \eta_n|^2\  \ \geq \ \, \
\int \rho^{-\de}|W {D_n^{*}\eta_n}|^2\;+|A\eta|^2+|B\eta|^2.
\label{evalueinequalitytobecontradicted}
\eear
where $A, B$ as in  (\ref{def.D*}). We may normalize the
$\eta_n$ so that the lefthand side of
(\ref{Fourierbound2'}) is one. By the Bubble Tree Convergence Theorem
there is a subsequence of the $F_n$ that converges to a stable map
$F_0$ from $C_0=C_1\cup C_2$ into $Z_0$, and this convergence is in
$C^\infty$ away from the nodes. On small compact sets $K$ in
the complement of the nodes, the $L^{1,2}_s$ norm in the
cylindrical metric is uniformly equivalent to the usual $L^{1,2}$
norm.  Standard elliptic theory implies that there is a subsequence of
the $\eta_n$ that converges in $C^\infty$ on $K$ to an $L^{1,2}_s$
section with $D_0^*\eta=0$ along $F_0\setminus K$.  Doing this for the
sequence
$K_m=\rho^{-1}([\frac1m,\infty))$ and passing to a diagonal
subsequence yields a limit $\eta$ defined on
$C_0\setminus\{\mbox{nodes}\}$ with $L^{1,2}_s$ norm at most
one, and such that  $D_0^*\eta=0$ along $F_0$ outside the nodes. Moreover,
${\bf D}_0^*\eta =0$ weakly,
i.e. for all $\zeta\in L^{1,2}_s$, $a\in T_pV$ and $v\in T_C\M$
\best
\lg {\bf D}_0(\zeta, a, v),\eta\rg=0
\eest
on $C_0$. We show this for $a\in TV$, the other parts being similar. On
$C_{\mu_n}$
\best
\lg {\bf D}_n(a),\eta_n\rg=\lg a, A\eta_n\rg_{TV}\ra 0
\eest
and ${\bf D}_n(a)=D_n(\b a)=(\ov\partial \b)a+\b \nabla_a J\circ dF_n$. Off
each neck $N_k=\{\rho< 1/k\}$
\best
\int_{C_\mu\setminus N_k} \rho^{-\de}\lg W{D_n(\b a)},\;W \eta_n \rg
\ra
\int_{C_0\setminus N_k}\rho^{-\de}\lg W D_0(\b a),\;W\eta \rg
\eest
while on $N_k$ $D_0(\b a)=D_0(a)= (\nabla_a J)\circ dF_n\circ
j$. So
\best
|\lg D_0(\b a),\eta_n\rg_{N_k}|=
\l|\int_{N_k}\rho^{-\de}\lg W{D_0(\b a)},W\eta_n \rg\r|\le
\l(\int_{N_k}\rho^{-\de}|W(\nabla_a J)\circ d F_n\circ j|^2\r)^{1/2}
\cdot \|\eta_n\|_{2,s}
\eest
But $V$ is $J$-invariant, so
$W{D_0(\b a)}=(\nabla_a J)\circ d (W F_n)\circ j+O(\rho)$ and
$|\nabla_aJ|$ is bounded and $|d{W F_n}| \le \rho^{1/3}$, so
the  first factor on the right hand side of the last displayed
equation goes to zero as $k\ra \infty$.

This means that ${\bf D}_0^*\eta=0$ where ${\bf D}_0$ is the
operator defined in (\ref{def.DD0}).   As we observed after equation 
(\ref {def.DD0}),
for generic $(J,\nu)$ we have $\cok {\bf
D}_0=0$, so $\eta=0$. Therefore  $\eta_n \to 0$ in $L^{1,2}$ on the
complement of each neck $N_k$,
which contradicts  (\ref{Fourierbound2'}).  \qed

\bigskip

\begin{lemma}
There is a constant $C$ such that for all $\la$ sufficiently small and all
  approximate maps $F\in {\cal A}_\la^{\de_0}$, each
$\eta\in\Omega^{0,1}(F^*TZ_\la)$
satisfies $\|\eta\|_2\ \leq\ C\, \|{\bf D} {\bf D}^*\eta\|_0.$
\label{bootstraplemma}
\end{lemma}
\pf  Cover $C_\mu$ by disks of radius 1 in the cylindrical metric so that
each point lies in at most 10 disks. Since $\rho$ varies by a bounded
factor
across each unit interval in the neck we can applying the basic
elliptic estimate on each disk, multiply by $\rho^{-\delta/2}$ and
sum to get
\best
\|\eta\|_{2,p,s}\ \leq\ C_p\,\l(\|D_\mu D_\mu^*\eta\|_{p,s}+
\|\eta\|_{p,s}\r)
\eest
for a constant $C_p$ independent of $\mu$. Adding together the $p=4$ and
$p=2$ inequalities we get
\best
\|\eta\|_2\ \leq\ C_p\,\l(\|D_\mu D_\mu^*\eta\|_0+\|\eta\|_0\r)
\eest
Using Lemma \ref{lowevaluebound} and applying  Holder's inequality for the
weighted $L^2$ norm
$$
c\; \|\eta\|^2_{2,s}\ \leq\ \|D^*\eta\|^2_{2,s}\ = \ \langle
\eta,D_\mu D_\mu^*\eta\rangle\  \leq\ \|\eta\|_{2,s}\ \|D_\mu
D_\mu^*\eta\|_{2,s}  \leq\
\|\eta\|_{2,s}\ \|D_\mu D_\mu^*\eta\|_{p,s}.
$$
which combined with the previous inequality gives the desired inequality.
\qed

\bigskip

\bigskip

In the next section we will use $P_\mu$ to coordinatize the normal
direction to the space of approximate maps.

\vskip.4in


\setcounter{equation}{0}
\section{The Gluing Diffeomorphism}
\label{GluingDiffeomorphismsection}
\bigskip

\bigskip

The norm (\ref{defofnorms}) induces a topology on the space
$\mbox{Maps}_s(C,Z_\la)$. Specifically, for
$C^0$ close maps with the same label $s$ we can write 
$(f',j')=\exp_{(f,j)}\,(\xi,h)$ and set
\bear
\mbox{dist}\l((C,f), (C',f')\r) \ =\ \|(\xi, h)\|_1
      \label{distancebetweenmaps}
\eear
This defines a topology and a distance (the inf of the lengths over all
paths piecewise of the above
type) on $\mbox{Map}_s(C,Z_\la)$.  Using this distance, we will show that
the moduli space of stable maps into $Z_\la$ is close to the space of
approximate maps, and that those spaces are in fact  are isotopic.

\bigskip

We start by  describing a parameterization for a neighborhood of
${\cal A}_\la^{\de}$
inside the space of maps (${\cal A}_\la^{\de}$ is the compact set
(\ref{8.delta}) of approximate maps).  Consider the Banach space
bundle $\Lambda^{01}\to{\cal AM}_{s}(\la)$  over
${\cal AM}_{s}^\de(\la)$ (the model space for approximate maps)
whose fiber at an approximate map $F_\mu:C_\mu\to Z_\la$ is
$\Lambda^{0,1}(F^*_\mu TZ_\la)$ with the norm (\ref{defofnorms}). Write
elements of $\Lambda^{01}$ as triples $(f,\mu, \eta)$, with
$f\in \M_s^V(X)\ti_{ev} \M_s^V(Y)$.  The map
\bear
\Phi_\la:\Lambda^{01}(\ep)\ra \mbox{Maps}_s(C,Z_\la\times {\cal U})
\qquad{\mbox{by}}\qquad
\Phi_\la(f,\mu,\eta)\ =\ \exp_{F_{f,\mu},C_\mu}(P_\mu\eta)
\label{defPhi}
\eear
defined on an $\ep$ neighborhood of the zero section of $\Lambda^{01}$
agrees with the gluing map $\Gamma_\la$ along the zero section.
The following lemma shows that $\Phi_\la$ coordinatizes a
neighborhood of ${\cal A}_\la^{\de}$.

\medskip

\begin{prop} There is a constant $c>0$ so that for all small
$\la$ $\Phi_\la$ is a diffeomorphism from an $\ep$- neighborhood of
the zero  section in $\Lambda^{01}$ onto a neighborhood of
${\cal A}_\la^\de$ in $\mbox{Maps}_s(C,Z_\la\times {\cal U})$ that contains
at least a $c\ep$ neighborhood of ${\cal A}_\la^\de$.
\label{exp map}
\end{prop}
\pf  By Lemma \ref{lemma2.4} $T_{F_\mu}{\cal A}_\la$
has the same dimension as $\ker {\bf D}_\mu=\l(\mbox{Im}\,P_\mu\r)^\perp$.
In fact,
\bear
T_{F_\mu}\Lambda^{01}\ =\ T_{F_\mu}{\cal A}_\la \oplus \mbox{Im}\,P_\mu
\label{ETAplusL}
\eear
because any  $P_\mu \eta$ which lies in $T_{F_\mu}{\cal A}_\la$ satisfies,  by
(\ref{Pmu.bounded}), Lemma \ref{lambdazlemma} and Lemma
\ref{TAalmostkernel} below,
\best
\|P\eta\|_1\ \leq\  E\,\|\eta\|\ =\ E\,\|{\bf D}_\mu P\eta\|\ \leq\
CE\,|\la|^{1/8|s|}\,\|P\eta\|_1,
\eest
so, for small $\la$, $P\eta$ is zero.
\smallskip

Next fix a path $(f_t,\mu_t)$ in ${\cal AM}_s(\la)$ starting at
$(f_0,\mu_0)$ and let $\xi\in T{\cal A}_\la$ the tangent vector at
$t=0$ of the corresponding path of approximate maps
$F_t=\Phi(f_t,\mu_t,0)$.  Each element $\tau$ in
the fiber of $\Lambda^{01}$ over
$(f_0,\mu_0)$ determines a vector field $P\tau$ along the image of  $F_0$
in $TZ_\la$. After  extending
$P\tau$ along $F_t$ by parallel translation we calculate
$$
d\Phi_{(f,\mu,\eta)}(\xi,h,\tau)\ =\ \l.\frac{d\,}{dt}\
\exp_{(F_t, j_\mu)}(tP\tau)\phantom{\int}\hspace{-.5em}\r|_{t=0}
\ =\ \xi+P\tau.
$$
Thus $d\Phi_\la$ is an isomorphism by (\ref{ETAplusL}), so $\Phi_\la$ is a
local diffeomorphism near the zero section of $\Lambda^{01}$.

To show injectivity, let $\Lambda^{01}(\ep)$ be the subset of
$\Lambda^{01}$ with $\|\eta\|\leq \ep$ and suppose that injectivity
fails on each $\Lambda^{01}(\ep)$.  Then for each $n$ there exist
elements $(f_n,\mu_n,\eta_n)\neq (f_n',\mu_n',\eta_n')$ in
$\Lambda^{01}(1/n)$ which have the same image under $\Phi_\la$.  After
passing to subsequences, we can assume that the $\{(f_n,\mu_n)\}$ and
$\{(f_n',\mu_n')\}$ converge in the stable map topology to limits
$f: C\to Z_0$  and $f': C'\to Z_0$ with
$f$ and $f'$ in ${\cal K}\subset \M_s^V\ma\ti_{\tiny
ev}\M_s^V$ and $C$ and $C'$ on the boundary of the
cylindrical end compactification of $\M_{g,n}$ defined at the end of Section
\ref {spaceofcurvessection}. (Thus  $C$ and $C'$ each consist of a
nodal curve together with an element of the real torus $T^\ell$).

Choose a compact region $R$ in $C$ which contains no nodes.  Then  for
small $\la$ we have
$F_n\ra f$ and $F_n'\ra f'$ in $C^1$ on $R$. Since our
$\|\cdot\|_1$ norm dominates both the $C^0$ norm on maps and, by Lemma
\ref{hdistancessamelemma}, the
cylindrical end metric on $\M_{g,n}$,
\best
\ma\lim_{n\ra \infty} \mbox{dist }(C_{\mu_n}, C_{\mu_n'}) \ + \
\sup_{x\in R}\  \mbox{dist }(f(x), f'(x)) &\leq&
\ma\lim_{n\ra\infty}\ (\|P\eta_n\|_1+\|P\eta_n'\|_1)
\\
&\le&
c \ma\lim_{n\ra\infty}\ (\|\eta_n\|+\|\eta_n'\|)\ =\ 0
\eest
using Lemma \ref{lowevaluebound}.  Thus (i) \ $C=C'$, and (ii)\ \
$f$ and $f'$ agree on $R$ and therefore, as in the argument after
(\ref{firstGamma}), agree everywhere. Consequently, for large $n$
$(f_n,\mu_n,\eta_n)$
and
$(f_n',\mu_n',\eta_n')$ lie in the region where $\Phi_\la$ is a local
diffeomorphism and  are therefore equal.  That establishes
injectivity. The surjectivity onto an $c\ep$ neighborhood folows from
the first inequality in (\ref{Pmu.bounded}). \qed

\bigskip

The norms (\ref{defofnorms}) for $\la=0$ induce a Banach manifold
structure on $\M_s^V(X)\ti_{ev}\M_s^V(Y)$, and hence on its cover the
model space ${\cal AM}_s(\la)$.  But the gluing map identifies
${\cal AM}_s(\la)$ with the space of approximate maps ${\cal A}_\la$, which
has a possibly different norm as a subset of the Banach space
$\mbox{Maps}_s(C,Z_\la\ti{\cal U})$. The next lemma shows that these two
norms on $T{\cal AM}_s$ are uniformly equivalent.

\begin{lemma} There are constants $c,C>0$, uniform on each compact
${\cal K}_\de\subset \M_s^V(X)\ti_{ev}\M_s^V(Y)$, so that for each tangent
vector $(\xi,h)$
to $\M_s^V(X)\ti_{ev}\M_s^V(Y)$ each of its images $(\xi_\mu,h_\mu)$
under the differential of the gluing map (\ref{firstGamma}) satisfy
\best
c\|(\xi,h)\|_1\le \|(\xi_\mu, h_\mu)\|_1\le C\|(\xi,h)\|_1.
\eest
\label{L.unif.equiv.metric}
\end{lemma}
\vskip-.2in
\pf Choose a path $(f_t,j_t)$ in $\M_s^V(X)\ti_{ev}\M_s^V(Y)$ with
tangent $(\xi,h)$ at $t=0$ and lift it to a path
$(F_t,\mu_t)\in {\cal AM}_s$ with initial tangent vector $(\xi_\mu,
h_\mu)$. By construction, 
the approximate maps agree with $(f_t,j_t)$ outside the region
$A_\mu=\{ \rho<2|\mu|^{1/4}\}$  and so $\xi_\mu=\xi$ and $h_\mu=h$ off
$A_\mu$. Moreover, $\xi$ and $\xi_\mu$ have the same
average value in $T_pV^{\ell}$, so we may assume without loss of
generality that this average value is 0. Then on $A_\mu$
$\xi_\mu=(1-\b_\mu)\xi$ while $h_\mu-h$ is of order $\dot \mu_t/\mu$
by Lemma \ref{hdistancessamelemma}. By differentiating the relation
$a_tb_t\mu_t^s=\la$ we see that
$\dot \mu_t/\mu$ is of order $\dot a_t/a+\dot b_t/b$. Integrating on
$A_\mu$ and using Lemma  \ref{L.Ms.Banach} gives
\best
\|(\xi_\mu-\xi, h_\mu-h)\|_{1,A_\mu}\le C|\mu|^{1/8}\|(\xi,h)\|_1
\eest
uniformly on the compact ${\cal K}$ (when $\delta<1/2$). \qed

\begin{lemma}
There is a constant $C$, uniform for $f_0$ in ${\cal K}_\de\subset
\M_s^V(X)\ti_{ev}\M_s^V(Y)$, such that for $\la$ small
enough the tangent vectors
$(\xi,h) \in T_{F}{\cal A}_\la$ at the approximate map $F=F_{f,\mu}$  satisfy
$$
\| {\bf D}_{\mu}(\xi,h)\|_0 \le C |\la|^{1/8s} \l(\|\xi\|_1+\|h\| \r).
$$
\label{TAalmostkernel}
\end{lemma}
\vskip-.2in

\pf This time,  choose a path  $(f_t,j_t)\in \M_s^V\ti_{ev}\M_s^V$ with
initial
tangent vector $(\xi_0,h_0)\in \mbox{\rm ker}\ {\bf D}_{f}$, lift to a path
$(F_t,\mu_t)\in {\cal AM}_s$,
and let $(\xi, h)$ be the initial tangent vector to the lifted path. On the
neck $\rho<|\mu|^{1/4}$ we
again have $\xi=(1-\b_\mu)\xi_0$.  The lemma follows  from the pointwise
estimates
\best
&&|{\bf D}_\mu(\xi,h)-{\bf D}_0(\xi,h)|\le |dF_\mu-df_0|\cdot (|\xi|+|h|)
\\
&& |{\bf D}_0(\xi,h)|\le |{\bf D}_0(\xi,h)-{\bf D}_0(\xi_0,h_0)|\le
|\nabla\b_\mu\cdot \xi_0|+|\nabla_{\xi_0} J \circ df_0\circ j|+
|J\circ df_0\circ (h-h_0)|
\eest
combined with Lemma \ref{L.unif.equiv.metric}.\qed
\bigskip

\begin{prop}
\label{nearapprox}
For each $\ep>0$, $\M_s^{flat}(Z_\la)$ lies in an $\ep$-neighborhood
of ${\cal A}_\la^\de$ for all $\la<\la_0(\ep)$.
\end{prop}
\pf As $\la_n\to 0$, any sequence $(f_n,j_n)\in
\M_s^{flat}(Z_{\la_n})$ has a subsequence which converges as in
(\ref{4.sequence}) to a limit $f_0$ from an $\ell$-nodal curve
$C_0$. Write $(\Si,j_n)=(C_n,\mu_n)$ where $C_n$ is an $\ell$-nodal
curve close to $C_0$ and choose $\mu_n'=(\mu_{n,1}',
\dots, \mu_{n,\ell}')$ with $\la_n=a_kb_k(\mu_{n,k}')^{s_k}$ for each
$k$; there are $|s|$ choices for each $\mu'_n$ which differ by roots
of unity. That data defines corresponding approximate maps
$F_n=F_{f_0,\mu_n'}:(C_0,{\mu'_n})\to Z_{\la_n}$ in ${\cal A}_{\la_n}$
via (\ref{defapproxmaps}).  We will show that for some choice of
$\mu_n'$
$$
\mbox{dist}((C_n,\mu_n),\;(C_0, \mu_n'))\,+\,\mbox{dist}
(f_n, F_n)\ <\ \ep \qquad \qquad \mbox{for large $n$.}
$$
Lemma \ref{lambdazlemma} shows that $\l({\mu_n}/{\mu_n'}\r)^s\to
1$ at each node.  After passing to a subsequence and modifying our
choice of $\mu_n'$ we have ${\mu_n}/{\mu_n'}\to 1$.  But then
(\ref{firsthdistance}) shows that $\mbox{dist}((C_n,\mu_n),\;(C_0,
\mu_n')) \to 0$.

For any $r_0<1/2$ the maps $F_n$ and $f_0$ agree on the sets $\{\rho
\geq r_0\}$ for all large $n$.  Inside the region $A(r_0)$ near each
node where $\rho\le r_0$, the vector
$\xi_n=F_n-f_n=\beta_\mu(f_0-f_n)$ can be written as
$\xi_n=(\zeta_n,\bar{\xi}_n)$ as in (\ref{def.zeta}).  But
$\bar{\xi}\to 0$ because $f_n\to f_0$ in $C^0$, and Lemma
\ref{1pconvergence} implies that $\|\zeta\|_{1, A(r_0)}\leq
cr_0^{1/6}$. Taking $r_0$ small enough and using the fact that outside
the neck $A(r_0)$ we have uniform convergence implies $\mbox{dist}
(f_n, F_n)=c(\|\zeta\|_1+|\ov\xi|)<\ep$ for large $n$.  \qed

\bigskip

The next step is to correct each approximate map $F_{f,\mu}\in {\cal
A}_\la$ to get $(F'_\mu,j_\mu')$ a true $(J, \nu)$-holomorphic map.
More precisely, $(F'_\mu,j_{\mu}')$ will
be a solution of the equation
\bear
\overline\partial_{j} f=\nu_{f} \quad \mbox{ where }\quad
(f, j)=\exp_{F_\mu,j_\mu} ( P_\mu \eta)
\label{del.eq}
\eear
and $\eta \in \Lambda^{0,1}(F_\mu^* TZ_\la)$.

\begin{prop} There are  constants $\ep, \la_0$ and $C$ (uniform on
${\cal K}_\de\subset \M_s^V(X)\ti_{ev}\M_s^V(Y)$) such that for
each $f\in{\cal K}_\de$ and $0<|\la|<\la_0$  equation  (\ref{del.eq})
has a unique solution
$\eta\in \Lambda^{0,1}(F^*_\mu TZ_\la)$ in the ball
$\|\eta\|\le \ep$, and that solution is smooth and  satisfies
$\| \eta \|\le C|\la|^{1\over 8s}$.
\label{exist.sol}
\end{prop}
\pf If  we write $(f,j)=\exp_{F_\mu,j_\mu} (\zeta )$ where $\zeta =(\xi,h)$
then
\bear
\overline\partial_j f-\nu_{f} = \overline\partial_{j_\mu}
F_\mu-\nu_{F_\mu} + {\bf D_\mu} (\zeta ) +Q_\mu(\zeta )
\label{del0}
\eear
where ${\bf D_\mu}$ is the  linearization at $(F_\mu,j_\mu)$ and  the
quadratic $Q_\mu$
satisfies (c.f. \cite{f})
\bear
\| Q_\mu(\zeta _1)- Q_\mu(\zeta _2) \|  &\le &  C \; (\; \|\zeta _1 \|_{1} +
\|\zeta _2\|_{1} ) \| \zeta _1-\zeta _2 \|_{1}
\label{floer1}
\eear
Taking $\zeta =P_\mu\eta$ and noting that $D_\mu P_\mu\eta= \eta$, equation
(\ref{del.eq}) becomes
\bear
\eta + Q ( P_\mu \eta) = v \qquad \mbox{where}\qquad v =
\nu_{\mu}-\overline\partial_{\mu}
f_\mu
\label{eq1}
\eear
Define an operator $T_\mu$ on the Banach space  obtained by  completing
$\Omega^{0,1}(f^*_\mu TX)$ in our norm (\ref{defofnorms2})   by
$$
      T_\mu \eta = v -   Q_{\mu} ( P_\mu \eta).
$$
Using (\ref{floer1})  and (\ref{Pmu.bounded})
\best
      \| T_\mu \eta_1 - T_\mu \eta_2 \| & \le&
C \;(\; \| P_\mu \eta_1 \|_{1} + \| P_\mu \eta_2 \|_{1})
      \|  P_\mu (\eta_1-\eta_2) \|_{1}          \\
      & \le &
      C\, E^2 \;( \;\| \eta_1 \| + \| \eta_2 \| ) \cdot \| \eta_1-\eta_2 \|.
\eest
Choosing $\ep < 1/(4C E^2)$, when $ \| T_\mu (0) \| \le \ep/2 $ then $
T_\mu :B(0,\ep) \ra B(0,\ep)$ is a contraction on the ball of radius
$\ep$.  Therefore $T_\mu$ has a unique fixed point $\eta$ in the ball,
and $\|\eta \|\le 2\| T_\mu(0)\|$.  Finally, since $\eta\in L^4_{loc}$
we have $\zeta =P\eta\in L^{1,4}_{loc}$ with $D\zeta +Q(\zeta )=Pv\in
C^\infty$.  Elliptic regularity then shows that $\zeta $ and $\eta$
are smooth.  \qed

\vskip.4in


\setcounter{equation}{0}
\section{Convolutions and the Sum Formula for Flat Maps}
\bigskip

We can now assemble the analysis of the previous several sections to
show that the approximate
moduli space,  which is built from maps into $Z_0$, is a good model
of the moduli
space of stable maps  into the symplectic sum $Z_\la$.  Recall that  in
    Sections 3 and 5 we showed that as $\la\to 0$ stable maps into
$Z_\la$ limited to maps into
$Z_0$ and that the complex structure $\mu$ on their domains are determined by
the limit map up to a finite ambiguity corresponding to the different
solutions of the
equation
$ab\mu^s=\la$.  That led to the definition of the model moduli space
${\cal AM}_s$ in Section 6.
On the other hand, each element of ${\cal AM}_s$ defines an
approximate holomorphic map by
equation (\ref{defapproxmaps}); for each $\la$ this gives  the gluing map
\bear
\Gamma_\la: {\cal AM}_s\ \overset{\approx}\longrightarrow\,  {\cal
A}_s(\la) \, \subset\
\mbox{Maps}_s(C,Z_\la\ti {\cal U})
\label{secondGamma}
\eear
whose image ${\cal A}_s(\la)$ we call the space of approximate maps.
And indeed,  Proposition
\ref{exist.sol} shows that each such approximate map   can be
uniquely perturbed to be true
$(J,\nu)$-holomorphic map.

    \smallskip

    In  this section we will show that ${\cal A}_s(\la)$ is  isotopic to
$\M_s(Z_\la)$  through an isotopy compatible with the evaluation
maps.   Thus ${\cal
AM}_s(\la)$ keeps track of the  fundamental homology class
$[\M_s(Z_\la)]$ which defines the GW and TW invariants of $Z_\la$ (we continue
to assume that all maps  have been stabilized as in Remark
\ref{R.exp}).  Passing to homology, we
then   define a ``convolution'' operation  and establish a formula of the form
\bear
TW_X^{V}\, *\, TW_Y^{V}\, = \, TW_Z
\label{3.1}
\eear
under the assumption that all curves contributing to the invariants are
$V$-flat (this condition will be eliminated in  Section
\ref{The General Sum Formula}).

\bigskip

     We noted in (\ref{def.K}) that as $\la\to 0$ the limits of the
$\delta$-flat  maps into
$Z_\la$ lie in the compact set ${\cal K}_\de$ of $ \M^V(X) \ti_{ev}
\M^V(Y)$.  We will work on
the corresponding compact sets ${\cal AM}_s^\de$ and
${\cal A}_s^\de(\la)$ defined in (\ref{8.delta}).

\begin{theorem}   Fix an ordered sequence $s$ and write  $|s|=\prod s_i$.
For generic $(J,\nu)$ and small $|\la|$, there is an $|s|$-fold cover
${\cal AM}_s^\de$ of ${\cal K}_\de$ and a diagram
\bear
\begin{array}{ccccc}
\ma\bigsqcup_s\,{\cal AM}_s^\de &  \ &
\underset{\Phi^1_\la}\longrightarrow & \ &
\ma\bigsqcup_s\,\Mf_s(Z_\la)\\
\Big\downarrow\vcenter{\rlap{$st$}} &\   &  &\  & \Big\downarrow
\vcenter{\rlap{$st$}}  \\
\wt{\M}\ti \wt{\M} &\ &
\overset{\xi}\longrightarrow &\ &  \wt{\M}
\end{array}
\label{10.diagram1}
\eear
where the top arrow is a diffeomorphism onto its image and is isotopic
to the restriction of (\ref{secondGamma}) to ${\cal AM}_s^\de$. The
diagram commutes up to homotopy. Furthermore, there is
a constant $c=c(\delta)$ so that the image of $\Phi^1_\la$ consists
of maps which are
$(\de-c\la)$-flat, and the image contains all $(\de+c\la)$-flat maps
in $\M_s(Z_\la)$.
\label{gluing thm}
\end{theorem}
\pf For each $(f_0,\mu)\in{\cal AM}_s$ the gluing map $\Gamma_\la$
associates a smooth curve
$C_\mu$ and an approximate map $F_\mu:C_\mu\ra Z_\la$.    By
Proposition  \ref{exp map} any pair $(f',C_{\mu'})$ that is $L^1_s$
close to $\Gamma_\la(f,\mu)$
can be uniquely written as
\bear
\Phi_\la(f,\mu,\eta)\ =\ \exp_{F_{f,\mu},C_\mu}(P_\mu\eta)
\label{10.eta}
\eear
for some $L^0_s$ section $\eta$ of the bundle $\Lambda^{0,1}$ with
$\|\eta\|<\ep$.
    Proposition  \ref{exist.sol} then used a fixed
point theorem to show that for small $|\la|$ there is a unique such
$\eta=\eta(f,\mu)$ such that
(\ref{10.eta}) is $(J,\nu)$-holomorphic.  Then
$$
\Phi^t_\la(f,\mu,\eta)\ =\ \exp_{F_{f,\mu},C_\mu}\l(tP_\mu\eta(f,\mu)\r)
$$
is a smooth 1-parameter family of maps from ${\cal AM}_s^\de$ to
$\mbox{Maps}_s(C,Z_\la\ti {\cal
U})$ with $\Phi^0_\la=\Gamma_\la$  and  the image of $\Phi^1_\la$
lying in the  $(\de-c\la)$-flat
maps in $\M_s(Z_\la)$.  The uniqueness of $\eta$ in the fibers of
$\Lambda^{0,1}$, combined with
Proposition \ref{exp map} implies that the $\Phi^1_\la$ is  injective.

\smallskip

It remains to show  that $\Phi^1_\la$ is surjective.  But
Proposition \ref{exp map} shows that
(\ref{10.eta}) is onto at least a $c\ep$ neighborhood of ${\cal
A}^{2\de}_\la$ and Proposition
\ref{nearapprox} implies that $\Mf_s(Z_\la)$ lies in that
neighborhood when $|\la|$ is small
enough.  Hence  for $|\la|$ small, each element of $\Mf_s(Z_\la)$ can
be written in the form
(\ref{10.eta}) with $(F,C_\mu)\in {\cal A}^{\de+c\la}_\la$ and
$\|\eta\|_0\le \ep$; this
$\eta$ must then  be the unique fixed point $\eta(f,\mu)$ of
Proposition \ref{exist.sol}.  Thus $\Phi^1_\la$ is surjective.
\qed

\vskip.3in

Diagram \ref{10.diagram1} leads to our first formula expressing the
absolute invariants of a symplectic sum $Z= Z_\la$ in terms of the
relative invariants of $X$ and $Y$.  Recall that the relative
invariant $GW^V_X$ is obtained by forming the space $
\ov{\M}_{\chi,n,s}^V(X,A)$ of relatively stable maps and pushing
forward its fundamental homology class by the map
\bear
\ep_V:   \ov{\M}_{\chi,n,s}^V(X,A) \to \wt{\M}_{\chi,n}  \ti X^n \ti
{\cal H}_{X,A,s}^V.
\label{10.MHVSV}
\eear
We can also consider the space of stable maps from compact, not
necessarily connected domains by
taking  the union of products of $\ov{\M}_{\chi,n,s}^V(X,A)$ and
again pushing forward in
homology. The resulting class in the homology of $\wt{\M}_{\chi,n}
\ti X^n \ti {\cal
H}_{X,A,s}^V$ is the relative TW invariant (\ref{2.relTW}).  As we
observed in the
introduction (see Figure 1), it is the TW invariant that will appear
in the symplectic sum
formula.

\bigskip

To proceed, then we should replace the vertical arrows in Diagram
\ref{10.diagram1} by the
above maps $\ep_V$ and  pass to homology.  We will do that in two
steps, first incorporating
the   spaces ${\cal H}_{X}^V$ and then including the $X^n$.  In each
case we will see
that the operation of gluing maps defines an extension of the
bottom arrow in Diagram \ref{10.diagram1}, which we examine in homology.

\bigskip

\noindent {\bf The Convolution Operation}\ \  We can glue a map $f_1$
into $X$ to a map $f_2$
into $Y$ provided the images  meet $V$ at the same points with the same
multiplicity.  The  domains of $f_1$ and $f_2$ glue according to the
attaching map $\xi$ of
(\ref{def.xi}), while the images determine elements of the
intersection-homology spaces ${\cal
H}_{X,A,s}^V$ and ${\cal H}_{Y,A,s}^V$ which glue according to the
map $g$ of (\ref{def.g}).
The convolution operation records the effect of these gluings at the
level of homology.

\medskip

For each $s$ the attaching map (\ref{def.xi}) defines a  bilinear form
$$
\left(\xi_{\ell}\right)_*:\ H_*(\wt M;\Q)
\otimes H_*(\wt M;\Q)
\longrightarrow H_*(\wt M;\Q)
$$
  for $\ell=\ell(s)$.  Similarly, for each $s$  the map
$g$ from (\ref{def.g}) induces a bilinear form on the
homology of ${\cal H}_Y^{V}\ti {\cal H}_Y^{V}$ with values in $RH_2(Z)$,
the (rational) group ring of $H_2(Z)$,  namely
\best
&&\langle\ ,\ \rangle:\ H_*({\cal H}_X^V;\Q)
\otimes H_*({\cal H}_Y^{V};\Q)
\longrightarrow  RH_2(Z)\\ \\
&&\langle  h\; , \, h'  \rangle_s=g_*\left[
      \left. h\ti h'\right|_{\ep^{-1}\left({\bold \Delta}_{s}\right)}\right]=
      \sum_{A\in H_2(Z)}\ g_*[{\bold \Delta}_{A,s}\cap (h\ti h')]\ t_A .
\eest
  This last equality holds because
$\ep^{-1}\left({\bold \Delta}_s\right)$ is the union of components ${\bold
\Delta}_{A,s}=\ep^{-1}\left({\bold \Delta}_s\right) \cap
g^{-1}(A)$.

\smallskip

Combining the two bilinear forms gives
the convolution operator that describes how homology
classes of maps combine in the gluing operation.

\begin{defn}\ \ The convolution operator
$$
*: \ H_*(\wt{\M}\ti {\cal H}_X^{V};\,\Q[\la])
\otimes H_*(\wt{\M}\ti {\cal H}_Y^{V};\,\Q[\la])
\ \longrightarrow\  H_*(\wt{\M};\, RH_2(Z)[\la])
$$
is given  by
\begin{equation}
(\kappa\otimes h)\,*\,(\kappa'\otimes h') \ =\
\sum_{s} {|s| \over \ell(s)!}\; \la^{2\ell(s)}\;
\left(\xi_{\ell(s)}\right)_* (\kappa \otimes\kappa')\ \  \langle  h \; , \,h'
\rangle_s
\label{wteddiagonal1}
\end{equation}
\end{defn}

\medskip

The right hand side of (\ref{wteddiagonal1}) includes three numerical 
factors which keep track of
how maps glue when we form the symplectic sum.  Recall that the 
powers of $\la$ record the
euler characteristic in the generating series of the invariants 
(\ref{defn.7}) and
(\ref{defnRelInvt2}); the factor $\la^{2\ell(s)}$ in 
(\ref{wteddiagonal1}) the reflects the
relation (\ref{2.sumofchis}) between the euler characteristics when we 
glue along $\ell(s)$
points.  The factor $|s|$ is the degree of the covering in Theorem 
\ref{gluing thm}; this
reflects the fact that each stable map into $Z_0$ can be smoothed in 
$|s|=s_1\cdots s_\ell$
ways.  Finally,  note that elements in the space $\Mf_s(Z_\la)$ in Diagram
\ref{10.diagram1} are {\em labeled} maps, i.e. they have $\ell(s)$ 
numbered curves on their domains as explained at the end of section
3.  But the GW and TW invariants of $Z_\la$
are defined using the space of {\em unlabeled} stable maps, which is 
the quotient of the space of labeled maps by the action of the 
symmetric group.  That accounts for the factor in $1/\ell(s)!$ is 
(\ref{wteddiagonal1}).

\medskip

Since ${\cal H}^V_X$ is the disjoint union of components ${\cal
H}^V_{X, A,s}$ with $A\in H_2(X)$ and $\deg s=A\cdot V$, there is an 
isomorphism
\best
H_*({\cal H}^V_X)\ \cong\  \sum_{A}\sum_{\deg s=A\cdot V} H_*({\cal 
H}^V_{X,A,s})\ t_{A}.
\eest
Below,  we will identify $h\in H_*({\cal H}^V_X)$ with $\sum_A h_A
t_A$, where $h_A$ are its components in $H_*({\cal H}^V_{X,A})$.

\bigskip

\begin{ex}
\label{conv.coord}
The formula for the convolution  simplifies  when there are no rim 
tori in $X$ and $Y$,
and therefore in $Z$  (c.f. (\ref{def.Rim})). Then  (i) the relative 
invariants have an expansion
of the form (\ref{2.last}), (ii) the map $g$ of
(\ref{def.g}) is the restriction to the diagonal $\Delta_s\subset
V^s\ti V^s$, and (iii) the $h$ part of the convolution
(\ref{wteddiagonal1}) is then given by the cap product with the
Poincar\'{e} dual of the diagonal:
$$
g_*\left[ \left.h\ti h'\right|_{\bold {\Delta}_{s}}\right]\ =\;
{\mbox{PD}}\left(\bold {\Delta}_{s}\right)\cap (h\ti h').
$$
We can then `split the diagonal' by fixing a basis $\{C^{p}\}$ of
$H^*(\bigsqcup_s V^s)$ and writing
$$
{\mbox{PD}}\left(\bold {\Delta}_{s}\right)\ =\
\sum_{p,q} Q^V_{p,q}\ C^p\ti C^q= \
\sum_{p} \ C^p\ti C_p
$$
where $ Q^V_{p,q}$ is the intersection form of $V^s$ for the  basis
$\{C^p\}$ and $C_p=\sum  Q^V_{p,q}\ C^q$ is the dual basis. If
$\{\gamma^{i}\}$ is a basis of $H_*(V)$, let  $\{{\bf C}_m\}$ be
the basis (\ref{last.basis2}) of $H^*(\bigsqcup_s V^s)$ corresponding to
$\{\gamma^i\}$ and let  $\{{\bf C}_{m^*}\}$ be the one corresponding
to the dual basis $\{\gamma_i\}$ (with respect to $Q^V$). The
convolution then has the more explicit form
\bear
(\kappa\otimes h)\,*\,(\kappa'\otimes h') \ =\
\sum_{m} {|m|\over m!}\;\la^{2\ell(m)}\;
\left(\xi_{\ell(m)}\right)_*(\kappa\otimes \kappa')  \
{\bf C}^*_{m}(h)\,  {\bf C}^*_{m^*}(h').
\label{*.m}
\eear
In passing from $s$ to $m$, we
used the fact that each fixed sequence $m$
corresponds to $\displaystyle
{{\ell(s)} \choose {(m_{a,i})}}= {\ell(s)!\over m!}$
{\em ordered} sequences $s$.
\label{ex.10.3}
\end{ex}

\bigskip

More generally, let $X$ be a symplectic manifold with two disjoint
symplectic  submanifolds $U$ and $V$ with real
codimension two.  Suppose that $V$ is symplectically identified with a
submanifold of
similar triple  $(Y,V,W)$ and that the normal bundles of
$V\subset X$ and
$V\subset Y$ have opposite chern
classes.  Let $(Z,U,W)$ be the resulting symplectic  sum. In this case,
(\ref{def.g}) is replaced by
\bear
g:{\cal H}_X^{U,V}\ti_{\ep} {\cal H}_Y^{V,W}\ra {\cal H}_Z^{U,W}
\label{gmap2}
\eear
which combines with the map $\xi_{\ell(s)}$ to give the convolution
operator
\bear
*: \ H_*(\wt{\M}\ti {\cal H}_X^{U,V};\Q[\la])
\otimes H_*(\wt{\M}\ti {\cal H}_Y^{V,W};\Q[\la])
\ \longrightarrow\  H_*(\wt{\M}\ti {\cal H}_Z^{U,W};\Q[\la])
\label{10.convolution}
\eear
as in (\ref{wteddiagonal1}). It  describes how homology
classes of maps combine in the gluing operation for the symplectic sum.

\bigskip

Finally, we include the   evaluation maps which record the images of
the  $n$ marked
points. These combine with
the projections from (\ref{2.collapsingmaps}) to give the diagram
\begin{equation}
\begin{array}{ccc}
\ma\bigsqcup_s \,{\cal AM}_s & \longrightarrow &
\ma\bigsqcup_s \,\Mf_s(Z_\la)\\
\Big\downarrow\vcenter{\rlap{{\scriptsize \mbox{$ev$}}}} &  & \Big\downarrow
\vcenter{\rlap{{\scriptsize \mbox{$ev$}}}}  \\
\ma\bigsqcup_n \;(X\sqcup Y)^n&& \ma\bigsqcup_n \;(Z_\la)^n
\\
{\scriptsize \mbox{$\pi_0$}}\searrow&&\swarrow{\scriptsize
\mbox{$\pi_\la$}}\quad
\\
&\ma\bigsqcup_n \;(Z_0)^n&
\end{array}
\label{big.diagr}
\end{equation}
which commutes up to homotopy. We can also include the spaces
$\wt{\M}$ of curves from Diagram \ref{10.diagram1}.  Pushing
forward then gives $\pi_{0*} (TW_X^V*TW_Y^V)=\pi_{\la*}( TW(Z_\la))$. 

\medskip

\begin{theorem} Assume that all curves
contributing to the invariants are flat along $V$. Then  (\ref{3.1}) holds in
the sense that  for any $\al_0\in {\Bbb T}(H^*(Z_0))$
\bear
TW_Z^{U\cup W}(\pi^*\al_0)\, = \,\l(TW_X^{U\cup V}\, *\, TW_Y^{ V\cup
W}\r)
(\pi_0^*\al_0).
\label{mainformflatcase}
\eear
\label{mainthmflatcase}
\end{theorem}
\pf  It suffices to verify this for decomposable elements  $\al_0 
=\al_0^1 \otimes \dots \otimes
\al_0^n$. Let $\al_V^k$, $\al_X^k$, $\al_Y^k$
denote the restriction of $\al_0^k$ to
$V$, $X$ and respectively $Y$. We can then choose geometric
representatives $B_V^k$ of the Poincar\'{e} dual of $\al_V^k$ in $V$ and
Poincar\'{e} duals $B_X^k$ of $\al_X^k$ in $X$ and $B_Y^k$ of $\al_Y^k$ in $Y$
which intersect $V$ transversely such that moreover $B_X^k\cap V=B_Y^k\cap
V=B_V^k$. Then the inverse image under $\pi_\la$
of $B_X^k\ma\cup_{B_V^k}B_Y^k$ gives a continuous family of geometric
representatives $B_\la^k$ of the Poincar\'{e} dual of $\pi^*\al_0^k$ in
$H^*(Z_\la)$.  The theorem then follows from Theorem \ref{gluing thm}
by cutting down the moduli spaces on the left of Diagram \ref{big.diagr} by 
$(B_X, B_Y)$ and the ones on the right by $B_\la$. Constraints in
$H^*(\wt\M)$ are handled similarly. The details of such arguments are standard 
(c.f.  \cite{rt1}).
\qed

\bigskip

We should comment on how the assumption that all maps are $\de$-flat
enters the above proof.
Notice that in the statement of Theorem \ref{gluing thm} the
$\de$-flat  maps  in ${\cal AM}_s$
are  paired with maps in $\M_s(Z_\la)$ which are not exactly
$\de$-flat --- there is a slight
variation in $\delta$.    But when all contributing maps  are flat,
the cut-down moduli
space $\ev^{-1}(B_\la)\subset \ov \M(Z_\la)$ limits as $\la\to 0$
to a compact subset of the
open set $\M_s\ti_{ev}\M_s$ as in (\ref{def.K}).  Hence for 
sufficiently small $\de$ the
set of elements of the  limit
set which are $\de$-flat is the same as the set of
$2\de$-flat elements, so the variation in $\de$ is inconsequential.

\bigskip

Theorem \ref{mainthmflatcase} is  a formula for the TW invariants
evaluated on only certain
constraints in $H^*(Z_\la$) --- those of the form  $\pi^*(\al_0)$.
The following definition
characterizes those constraints.   It  is based on the diagram
induced by the collapsing maps of
(\ref{2.collapsingmaps})

\smallskip
\begin{equation}
\begin{CD}
@.@.@.{\Bbb T}(H^*(Z_0))@.
\\
@.@.{\pi^*}\swarrow@.@.\searrow{\pi_0^*}
\\
@.{\Bbb T}(H^*(Z))@.@.@. @. {\Bbb T}(H^*(X)\oplus H^*(Y))
\end{CD}
\label{al.separates}
\end{equation}

\smallskip

\begin{defn} We say that a constraint  $\al\in {\Bbb T}(H^*(Z))$ separates
as
$(\al_X, \al_Y)$ if there
exists an $\al_0\in {\Bbb T}(H^*(Z_0))$ so that
$\pi^*\al_0=\al$ and $\pi_0^*(\al_0)=(\al_X,\al_Y)\in 
{\Bbb T}(H^*(X)\oplus H^*(Y))$.
\label{def.al.separates}
\end{defn}

\medskip

Here are three observations to help clarify which classes $\al\in
H^*(Z)$ separate.  These follow by combining  the Mayer-Vietoris
sequences for $Z_\la=(X\setminus V)\cup (Y\setminus V)$
\[\begin{CD}
H^{*-1}(S_V)@>{\de^*}>>H^*(Z)@>{i^*}>>
H^*(X\setminus V)\oplus H^*(Y\setminus V)@>{}>>H^*(S_V)@>{\de^*}>>
\end{CD}
\label{10.MV}
\]
and the similar one for $Z_0$ with the Gysin sequence for $p:S_V\ra V$
\bear
\begin{CD}
H^{*-2}(V)@>{\cup c_1}>>H^{*}(V)@>{p^*}>>H^{*}(S_V)@>{p_*}>>
H^{*-1}(V).
\end{CD}
\label{10.gysin}
\eear

\begin{enumerate}
\item[(a)]  When the  first map in (\ref{10.gysin}) is injective
then all classes $\al$ separate. In  dimension four, that occurs
       whenever the normal
bundle of $V$ in $X$ is topologically non-trivial.
\item[(b)] In general the separating classes are those $\al$ for which
$j^*(\al)\in H^{*}(S_V)$ is in the image of the second map in
(\ref{10.gysin}).
\item[(c)] the decomposition $(\al_X,\al_Y)$, if it exists, is unique only
up to elements in the image of $\de_X^*\oplus \de^*_Y:H^{*-1}(S_V)\to
H^*(X)\oplus H^*(Y)$ (the elements that can be ``pushed to either
side'').
\end{enumerate}

\medskip

Using Definition \ref{def.al.separates} and for simplicity taking $U$
and $W$   to be  empty,
Theorem
\ref{mainthmflatcase} becomes:
\begin{theorem}
\label{mainthoeremflatcase}  Suppose that all curves
contributing to the invariants are flat along $V$ and $\al$ separates
as $(\al_X,\al_Y)$.  Then
\bear
TW_Z(\al)\, = \,\l(TW_X^{V}\, *\, TW_Y^{ V}\r) (\alpha_X, \alpha_Y).
\label{mainformulaflatcase}
\eear
\end{theorem}

\medskip

Note that when $(\al_X,\al_Y)$ decomposes as $\al=\al_X\otimes \al_Y$ the
right hand side is
$TW_X^{V}(\alpha_X)\, *\, TW_Y^{ V}(\alpha_Y)$, but in general
$(\al_X,\al_Y)$ is a sum of tensors of the form  $(\al_X^1 +
\al_Y^1)\otimes \cdots\otimes
(\al_X^k + \al_Y^k)$ and the right hand side of
(\ref{mainformulaflatcase}) is the corresponding
sum.

\medskip

To focus on the decomposable case we make another definition:  we say
$\al$ is {\em supported
off the neck} if the restriction $j^*(\al)\in H^*(S_V)$ vanishes.  In
that case $\al$
separates into relative classes $\al_X\in H^*(X,V)$ and $\al_Y\in
H^*(Y,V)$, generally in several ways.  For each such decomposition Theorem
\ref{mainthoeremflatcase} gives
\bear
TW_Z(\al_X,\al_Y)\, = \,TW_X^{V}\,(\al_X) *\, TW_Y^{V}
(\alpha_Y).
\label{10.lastfrom}
\eear
This was the formula described in \cite{ipa}.

\bigskip

\begin{ex}  Take $\al$ to be  the Poincare dual of a point in $Z$.
This constraint is
supported off the neck and has two independent decompositions
depending whether the point is in $X$ or $Y$.
\end{ex}

\bigskip

\begin{ex}  Suppose $\al=\al_X\otimes \al_Y$ is supported off the
neck and there are no rim tori in $(X,V)$ and $(Y,V)$ and that all
curves contributing to the invariants are $V$-flat.  Then we can
choose a basis of $H^*(V)$ and expand the relative TW
invariants as in Example \ref{ex.10.3}. Combining (\ref{10.lastfrom})
with  (\ref{*.m}) gives the explicit   formula
\best
TW_{\chi,A,Z}(\al_X,\al_Y)=\sum_{A=A_1+A_2\atop
\chi_1+\chi_2-2\ell(m)=\chi}\sum_{m} \la^{2\ell(m)}{|m|\over m!}\
TW^V_{\chi_1,A_1,X}\left(\alpha_X;C_m\right)
\cdot TW^V_{\chi_2,A_2,Y}\left(C_{m^*} ;\alpha_Y\right).
\eest
Note that from the definition of relative invariants, the only terms
contributing are those for which $A_1\cdot V=\ell(m)= A_2\cdot V$.
E. Getzler has pointed out that the formula above can be neatly expressed
in terms of the generating
series (\ref{TW.gen.fcn}) and the intersection matrix $Q^V$ of $V$,
specifically
\best
TW_{Z}(\al_X,\al_Y)&=&\l.\exp \l( \sum_{a,i,j} a \la^{2}\; Q^V_{ij}\;
{\partial\over\partial z_{a,i}}\;{\partial\over\partial w_{a,j}} \r)
\l(TW^V_X(\al_X)(z)\cdot TW^V_Y(\al_Y)(w)\r)\r|_{z=w=0} .
\eest
\label{al.off.neck}
\end{ex}

\medskip

Because the decomposition of separating constraints $\al$ is not
unique, we can often
choose several different decompositions, and use Theorem
\ref{mainformulaflatcase} to get several expressions for the same TW
invariant.  That
yields relations among relative TW invariants.  In  Section
\ref{section15} we will use that
idea to derive recursive formulas which determine the relative
invariants in some interesting
cases.

\vskip.4in


\setcounter{equation}{0}
\section{The space ${\Bbb F}$ and the S-matrix}
\bigskip

Starting from the normal bundle $N_XV$
of $V$ in $X$, we can form the $\P^1$ bundle
$$
{\Bbb F}\ =\ {\Bbb F}_V\ =\  \P(N_X V\oplus \cx)
$$
over $V$ by projectivizing the sum of the normal bundle $N_X V$ and
the trivial complex line bundle. Let $\pi:{\Bbb F}\to V$ be the
projection map. In ${\Bbb F}$, the zero section $V_0$ and the infinity
section $V_\infty$ are disjoint symplectic submanifolds, both
symplectomorphic to $V$. Moreover, note that
${\Bbb F}\#_{V}{\Bbb F}={\Bbb F}$.
\medskip

Under the natural identification of $V_0$ with $V_\infty$, the
convolution operation
(\ref{10.convolution}) defines an algebra structure on
$H_*(\wt{\M}\ti {\cal H}_{\Bbb
F}^{V,V};\Q[\la])$.  That allows us to multiply by TW invariants.  Of
particular interest are the invariants with no constraints on the
image, that is $TW^{V,V}_{\Bbb F}(\al)$ with $\al=1$, which give an
operator
\bear
\l[TW^{V,V}_{\Bbb F}(1)\r]\,*\,: \  H_*(\ov{\M}\ti{\cal H}^{V}_{\Bbb
F};\,\Q[\la])\ \to\ H_*(\ov{\M}\ti{\cal H}^{V}_{\Bbb F};\,\Q[\la])
\label{11.defTWVV1}
\eear
defined by  a power series as in (\ref{defnRelInvt2}).
This operator is key to the general symplectic
sum formula given in the next section.  In this section we describe
(\ref{11.defTWVV1}) and its inverse
and develop  some examples.

\medskip

      Each  $(J,\nu)$-holomorphic bubble map $f$ into ${\Bbb F}$  projects to a
map $f_V=\pi\circ f$ into $V$.   Although $f_V$ may not be
$(J,\nu)$-holomorphic, we
can still ask whether $f_V$ is stable, using the second definition of
stability given after (\ref{2.defstable}),
namely $f$ is stable if its restriction to each unstable domain component
is non-trivial in homology.

\begin{defn}
A $(V_0,V_\infty)$-stable map $f:C\to {\Bbb F}$  is 
{\em ${\Bbb F}$-trivial} if each of its components is an
unstable rational curve whose image represents a multiple of the fiber $F$
of ${\Bbb F}$.
\label{5.relativelystable}
\end{defn}

Thus the ${\Bbb F}$-trivial curves are rational curves representing $dF$
with one marked point on the zero section
and one on the infinity section, both intersecting with multiplicity $d$.
Let $\M_{\Bbb I}$ denote the set of
${\Bbb F}$-trivial maps in
$\M^{V_0,V_\infty}_{\Bbb F}$ and consider the disjoint union
\bear
\M^{V_0,V_\infty}_{\Bbb F}\ =\ \M_{\Bbb I}\ \cup\ \M_{R}
\label{MI+MRdecomp}
\eear
where $M_{R}$ is the set of non-${\Bbb F}$-trivial maps.

\bigskip

For the next lemma we fix a metric $g'$ on ${\Bbb F}$ for which $\pi:{\Bbb
F}\to V$ is a Riemannian submersion.   The
procedure described in the appendix of \cite{ip4} then constructs a
compatible triple $(\w,J,g)$ on ${\Bbb F}$
for which $\pi$ is holomorphic and is a Riemannian submersion.  Using this
metric, each  perturbation term $\nu_V$
on $V$  has a horizontal lift  $\pi^*\nu_V$ in $\Omega^{0,1}(T{\Bbb F})$.
We will call such a structure $(\w,J,g,
\pi^*\nu_V)$ a {\it submersive structure}.  For submersive structures,
each
$(J,\pi^*\nu_V)$-holomorphic map $(f,j)$ into ${\Bbb F}$ projects to a
$(J,\pi^*\nu_V)$-holomorphic map $(\pi\circ f,j)$ into $V$.

\begin{lemma} (a) $\M_{\Bbb I}$ is both open and closed.  The
corresponding
decomposition of (\ref{11.defTWVV1}) is
\bear
TW_{\Bbb F}^{V,V}(1)={\Bbb I}+ R^{V,V}
\label{TW=1+R}
\eear
that is, the ${\Bbb F}$-trivial maps contribute the identity to the TW
invariant.

\smallskip

(b) The non-${\Bbb F}$-trivial maps have $E(f_V)\geq \alpha_V$, where
$\alpha_V$ is the constant of Definition \ref {defnFlat}.

\smallskip

(c)  For each fixed $A$, $n$ and $\chi$, the corresponding term in
the convolution $R^m=R*\cdots *R$ vanishes for $m$ large enough.
Therefore, the inverse of $TW$ is well defined by:
\bear
\label{TWinverse}
\left(TW_{\Bbb F}^{V,V}(1)\right)^{-1}\ =\ \sum_{m=0}^{\infty} (-1)^m
R^m.
\eear
\label{TW=1+Rlemma}
\end{lemma}
\pf (a)  Clearly  $\M_{\Bbb I}$ is closed.   To show that the complement
of
$\M_{\Bbb I}$ is closed, suppose that a sequence $(f_i)$
      in the complement converges to a trivial map $f$ in the topology of the
space of stable maps.  Then the homology
classes converge so, after passing to a subsequence, we can assume that
each $f_i$ represents $dF$.  Similarly,
the stabilizations of the domains converge in the Deligne-Mumford space,
so
we can assume that all domain
components of each $f_i$ are unstable.  But then the $f_i$  lie in
$\M_{\Bbb I}$.   We conclude that $\M_{\Bbb I}$ is both open and
closed. Finally, the
decomposition (\ref{MI+MRdecomp})
gives splitting   (\ref{TW=1+R}) of the TW invariant because  convolution
by
elements of $\M_{\Bbb I}$
is the identity.
\smallskip

(b)  If $E(f_V)< \alpha_V$ then, as in the proof of Lemma 1.5 of
\cite{ip4}, every component of the domain is
unstable and $f_V$ is trivial in homology and therefore $f$ represents a
multiple of $F$.

\smallskip

(c) For each $(J, \nu)$ , we shall bound the number $N$ for which
there are maps   in the moduli
space defining the convolution $R^N$.  That moduli
space consists of maps $f$ from a domain $C$ (whose Euler class $\chi$ and
number $n$ of marked points is fixed)
to the singular manifold
${\Bbb F}\#\cdots \# {\Bbb F}$ obtained from $N$ copies of ${\Bbb F}$ by
identifying the infinity section of
one with  the zero section of the next.  Furthermore, these $f$ decompose
as $f=\bigcup f^j$ where $f^j$ is a map
from some of the components of $C$ into the $j^{\mbox{th}}$ copy of
      ${\Bbb F}$.

Fixing such an $f$, let $N_1$ be the number of $f^j$
whose domain has at least one stable component $C_j$.  These components
appear in the stabilization $st(C)$.  But
$st(C)$ lies in the space $\M_{\chi,n}$ of stable curves, and hence has at
most
$\mbox{dim}\ \M_{\chi,n}$ components.  This gives an explicit bound for
$N_1$ in terms of $\chi$ and $n$.

The remaining $N_2=N-N_1$ of the $f^j$ each have a domain component with
$\pi_*[f^j(C_j)]\in H_2(V)$ non-trivial, so
satisfy $E(\pi\circ f^j)>\alpha_V$ by (b) above. We therefore have
$$
N_2\alpha_V \ \leq \  \sum E(\pi\circ f^j) \  \leq \ E(\pi\circ f)  \
\le \ C\left[A(\pi(f))+C_\nu\right],
$$
where the first sum is over those $j$ contributing to $N_2$ and the last
inequality is as in the proof of Lemma \ref{5.energybound}. Since the
symplectic area $A(\pi(f))$ of the projection is a topological
quantity, this bounds $N_2$ and hence $N$. \qed

\bigskip

\begin{defn} The $S$-matrix is defined to be the inverse of the TW
invariant of  Lemma
\ref{TW=1+Rlemma}:
\best
S_V=\left(TW_{\Bbb F}^{V,V}(1)\right)^{-1}.
\eest
(Note that this depends not just on $V$ but on $N_V$ and the 1-jet of
$(J,\nu)$ along $V$.)
\label{defSmatrix}
\end{defn}

\medskip
The symplectic  sum of $(X,U,V)$ and $({\Bbb F}, V_\infty,V_0)$  along
$V=V_\infty$ is a
symplectic deformation of $(X,U,V)$,
so has the same  $TW$ invariant. The convolution then defines a operation
$$
      \ H_*(\wt{\M}\ti {\cal H}_X^{U,V};\Q[\la])
\otimes H_*(\wt{\M}\ti {\cal H}_{{\Bbb F}}^{V,V};\Q[\la])
\ \longrightarrow\  H_*(\wt{\M}\ti {\cal H}_X^{U,V};\Q[\la]).
$$
Thus for each choice of constraints $\alpha\in {\Bbb T}({\Bbb
F},V_\infty\cup
V_0)$, the
$TW$ invariant of ${\Bbb F}$ relative to its zero and infinity section
defines an endomorphism
\bear
\label{3.Endo}
TW_{{\Bbb F}}^{V_\infty,V_0}(\alpha)\,\in\, \mbox{End}\
\left(H_*(\wt{\M}\ti {\cal H}_X^{U,V};\Q[\la])  \right)
\eear
which describes how families of curves on $X$
are modified --- ``scattered''--- as they pass through a neck modeled on
$({\Bbb F},V_\infty,V_0)$  containing the constraints $\alpha$.

\medskip

The identity endomorphism in  (\ref{3.Endo}) is always
realized
as  the convolution
by the element
\best
{\Bbb I}\in H_*(\wt{\M}\ti {\cal H}_{{\Bbb F}}^{V,V};\Q[\la])
\eest
corresponding to that part of $TW$ coming from ${\Bbb F}$-trivial maps.
Thus the statement that
$S_V=\mbox{Id.}$ means that the only curves present are those which are
irreducible fibers of ${\Bbb F}$.

\vskip.3in

\begin{ex} When $V=\P^1$,   ${\Bbb F}\to V$ is one of the rational ruled
surfaces with its standard symplectic structure.  If we wish to count all
pseudo-holomorphic maps, without constraints on the genus
or the induced complex structure, the relevant S-matrix is  the relative
$TW$ invariant  with
$(\kappa,\alpha)= (1,1)$.  This  case
works out neatly:  Lemma \ref{CHproposgator0} implies that
$S_V=\mbox {Id.}$
\label{ex11.1}
\end{ex}

\vskip.3in

\begin{ex}  When we put no constraints on either the domain or the
image $S_V$  is an
operator given in terms of  $TW^{V,V}_{{\Bbb F}}$ by the $S$-matrix expansion
(\ref{TW=1+Rlemma}).  In  cases where there are no rim tori in ${\Bbb
F}$, we can expand the TW
invariants in the power series (\ref{TW.gen.fcn}) of the appendix.  Letting
$TW_{\chi, A}(C_{m};C_{m'})$
denote the relative invariant of ${\Bbb F}$ satisfying the contact
constraints $ C_{m}$ along
$V_\infty$ and $C_{m'}$ along $V_0$,  the $S$-matrix expansion shows
that $S_V$ has an expansion like (\ref{TW.gen.fcn}) with coefficients
\best
S_{\chi, A}(C_m;\; C_{m'})
      & = & \de_{m,m'}\ -\ TW^{V,V}_{{\Bbb F},\chi, A}(C_{m};\
C_{m'}) \\
& + &\hskip-.2in
       \sum_{A_1+A_2=A\atop\chi_1+\chi_2-2\ell(s_1)=\chi} \ma\sum_{m_1}
\la^{2\ell(m_1)}{|m_1|\over m_1!}\
TW^{V,V}_{{\Bbb F},\chi_1, A_1}(C_{m}; C_{m_1})\,
TW^{V,V}_{{\Bbb F},\chi_2, A_2}(C_{m_1^*};\  C_{m'} )-\dots
\eest
\label{ex11.2}
\end{ex}

\vskip.2in

\medskip
\setcounter{equation}{0}
\section{The General Sum Formula}
\label{The General Sum Formula}
\bigskip

In all of our work thus far we have assumed that the
$(J,\nu)$-holomorphic maps we are gluing are
$\delta$-flat as in Definition \ref{defnFlat}. In this section we
remove this flatness assumption and
prove the symplectic sum formula in the general case.

The idea is to reduce the general case to the flat case by  degenerating
       along many parallel copies of $V$.   Thus instead of viewing
$Z_\la$ as the
symplectic sum  $X\ \#_VY$  along
$V$ we regard it  as the symplectic sum  of  $2N+2$ spaces:  $X$ and
$Y$ at the ends and $2N$ middle
pieces each of which is a copy of the ruled space ${\Bbb F}$ associated to
$V$ --- see Figure 2 of the introduction.  The  pigeon-hole principle 
then implies that for large $N$ all
holomorphic maps into
$Z_\la$ are close   to maps which are flat along each `seam' of the
$2N$-fold  sum.

\begin{lemma}
There is a constant $E=E_{\chi,n,A}(J,\nu)$ such that every
$(J,\nu)$-holomorphic map into $Z$
representing a class $A\in H_2(Z)$ has energy at most $E$.
\label{5.energybound}
\end{lemma}
\pf In an orthonormal frame $\{e_1,e_2=je_1\}$ on the domain, the
holomorphic map equation is $f_*e_1+Jf_*e_2=2\nu(e_1)$. Taking the
norm squared and noting that
$\langle f_*e_1,Jf_*e_2\rangle=f^*\w (e_1,e_2)$ gives $|df|^2=2|\nu|^2+2f^*\w
(e_1,e_2)$. The energy is therefore the $L^2$ norm of $\nu$ plus the
topological quantity
$\langle\w,A\rangle$. The lemma follows.
\qed

\medskip

For the remainder of this section we fix the data $\chi,n,A,J,\nu$ which
determined the constant $E$ of
Lemma \ref{5.energybound} and fix an integer $N$ with
\bear
N\alpha_V\ >\ E
\label{5.defN}
\eear
where  $\alpha_V<1$ is the constant of Definition \ref {defnFlat}.

\medskip

Fixing $\la$, we partition the neck
of $Z=Z_\la$ into $2N$ segments $Z^j$ using the coordinate $t$ from
(\ref{3.momentmap}):
$$
Z^j\ =\ \left\{z\in Z_\la\ |\ (j-N-1)\ep\leq  t(z)\leq \ (j-N)\ep\
\right\} \qquad j=1,
\dots, 2N
$$
where $\ep$ is as in Figure 3. Squeezing the neck at the
midpoints $t_j(z)=j-N-\frac12$ of each of these
segments
defines a family
\bear
{\cal Z} \to D\subset \cx^{2N+1}
\label{5.bigfamily}
\eear
as in Theorem \ref{thm1.1} but with many `necks'. Thus the fiber over
$(\mu_1, \dots,\mu_{2N+1})$, defined for $|\mu|<<|\la|$,  is a space
$Z_\la(\mu_1, \dots,\mu_{2N+1})$ with
a neck of size $\mu_j$ inside each
$Z^j$ and the fiber over
$\mu=0$ is the singular space obtained by connecting $X$ to $Y$ through a
series of $2N$ copies of the rational ruled manifold ${\Bbb F}$
associated with $V$.  One such space is depicted in Figure 2 of the 
Introduction.

\medskip

Fix $\de>0$ such that $\de\le {\ep\over 10 N}$ and consider the space $\M=
\M_{\chi,n,A}(Z_{\la})$ of holomorphic maps into $Z_\la$.  Let $f^j$
denote the restriction of $f\in\M$ to $f^{-1}(Z^j)$.  We can then
define an open cover of $\M$ that keeps track of the values of $j$ for
which the energy $E_\de(f^j)$ on the $\de$ neck around the cut is small as in
equation (\ref{defnFlateq}). Specifically, to each subset $\{
i_1,\dots, i_k\}$ of $\{1,\dots, 2N\}$ we associate the open subset of
$\M$
\bear
\M^{i_1,\dots i_k}\ =\ \left\{f\in \M\ |\ E_\de(f^j)<\alpha_V/2\mbox{ for }
j=i_1,\dots, i_k\; \right\}.
\label{5.defMVjm}
\eear
\medskip

\begin{lemma}
The $\M^{i_1,\dots i_k}$ cover $\M=\M_{\chi,n,A}(Z_{\la},A)$ and
set theoretically
\bear
\M\ =\ \bigcup \M^{i}\ -\  \bigcup\M^{i_1,i_2}\ +\  \bigcup
\M^{i_1,i_2,i_3}\ -\ \dots
\label{5.alternatingset}
\eear
\label{5.alternatingsetlemma}
\end{lemma}

\pf  Each $f\in\M$ has $\sum_jE(f^j)\leq E(f)< E$, so
(\ref{5.defN}) implies that $f\in \M^{i}$ for at least one $i$. If
$E_\de(f^j)<\alpha_V/2$ for exactly  $\ell$ of the
$j$, then  $f$ is counted
$$
\ell-{\ell\choose 2}+{\ell\choose 3}-\cdots \pm {\ell\choose \ell}\ =\ 1
$$
times on the right hand side of (\ref{5.alternatingset}).
\qed

\bigskip

Now every $f\in\M^{i_1,\dots i_k}$ has small energy in the segment $Z^j$
for $j=i_1,\dots, i_k$.    Replacing these
$\la_j$ by $\mu_j=\mu\la_j$ for those values of $j$ (and keeping the
remaining $\la_j$ fixed) defines a
1-parameter subfamily $Z_\mu$ of (\ref{5.bigfamily}).  That family
degenerates in the middle of exactly
$k$ of the segments $Z^j$.  At each of those degenerations  $f$  is
$\delta$-flat in the sense of Definition \ref{defnFlat}. Hence
\bear
\M^{i_1,\dots i_k}=\M_X^V \times_{ev} (\M_{\Bbb
F}^{V,V})^{k-1}\times_{ev} \M_Y^V
\label{5.piece}
\eear
We can therefore apply the  sum formula
(\ref{mainthoeremflatcase}), obtaining, for a
fixed $A$ and $\chi$,
\bear
TW_{X\# Y}= TW_X^V *\left[ \ma\sum_{k=1}^{2N} (-1)^{k-1}{2N\choose k}
(TW_{\Bbb F}^{V,V})^{k-1}\right] * TW_Y^V.
\label{eq.5.3}
\eear
This formula appears to be dependent on the number of cuts $2N$.
However, there is
a way to rewrite it to see
that  it is independent of $N$.   Note that after multiplying by $TW$ the
middle sum is a binomial expansion, in
fact, using Lemma \ref{TW=1+Rlemma}c,
$$
\ma\sum_{k=1}^{2N} (-1)^{k-1}{2N\choose k} (TW)^{k-1}\ =\
\frac{1-(1-TW)^{2N}}{T}\ =\ \frac{1-(-R)^{2N}}{TW}\ =\ TW^{-1}.
$$
Thus the middle part of (\ref{eq.5.3}) is exactly the $S$-matrix of
Definition (\ref{defSmatrix}). This gives the symplectic sum formula
in the general case.

\bigskip

\begin{theorem}[Symplectic Sum Formula]\   Let $(Z,U,W)$ be the symplectic
sum of $(X,U,V)$ and  $(Y,V,W)$ along $V$.  Suppose that
$\alpha\in{\Bbb T}(Z)$ is supported off the neck as in Example
\ref{al.off.neck}. For any fixed decomposition $(\al_X,\;\al_Y)$ of
$\al$ the relative TW invariant of $Z$ is given
in terms of the  invariants of $(X,U,V)$ and $(Y,V,W)$   and the $S$-matrix
(\ref{defSmatrix})
by
\bear
TW_Z^{U,W}(\al)\ =\   TW_X^{U,V}(\al_X) \, * \,
S_V\, * \, TW_Y^{V,W}(\al_Y).
\label{mainthmformula}
\eear
\label{mainthm}
\end{theorem}
In fact, the Theorem holds more generally when $\al$
separates as in  Definition (\ref{al.separates}), except that the
definition of the  $S$-matrix needs to be enlarged. Instead of restricting
$TW_{\Bbb F}^{V,V}$ to $\al=1$ we restrict it to the subtensor
algebra ${\Bbb T}_V$ of ${\Bbb T}({\Bbb F})$ generated by the kernel
of the composition
\best
H^*({\Bbb F})\ma\longrightarrow^{i^*}
H^*(S_V)\ma\longrightarrow^{p_*} H^*(V)
\eest
where $S_V$ is the circle bundle on $N_V$, $p_*$ is the integration
along its fiber and $i:S_V\ra {\Bbb F}$ is the
inclusion. In that case we get an $S$-matrix defined by
\bear
S_V=(TW^{V,V}_{\Bbb F}|_{{\Bbb T}_V})^{-1}
\label{S.gen}
\eear

\bigskip

In the important   case  when $U$ and $W$
are empty   Theorem \ref{mainthm}  expresses the absolute
invariant of $Z$ in terms of the relative invariants of $X$ and $Y$.

\medskip

\begin{theorem}  Let $Z$ be the symplectic  sum of $(X,V)$ and
$(Y,V)$ and suppose that $\alpha\in{\Bbb T}(Z)$ separates as
$(\al_X,\al_Y)$ as in Definition (\ref{al.separates}). Then
\bear
TW_Z(\al)\ =\ (TW_X^V*S_V*TW_Y^{V})(\al_X,\al_Y).
\label{full.deg.form}
\eear
where $S_V$ is the $S$-matrix (\ref{S.gen}).
\label{cor4.2}
\end{theorem}
If moreover $\al$ decomposes as $\al=\al_X\otimes \al_Y$ then 
(\ref{full.deg.form}) becomes 
\best
TW_Z(\al)\ =\ TW_X^V(\al_X)*S_V(\al_V)*TW_Y^{V}(\al_Y)
\eest
where $\al_V\in {\Bbb T}_V$ is the pullback to ${\Bbb F}$ of the 
restriction of $\al$ to $V$.
\medskip

As a check, it is interesting to
verify the  symplectic sum formula in one very simple case where  the
GW invariant is simply the euler
characteristic.

\medskip

\begin{ex} Consider the $(J,\nu)$-holomorphic maps from an
elliptic curve $C$ with fixed complex structure representing the class
0.  When $\nu=0$ all such maps are maps to a single point, so the
moduli space is $X$ itself.  Furthermore, the fiber of the obstruction
bundle at a constant map $p$ is $H^1(T^2, p^* TX)$, which is naturally
identified with $T_pX$.  The (virtual) moduli space for $\nu\neq 0$
consists of the zeros of the generic section $\ov\nu=\int_C \nu$ of
this obstruction bundle $TX\to X$.  Thus this particular GW invariant
is $\chi(X)$.

Similarly, when $\nu=0$ the moduli space of $V$-regular curves is
$X\setminus V$ and its $V$-stable compactification, defined in
\cite{ip4}, is $X$. To compute the GW invariant relative to $V$, we
need to know how many of these point maps become $V$-regular after we
perturb to a generic $V$-compatible $\nu\ne 0$.  Because any
$V$-compatible $\nu$ is tangent to $V$ along $V$ the corresponding
section $\ov \nu$ has $\chi(X)$ zeros on $X$, out of which $\chi(V)$
lie on $V$. Thus the relative invariant is $GW_X^V=\chi(X)-\chi(V)$.
Note that $\chi({\Bbb F_V})=2\chi(V)$, so the $S$-matrix is the
identity in this case.  The symplectic sum formula therefore reduces
to the formula 
\best 
\chi(X)+\chi(Y)-2\chi(V)=\chi(X\#_V Y).  
\eest
\end{ex}
Much more interesting examples will be given in Section 15.

\bigskip

     Finally,  can also include   $\psi$  and $\tau$ classes as
constraints.  Recall that $\phi\in H^2(\ov\M_{g,n})$ is the first
chern class of ${\cal L}_i$, the relative cotangent bundle over at the
$i$th marked point. There is similar bundle $\wt{\cal L}_i$ over the
space of stable maps whose fiber at a map $f$ is the cotangent space
to the (unstabilized) domain curve, and whose chern class is denoted by
$\psi_i$.  It is also useful to pair each $\psi_i$ class with an
$\alpha_i\in H^*(Z)$ and consider the `descendent' $\tau_k (\al_i)=
ev_i^*(\alpha_i)\cup  \psi^k_i$.  It is a straightforward exercise,
left to the reader, to incorporate these constraints into Theorems
\ref{mainthm} and \ref{cor4.2}.

\vskip.4in


\setcounter{equation}{0}
\section{Constraints Passing Through the Neck}
\bigskip

Not every constraint class $\al\in H^*(Z)$ separates as in Definition
\ref{al.separates}. Yet for applications it is useful to have a
version of the symplectic sum formula for more general constraints ---
ones whose Poincar\'{e} dual cuts across the neck.  Since the
Poincar\'{e} dual of $\al\in H^*(Z)$ restricts to a class in
$H_*(X,V)$ such a general symplectic sum formula will necessarily
involve relative TW invariants of classes $\al\in H^*(X\setminus V)$.
That requires generalizing the relative invariant $TW^V_X$, which was
defined in \cite{ip4} only for constraints in $H^*(X)$.

\medskip

We begin by recalling the `symplectic compactification' of $X\setminus
V$ which was used in \cite{ip4}.  Let $\hat{X}$ be the manifold
obtained from $X\setminus V$ by attaching as boundary a copy of the
unit circle bundle $p:S_V\to V$ of the normal bundle of $V$ in $X$,
and let $p:\hat{X}\ra X$ the natural projection.  Suppose that $Z$ is
a symplectic sum obtained by gluing $\hat{X}$ to a similar manifold
$\hat{Y}$ along $S$.  We can then consider stable maps in $Z$
constrained by classes $B$ in $H_k(Z)$, i.e. the set of stable maps
$f$ with the image $f(x)$ of a marked point lying on a geometric
representative of $B$.  Restricting to the $\hat{X}$ side, such a
geometric representatives define constraints associated with classes
in $H_*(\hat{X},S)$.

Specifically, given a class $B\in H_*(\hat{X},S)$, we can find a
pseudo-manifold $P$ with boundary $Q$ and a map $\phi:P\to X$ so that
$\phi(Q)\subset S$ that represents $B$ and use this to cut-down the moduli
space.  Thus for generic $(J,\nu)$
$$
\ep_V\left(\ov\M_s^V(X,A)\right)\cap p(\phi(P))
$$
defines a orbifold with boundary that we denote  by
\bear
{TW}_{X,A,s}^V(\phi).
\label{phidependentGW}
\eear 
After cutting down by further constraints of the appropriate
dimension, this reduces to a finite set of points, giving numerical
invariants constructed using $\phi$. This is particularly simple when $B\in
H_*(X\setminus V)$, i.e. when $B$ can be represented by a map into
$\hat{X}\setminus S$.  The cobordism argument of Theorem 8.1 of
\cite{ip4} then shows that the relative invariants
(\ref{phidependentGW}) are well-defined. Note that these relative
invariants depend on $B\in H_*(X\setminus V)$ not on its inclusion
$B\in H_*(X)$. For example, rim tori and the zero class in
$H_2(\hat{X},S)$ have the same image under $p:\hat X\ra X$, but might
have different invariants (\ref{phidependentGW}).

\medskip

In general the constrained invariant (\ref{phidependentGW}) will not be
well-defined but will depend on the choice of $\phi$. The space
\bear
{\cal J}^V \ti \mbox{Maps}((P,Q),\,(\hat X,S))
\label{chamberedspace}
\eear
has a subset
$$
W\ =\ \bigcup_{i=1}^n\ \{ (J,\nu,\phi)\ |\ \mbox{there is a $V$-stable
$(J,\nu)$-holomorphic map $f$ with }
f(x_i)\in p(\phi(Q))\subset V\
\}
$$
where for some map one of the marked points $x_i$ lands on the
projection of $\phi(Q)$ into $V$.  Except in special cases, $W$ will
have codimension one, and thus will form walls which separate
(\ref{chamberedspace}) into chambers.

\begin{lemma}
The number (\ref{phidependentGW}) is constant within a chamber.  When
$B=[\phi]$ satisfies $p_*[\bd B]=0$ then
there is only one chamber, and therefore  (\ref{phidependentGW}) depends
only on $B$.
\label{chamberlemma}
\end{lemma}

\pf Any two pairs $(f,\phi)$ that lie in the same chamber can be
connected by a path $(f_t,\phi_t)$ with $f_t(x_i) \in
\phi_t(P\setminus Q)$.  The cobordism argument of Theorem 8.1 of
\cite{ip4} then proves the first statement.

Each $B$ in the kernel of $p_*\bd$ can be represented by a map $\phi$
as above with $\phi(Q)$ of the form $p^{-1}(R)$ for some $k-2$ cycle
$R$ in $V$. After restricting the last factor of
(\ref{chamberedspace}) to such $\phi$, the wall $W$ has codimension
two, giving the second statement.  \qed

\medskip

The following lemma relates the invariants associated with different chambers.

\medskip

\begin{lemma}
\begin{enumerate}
\item If $\phi_1,\phi_2:P\to X$ are two maps that agree on $\partial P$
then
\best
{TW}_X^V(\phi_1)={TW}_X^V(\phi_2)+TW_X^V(a)
\eest
where $a=[\phi_1\#(-\phi_2)]\in H_*(X\setminus V)$.
\item If $\phi_1,\phi_2$ define the same class in
$H_*(X,V)$ then we can find $\phi':R\to S$ where $\bd R=Q_1\sqcup (-Q_2)$
such that $\phi'$ agrees with $\phi_1$ on $Q_1$ and agrees with
$\phi_2$ on $Q_2$. Then $\phi_1$ and  $\phi_2\# \phi'$ have the same
boundary.
Moreover,
\best
{TW}_X^V(\phi_2\# \phi')={TW}_X^V(\phi_2)+TW_X^V\cdot {TW}_F^{VV}(\phi')
\eest
\end{enumerate}
\label{chamberlemma.2}
\end{lemma}

This actually means that in order to extend the definition of the relative
invariants from \cite{ip4}, we  only need to pick one geometric representative
$B$ ({\em any one})
such that $[B]\in H_*(X,V)$, $[\bd B]=\beta$ for each
$\beta\in \ker[H_{*-1}(S)\ra  H_{*-1}(X)]$.

Altogether, the invariants can be thought as giving (non-canonically)
a map
\bear
{TW}_X^V:{\Bbb T} (X\setminus V)\longrightarrow H_*(\M\ti {\cal H}_X^V)
\label{13.end}
\eear
although they depend on the actual representatives for the class $\al$ as
described in Lemma \ref{chamberlemma.2}.

\medskip

 With this extended definition of the relative invariants the  proof of Theorem
\ref{mainthmflatcase} carries through.  That proof began by choosing 
geometric representatives of
constraints $\al$ which separate.  For a general constraint $\al\in
H^*(Z)$ we can still choose a geometric representative $B$ of the
Poincar\'{e} dual, and consider its restrictions $B_X$ and $B_Y$ to
$(\hat X, S)$ and  $(\hat Y, S)$ respectively. The remainder of the 
proof still applies, giving
a sum formula relating the invariants $TW_Z(\al)$ of $Z$   to the relative TW
invariants (\ref{13.end}) of $X$ and $Y$ cut down by the constraints 
$B_X$ and $B_Y$.

\vskip.4in


\setcounter{equation}{0}
\section{Relative GW Invariants in Simple cases}
\label{Relative Invariants in Simple cases}
\bigskip

The symplectic sum formula of Corollary \ref{cor4.2}  expresses the
invariants of $X\# Y$ in terms of
the relative invariants of $X$ and $Y$.  In the next section we will apply
that formula to spaces that
can decomposed as symplectic sums where the spaces on one or both sides are
simple enough that
their relative invariants are computable.  That strategy can succeed only
if one has a collection of
simple spaces with known relative invariants.  This section provides four
families of such simple
spaces.

\medskip

     In some of the examples below  the set ${\cal R}$ of rim tori is
non-trivial.  In those
cases we will give formulas for  the invariants $\overline{GW}^V_X$
defined in the appendix although, as the
examples will show, it is sometimes possible to  compute the
${GW}^V_X$ themselves even though there are rim tori present.

\bigskip

\subsection{Riemann Surfaces}
\label{torus subsection}

For Riemann surfaces one can consider the
GW invariants as
absolute invariants or relative to a finite set of points.  These
invariants count coverings, and
the  homology class $A$ is simply the degree $d$ of the covering.

     In dimension two the symplectic sum is the same as the
ordinary connect sum --- one joins two Riemann surfaces by
identifying a point on one with a point on another, and then smooths.  Of
course, to apply the  sum formula one must first find $S_V$, which in  this
case is built from the
     relative invariants of $({\Bbb P}^1, V)$ where
$V=\{p_0,\;p_\infty\}$ two distinct points and where the constraints lie on
$V$.  In that context, we
fix  a nonzero degree $d$ and two sequences $s, s'$ that
describe the multiplicities of points  at the preimages of $p_0$ and
$p_\infty$ respectively.

\begin{lemma} The invariants $GW^V_{d,g,s,s'}$  with no
constraints except those on $V=\{p_0,\;p_\infty\}$ vanish
except  when $g=0$ and  $s$ and $s'$ are single points with
multiplicity $d$.  In that case
\best
GW^V_{d,0,s,s'}=1/d
\eest
     Moreover, in dimension two the $S$-matrix  is always the identity.
\label{P1proposgator0}
\end{lemma}
\pf This invariant is the oriented count of the 0-dimensional components of
$\ov{\M}^V_{d, g, s,s'}$.
But using (\ref{2.maindimformula})
\best
\dim \M^V_{d, g, s,s'}\ =\  2d+2g-2+\ell(s)-\deg s+\ell(s')-\deg s\ =\ 2g-2+
\ell(s)+\ell(s')
\eest
     is zero only if $g=0$ and $\ell(s)=\ell(s')=1$, i.e. $s$ and $s'$
specify single points with multiplicity $d$. If we stabilize, there is
only one such map, given by the equation $z\ra z^d$, so it's
contribution to $GW^V_{d,0,s,s'}$ is $1/d$. This map is
${\Bbb F}$-trivial, and hence doesn't contribute to the $S$-matrix. \qed

\medskip

The same dimension count gives the invariant with one constraint:

\medskip

\begin{lemma} The invariants $GW^V_{d,g,s,s'}(b)$  with one fixed branch
point and no other constraints except those on $V=\{p_0,\;p_\infty\}$ vanish
except  when $g=0$ and $\ell(s)+\ell(s')=3$, in which
case $GW^V_{d,0,s,s'}=1$.
\label{p1invt}
\end{lemma}

Perhaps the most interesting two-dimensional example  is the $g=1$
invariant of the torus
$T^2$.
\begin{lemma}
\label{8.toruslemma}
The $g=1$ invariants of the torus relative to a  set $V$ of $k\geq 0$
points form a series
$$
GW^V_{1}(T^2)\ =\ \sum GW^V_{d,1}(T^2)\,t^d
$$
that is equal to  the generating function for the sum of the divisors
$\sigma(n)=\sum_{d|n} d$, namely
\bear
G(t)\ = \ \sum_{n=1}^\infty \sigma(n)\,t^n\ =\  \sum_{d=1}^\infty
\frac{d\,t^d}{1-t^d}.
\label{defG(t)}
\eear
\end{lemma}

\medskip

\pf This is a matter of counting the (unbranched) covers of the torus.  That
was done in \cite{ip1} for $k=0$. In general, for  each degree $d$
cover each point  of $V$ has $d$ inverse images, each
with  multiplicity one.  Following the  notations of \cite{ip4} we order
the inverse images and divide by $d!$, leaving us with $G(t)$ again.
\qed

\bigskip

\subsection{$T^2\ti S^2$}
\label{torus-sphere subsection}

      Next we consider the $g=1$ invariants  of $X=T^2\ti S^2$.  Thinking of
this as an elliptic fibration over
$S^2$, we  fix a   a section $S$  and two disjoint fibers $F$ and
denote the corresponding  homology
classes  by $s$ and $f$.  Focusing on the classes $df$ and $s+df$ for
$d\geq 0$, we can
form generating functions for the absolute GW
invariants and the GW invariants relative to  one or two copies of the
fiber.

First consider  the classes $df$, where the invariants $GW_{df,1}$,
$GW^{F}_{df,1}$, and
$GW^{F,F}_{df,1}$ have dimension 0 by (\ref{2.maindimformula}).  There
are no rim tori in $X\setminus F$, and when $V$ is one or two copies
of the fiber we have $\ell=d\cdot f\cdot V=0$, so $V^\ell$ is a point in
(\ref{HVcover}).  Therefore $GW^F_{df,1}$  has values in $H_2(X)$ and
$GW^{F,F}_{df,1}$ has values in ${\cal H}^V=H_2(X)\times {\cal R}$.
Thus all three invariants can be written as  power series with
numerical coefficients.

\begin{lemma}  The genus one   invariants $GW$ and $GW^F$ in the
classes $df$ are given by
\best
\sum_d GW_{df,1}\,t_f^d \ =\ 2G(t_f) \qquad\mbox{and}\qquad \sum_d
GW^F_{df,1}\,t_f^d\ =\ G(t_f)
\eest
with $G(t)$ as in (\ref{defG(t)}) .  The corresponding relative
invariants $GW^{F,F}$ are indexed by
classes
$df+R$
for rim tori  $R$ and these all vanish:
\best
\sum_d GW^{F,F}_{df,1}\,t_{df+R}=0.
\eest
\label{L.14.4}
\end{lemma}

\pf  The generic complex structure on a topologically trivial  line bundle
over $T^2$ admits no non-zero holomorphic sections.  After  projectivizing,
we get a complex
structure on $T^2\ti S^2$ for which  the only holomorphic curves
representing $d f$ are multiple
covers of the zero section $F_0$ and the infinity section $F_\infty$. This
is a generic
$V$-compatible structure for $V=F_0$ or $F_0 \cup F_\infty$. As in
Lemma \ref{8.toruslemma} these contribute $G(t)$ to the power series for these
invariants. (Note that for the
relative invariant, we compute only the  contribution of curves that have
no components in $V$).
\qed

\medskip

The invariants for the classes $s+df$ are more complicated.  By
(\ref{2.maindimformula})
     the corresponding moduli spaces have dimension 4, so become points
in ${\cal H}_X^V$ after
imposing two point
constraints; these constraints can be either points $p\in X\setminus V$, or
$C_1(q)$, a contact of order 1 to $V$ at a fixed point $q\in V$.
Again  rim tori $R$ appear only for the invariant relative   to two copies
of a fiber.

\begin{lemma}  The genus one invariants $GW$ and $GW^F$ in the
classes $s+df$, $d>0$, are
\best
\sum_d GW_{s+df,1}(p^2)\,t_f^d\ =\  2G'(t_f) \qquad\mbox{and}\qquad
\sum_d GW^F_{s+df,1}(p;C_1(p))\,t_f^d\ =\  G'(t_f).
\eest
The corresponding relative invariants $GW^{F,F}$ can be indexed by  classes
$s+df+R$
for rim tori  $R$ and those with two point constraints on $V$ vanish:
\best
     \sum_d GW^{F,F}_{s+df+R,1}\left(\beta\right)\,t_{s+df+R} \ =\
\left\{\begin{array}{ll}
2G'(t) \qquad & \mbox{if } \beta=p^2 \mbox{ and } R=0,\\
G'(t) \qquad & \mbox{if } \beta=p;C_1(p) \mbox{ and } R=0,\\
0 \qquad & \mbox{if } \beta=C_1(p);C_1(p).
\end{array}\right.
\eest
\label{L.8.2}
\end{lemma}
\pf We can compute using the product structure $J_0$ on $T^2\ti S^2$.
Consider a $J_0$-holomorphic map representing $s+df$, passing through
generic points $p_1$ and $p_2$, and whose domain is a genus 1 curve
$C=\cup C_i$. The projection onto the second factor gives a degree 1
map $C\to S^2$, so $C$ must have a rational component $C_0$ which
represents $s$. The projection of the remaining components is zero in
homology, therefore they are multiple covers of the fibers. Because
the total genus is one there is only one such component.

Summarizing, for the product structure $J_0$ the only $g=1$
holomorphic curves representing $s+df$ have two irreducible
components, one of them a section $S$, and the other a multiple cover
of a fiber $F\notin V$.  The constraints require that $S$ pass through
$p_1$ and $F$ pass through $p_2$, or vice versa.  For each of those
two cases there are $d$ choices of the marked point on the domain of
$F$, so the count is the same as in   Lemma \ref{L.14.4}
with $G(t)$ replaced by $G'(t)$.  This gives the first formula.

The count for the second formula is similar.  Any $V$-regular genus 1
holomorphic map through an interior point $p$ and a point $q\in V$ has
two components: a section through $g$ and a $d$-fold cover of a fiber
$F$ through $p$.  The fiber domain can be marked in $d$ ways, giving
the count $G'(t)$.

For the invariant relative two copies of the fiber $F$, there are rim
tori, but the discussion above implies that for $J_0$ the only
holomorphic curves in the classes $s+df +R$ appear only for $R=0$
(where these curves {\em define} what $R=0$ means).  \qed

\bigskip

\subsection{Rational Ruled Surfaces}

Here let ${\Bbb F}_n$ be the rational ruled surface whose fiber $F$,
      zero section $S$ and  infinity section $E$ define homology
classes with $S^2=-E^2=n$.  We will compute some of the relative
invariants $GW^{V}$ with $V=S\cup E$ and with no constraint on the
complex structure of
the domain ($\kappa=1$).

\medskip

Fix a non-zero class $A=aS+bF$ and two sequences $s, s'$ of
multiplicities that describe the intersection with $S$ and $E$
respectively. The relative $GW$ invariant with no constraint on the
complex structure and $k$ marked points lies in the homology of the
moduli space in $X^k\ti S^\ell\ti E^{\ell'}$ with $\ell=\ell(s)$ and
$\ell'=\ell(s')$. After
imposing constraints $\alpha=(\alpha_1, \dots,\alpha_k)$
\best
GW^{S,E}_{A,g, s,s'}(\alpha)\in   H_*(S^\ell)\otimes H_*(E^{\ell'})
\eest
where $S\cong E \cong \P^1$.  Noting that the canonical class of
${\Bbb F}_n$ is $K=-2S+(n-2)f$ and
$\deg s =  E\cdot A =b$ and $\deg s'=  S\cdot A=b+na$,  we have
\best
{1\over 2} \mbox{dim } GW^{S,E}_{A,g, s,s'}(\alpha)
& = & (n+2)a+2b+g-1 -(\deg s -\ell(s))-(\deg s' -\ell(s')) -\deg \alpha\\
& = & 2a+g-1 +\ell+\ell' -\deg \alpha
\eest
      But $SV_s$ has dimension ${\ell(s)}$, so the moduli space
represents zero in
homology
unless $\mbox{dim } \M_{g,k, s,s'}({\Bbb F_n},A)\leq \ell + \ell'$,
so we always have
\bear
2a+g\ \leq\ 1+\deg\alpha.
\label{dimensionct1}
\eear

\begin{lemma} The invariants $\displaystyle GW_{A,g,s,s'}^{S,E}$ with no
constraints except those on $V=S\cup E$   vanishes except when
$A=bF$, $g=0$, and
$s$ and $s'$ are
single points with multiplicity $b>0$.  In that case
\best
      GW^{S,E}_{bF,0,s,s'}={1\over b} \l(S\otimes1\,+\, 1\otimes E\r).
\eest
Moreover, the $S$-matrix in ${\Bbb F}_n$ vanishes.
\label{CHproposgator0}
\end{lemma}
\pf It suffices to show  that the only contributions to $GW$ from classes
$A=aS+bF$  come from unstable rational domains with $a=0$, i.e. from
     ${\Bbb F}$-trivial maps.   Taking
$\kappa=\alpha=1$, (\ref{dimensionct1})   implies that
$A=bF$ and $g=0$ or $1$.  Moreover, because
every $bF$ curve intersects both $E$ and $S$, we have $\ell
+\ell'\geq 2$, and when $g=0$ stability of the domain requires that $\ell
+\ell'\geq 3$.  In these cases the moduli space $\M^V_{g,s,s'}({\Bbb
F},bF)$  is either empty or has
dimension $\geq 2$.

Suppose that the moduli space is non-empty and the above stability
conditions hold.  Since $E$ and $S$
are copies of $\P^1$, $H_*(S^\ell)\otimes H_*(E^{\ell'})$ is
generated by point or $[\P^1]$
constraints.  Then for each generic $(J,\nu)$ there are maps  $f$ in
the moduli space
whose images passes through at least two fixed points $p,q\in E\cup S$ in
generic position.    Take  $(J,\nu)\to (J_0,0)$ where $J_0$  is a
complex structure with a holomorphic projection
$\pi:{\Bbb F}_n\to \P^1$.  In the
limit we obtain  a connected stable map $f_0$ through $p$ and $q$ with
components representing $a_iS+b_if$ such that
$bF=\sum a_iS+b_if$.  But then each  $a_i=0$, so the
image of $\pi\circ f_0$ is a
single point containing $\pi(p)$ and $\pi(q)$.  This cannot happen for
generic $p,q$.

Thus $\M^V_{g,s,s'}({\Bbb F},bF)$ consists of  ${\Bbb
F}$-trivial maps (c.f. Definition \ref{5.relativelystable})
representing $A=bF$.  Such maps do
not appear in the
$S$-matrix. \qed

\bigskip

\begin{lemma} Fix a point $p\in {\Bbb F}\setminus V$ with $V=E\cup S$. Then
$\displaystyle GW^{V}_{A,g,s,s'}(p)$ vanishes except in the following cases:

\smallskip

(i) $\displaystyle GW^{S,E}_{bF,0,s,s'}(p)=1$ when $s$ and
$s'$ are single points with multiplicity $b>0$.

\smallskip

(ii)  $\displaystyle GW^{S,E}_{S+bF,0,s,s'}(p)= SV_s\ti
SV_{s'}$ whenever  $\deg s=b,\; \deg s'=b+n$.
\label{CHInvt1}
\end{lemma}
\pf From  (\ref{dimensionct1}) we have $GW^{V}_{aS+bF,g,s,s'}(p)=0$
     unless  $2a+g \leq 2$.  Thus either (i) $a=0$,  or (ii) $a=1$ and $g=0$.

\medskip

In case (i) each map contributing to the invariant represents $bF$,
passes through
$p$, and hits $E$
and $S$.   Hence  $\mbox{dim } \M^V_{g,s,s'}({\Bbb F},bF)=g-1+\ell(s)
+\ell(s')$ with $\ell(s)
+\ell(s')\geq 2$.  The limiting argument used in Lemma
\ref{CHproposgator0} then shows that
$GW^{V}_{bF,g,s,s'}(p)$ vanishes unless $g=0$ and $\ell(s) =\ell(s')=1$. Thus
$s$ and $s'$ are single points of multiplicity $b$, and the maps pass
through $p$. Moving to the fibered
complex structure, one sees that there is a unique such
stable map for each $b>0$.  This gives (i).

\medskip

In case (ii) the moduli space $\M^V_{0,s,s'}({\Bbb F},S+bF)$ has
dimension $\ell(s) +\ell(s')$ and is
empty unless  $b\geq 0$.  That means $\M^V_{0,s,s'}({\Bbb F},S+bF)$ is
a multiple of $SV_s\times
SV_{s'}$, so  invariant vanishes except when  all contact points on
$E$ and $S$ are fixed.  By the adjunction
inequality, any irreducible curve $C$ representing $S+bF$ is rational
and embedded,  so we can compute the invariant by intersections in
$\P(H^0({\Bbb F}_n,{\cal O}_{{\Bbb F}_n}(S+bF))$  (the
standard complex structure on ${\Bbb
F}_n$ is generic for these curves $C$ because $h^{1}(C;{\cal O}(S+bF)|_C)
= h^{1}(\P^1;{\cal O}(n+2b))=0$).  But
$h^0({\Bbb F}_n, {\cal O}(S+bF))=n+2+2b$, and each of the conditions
imposed (including multiplicities) are
linear conditions.    Thus the number of curves representing $S+bF$
passing through a point $p$ and meeting  $E$ and $S$ at fixed contact
points is 1.
\qed

\bigskip

\subsection{The Rational Elliptic Surface}
\label{rational elliptic subsection}

     As a final example we  consider the rational elliptic surface
$E$. Let $f$ and $f$ denote, respectively, the homology classes of a
fiber and a fixed section of an elliptic fibration $E\to\P^1$.
The following lemma describes the invariants  relative to a fixed
fiber $F$ in the classes $A=s+df$ where $d$ is an
integer.  In this case there
{\em are} rim tori in $E\setminus F$, suggesting that one use the average
invariant $\ov{GW}$ defined in
the appendix.  However, the lemma shows that the
average contains only  only one non-zero term (as happened in the last case
of Lemma \ref{L.8.2}).

\begin{lemma}  The genus $g$ relative and absolute invariants  of $E$
in the classes $s+df\in H_2(E)$ are related by:
\best
     GW_{s+df,g}(p^g)\ =\ \overline{GW}^F_{s+df, g}(p^g; C_1(f))\ =\
GW^F_{s+df, g}(p^g;
C_1(f))
\eest
where the second equality means that $GW^{F}$ can be  indexed by  classes
$s+df+R$ for rim tori  $R$ and
these vanish whenever $R\ne 0$.
\label{L.avg=abs}
\end{lemma}
\pf The first equality holds because generically all maps contributing to
the absolute invariant are
$V$-regular. That is true because if some component of a stable map is
taken into $V=F$, then that
component must have genus at least 1. But then the remaining components
have genus less than $g$, so
cannot pass through $g$ generic points.

The second equality follows from a projection argument
like the one used for Lemma \ref{L.8.2}. Consider a curve  be
     $C=\cup C_i$ representing $s+df$ which is  holomorphic for a fibered
complex structure $J_0$ on $E$.
Since the projection to $\P^1$ gives a degree one composition, $C$ must
have a rational component
$C_0$ that intersects each fiber in exactly
one point, while the other components are multiple covers of fibers,  so
represent $df\in H_2(E\setminus F)$. Moreover, $C_0$ is an embedded
section representing $s$. Since $s^2=-1$
then $C_0$ must be the unique holomorphic curve in the class $s$. Thus
the only curves in the class $s+df+R$ appear only for $R=0$.
\qed

\medskip

The invariants of Lemma \ref {L.avg=abs}  will be explicitly computed in
section \ref{CurvesinE}.

\bigskip

\subsection{\bf Rational relative invariants}
\label{Rational relative invariants}

Counting rational curves requires only the $g=0$ relative invariants and the
corresponding $S$-matrix.  The following two propositions show that
these are particularly simple: the $S$-matrix is the identity and the
relative invariant is the same as the absolute invariant in the
absence of rim tori.

\medskip

\begin{prop} When $g=0$, $s=(1,\dots, 1)$ and $A\in H_2(X)$, the
relative invariant (summed over rim tori as in (\ref{def.avgGW})) equals
the
absolute invariant:
\best
{1\over \ell(s)!}\; \ov{GW}^V_{A,0}(\al;C_s(\gamma))=GW_{A,0}(\al;
i_*(\gamma))
\eest
where $\al=(\al_1,\dots,\al_n)\in (H_*(X))^n$,
$\gamma= (\gamma_1,\dots,\gamma_\ell)\in (H_*(V))^\ell$ and
$i_*:H_*(V)\ra H_*(X)$ is the inclusion.
\end{prop}
\pf Fix a generic $V$-compatible pair $(J,\nu)$. Recall that $(J,\nu)$ is
generic for curves that have no components in $V$, and also its
restriction to $V$ gives a generic pair on $V$. However, for a curve
entirely contained in $V$, even though $(J,\nu)$ is generic when the curve
is considered in $V$, it might not be generic when the curve is considered
in $X$.

For any genus $g$ and ordered sequence $s$, consider the natural
inclusion:
\best
\M_{X,A, g,s}^V\hookrightarrow \M_{X,A,g}
\eest
where $A\in H_2(X)$, so on the left we took the union over all rim tori.
When $s=(1,\dots,1)$, any element in $\ov\M_{X,A,g}$ that has no
components in $V$ is in fact an element of $\M_{X,A,g,s}^V$.  We will show
that for generic $V$-compatible $(J,\nu)$, when $g=0$ the contribution of
the moduli space of curves with some components in $V$ to the absolute
invariant vanishes, and therefore the two invariants are equal.

For simplicity,  start with the case when $f$ has only one component, and
this is entirely contained in $V$. Then $A=i_*(A_0)$ with
$A_0\in H_2(V)$, and $\ell(s)=A\cdot V= c_1(N_XV )\cdot A$. Then the
moduli
space of such curves has
\best
\dim \M_{V,A_0,g}(i_*\gamma)&=&-K_V\cdot A_0+(\dim V-3)(1-g)-
\sum_{i=1}^\ell (\dim V-\dim \gamma_i)\\
&=&\dim \M^V_{X,A,g}(\gamma)-1+g
\eest
as in equation(6.4) of \cite{ip4}. This means that for genus $g=0$ the
dimension of the moduli space of curves entirely contained in $V$ is one
less
then the (virtual) dimension when considered as curves in $X$. Therefore
if
the virtual dimension in $X$ is 0, there are no curves in $V$ who could
contribute. The general case of a curve with some components in $V$ and
some off $V$ follows similarly. \qed

\bigskip

\begin{prop} The $g=0$ part of the $S$-matrix is the identity for any $V$
and
any normal bundle $N$.
\end{prop}
\pf  By (\ref{TW=1+R}) this statement is equivalent to showing that there
is no contribution to the $g=0$ GW-invariant coming from maps into
${\Bbb F}$ which are not ${\Bbb F}$-trivial. Consider the 0 dimensional
moduli space $\M_{{\Bbb F},A,0,s}^{V_0,V_\infty}(\gamma)$ constrained only
along $V_0$ and $V_\infty$, such that the corresponding $GW$ invariant is
not zero.  By Theorem 1.6 of \cite{ip4} the same moduli space would
be non-empty for the submersive structure associated with a generic
$\nu_V$ on $V$ (as defined before Lemma \ref{TW=1+R}).  Then each
$f\in \M_R$ would project to a map
$f_V$ in
$\M_{V,\pi_*A,0,s}$ that passes through the $\gamma$ constraints.
But  counting virtual dimensions
using  equation (6.4) of \cite{ip4}, we see that
\best
\mbox{dim}\ \M_{V,\pi_*A,0,s}(\gamma)\ =\ \mbox{dim}\
\M_{{\Bbb F},A,0,s}^{V,V}(\gamma) -
\mbox{index}\  D^N_s \ = \ 0+ g-1
\eest
is negative when $g=0$, so this moduli space is empty for generic $\nu_V$.
\qed

\vskip.4in

\setcounter{equation}{0}
\section{Applications of the Sum Formula}
\label{section15}
\bigskip

This last section  presents  three applications of
the sum formula:  (a) the Caporaso-Harris formula for the number of
nodal curves in $\P^2$, (b) the
formula for the Hurwitz numbers counting branched covers of $\P^1$,
and (c) the formula for the
number of rational curves representing a primitive homology class in
the rational elliptic surface.  These  formulas have all recently 
been established using Gromov-Witten invariants in some guise. Here we
show that all three follow rather easily from the symplectic sum formula.

\subsection{The Caporaso-Harris formula}
\label{CHsubsection}

Our first application is a derivation of the Caporaso-Harris recursion
formula for the number
$N^{d,\delta}(\al,\b)$ of curves in $\P^2$ of
degree $d$ with $\delta$ nodes, having a contact with $L$ of
order $k$ at $\al_k$ fixed
points, and at $\b_k$ moving points, for $k=1,2,\dots$ and passing through
the appropriate number $r$ of generic fixed points in the complement of $L$.

For this we consider the pair $(\P,L)$, which can be written as
a symplectic connect sum:
\bear
(\P^2,L)\ma\#_{L=E} ({\Bbb F}_1,E,L)\ = \ (\P^2,L)
\label{5.CH1}
\eear
where $({\Bbb F}_1,E,L)$ is the ruled surface with Euler class one with
its
zero section $L$ and its infinity section $E$.  We can then get a
recursive formula for the TW invariant of $(\P^2,L)$ by moving one point
constraint $pt$ to the ${\Bbb F}$ side, and then using the symplectic sum
     formula.

The splitting (\ref{5.CH1}) is along a sphere $V=E=L$, so there are no rim
tori. The relative invariant  therefore lies in the homology of $SV$
    and is invariant under the action of the subgroup of the symmetric group
that switches the order of points of same multiplicity. A basis for this
homology is given by (\ref{last.basis2}), where $\{\gamma_i\}$ with
$\gamma^1=p$ a point and $\gamma_2=[\P^1]$ is a basis of  $H_*(V)$.

To recover \cite{ch} notation, for each sequence $(m_{a,i})$, denote
$\al_a=m_{a,1}$ and $\beta_a=m_{a,2}$, and  let
     $\alpha=(\alpha_1,\al_2,\dots)$,   $\beta=(\beta_1,\b_2,\dots)$.
Then with this change of coordinates,
\best
N^{d,\delta}(\alpha,\beta)\ =\ TW_{\chi,dL,\P^2}^{L} (p^r, C_{m})
\eest
where  $\chi-2\delta=-d(d-3)$ is the "embedded euler characteristic" and
$r=3d+g-1-\sum\al_i-
\sum (\b_j-1)$, and we are imposing no constraints on the complex
structure of the curves. Similarly, let
\best
N^{a,b,\chi}(\al',\b';\; p;\; \al,\b)\ =\
TW_{\chi, aL+bF,{\Bbb F}}^{E,L} (C_{m}; p; C_{m'})
\eest
denote the number of curves of Euler characteristic $\chi$ in ${\Bbb F}$
representing $aL+bF$ that have contact described by
$(\al', \b')$ along $E$, $(\al, \b)$ along $L$ and pass through an extra
point $p\in {\Bbb F}$ (we prefer to label these numbers using $\chi$
rather then the number of nodes).

By Lemma \ref{CHproposgator0}  the $S$-matrix vanishes. The symplectic sum
theorem then implies:
\best
     N^{d,\chi}(\al,\b)=\ma \sum |\al'|\cdot|\b'|
\cdot N^{d',\chi'}(\al',\b') \cdot N^{d-d',b,\chi''}(\b',\al';p;\al,\b)
\eest
where the sum is over all $\al'$, $\b'$ and all
decompositions of $(dL, \chi)$ into $(d'L,\chi')$
and $((d-d')L+bF,\chi'')$ such that
$\chi=\chi'+\chi''-2\ell(\al')-2\ell(\b')$.
Combining Lemmas \ref{CHproposgator0} and \ref{CHInvt1} we see that there
are  exactly two types  of curves that contribute to the
relative GW invariant $TW_{\P^1}^{E,L} (C_{s,\gamma}; p; C_{s,\gamma})$ of
${\Bbb F}$ with one fixed point $p$.
\begin{enumerate}
\item several $g=0$ unstable domain multiple covers of the fiber, one of
them say of multiplicity $k$ passing through  the point $p$, corresponding to
the situation $d'=d$ and
\best
\b'=\b + \ep_k;\; \al'=\al-\ep_k
\eest
where $\ep_k$ is the sequence that has a 1 in position $k$ and 0
everywhere else.
\item  several $g=0$ unstable domain multiple covers of the fiber
together with one $g=0$ curve in the class $L+a F$ passing through $p$ and
having all contact points with $E$ and $L$ fixed say described by
$\al_0'$ and $\al_0$; this corresponds to $d'=d-1$ and
the situation
\best
\al=\al_0+\al';\; \b'=\al_0'+\b; \quad \mbox{ equivalently  }
     \b'\ge\b;\;\al\ge \al'
\eest
\end{enumerate}
In each situation above, the number of $V$-stable curves is 1.
In the second case, note that there are ${\al\choose \al'}$ choices of
$\al_0$ and ${\b'\choose \b}$ of $\al_0'$. Moreover, for each
${\Bbb F}$-trivial curve its invariant combines with its corresponding
     multiplicities in $|s'|\ell(s')!$ to give 1.
Therefore, the  remaining multiplicity in case 1 is $k$, while in case
2, is $|\al_0'|=|\b'-\b|$.  Putting all these
together, we get:
\best
     N^{d,\delta}(\al,\b)= \ma \sum
     k N^{d,\delta'}(\al-\ep_k,\b+\ep_k) +
\ma \sum |\b'-\b|{\al\choose \al'}{\b'\choose \b} N^{d-1,\delta'}(\al',\b')
\eest
where the last sum is over all $\b'\ge \b$, $\al'\ge \al$. This is exactly the
Caporaso-Harris formula.

\vskip.4in

     \subsection{Hurwitz numbers}

The method of section \ref{CHsubsection} can also be applied for maps
into $\P^1$.  In
that case the  symplectic sum formula yields  the
cut and paste formula for Hurwitz numbers that was first proven using
combinatorics by Goulden, Jackson and Vainstein in \cite{gjv}. (Recently
Li-Zhao-Zheng \cite{lzz} have derived a similar formula using  \cite{lr}).

\medskip

The Hurwitz number $N_{d,g}(\al)$ counts the number of genus $g$,
degree $d$ covers of $\P^1$
that have the branching pattern over a fixed point $p\in \P^1$ specified by
the unordered partition $\al$ of $d$, while the remaining branch points are
simple and fixed.  We can get at these numbers by regarding the pair
$(\P^1,p)$  as a symplectic  sum:
\bear
(\P^1,p)=(\P^1,x)\ma\#_{x=y}(\P^1, y,p)
\label{p1.split}
\eear
We then get a recursive formula for the GW invariant of $(\P^1,p)$ by moving
one simple branch point $b$ to the $(\P^1, y,p)$ side and applying
      the symplectic sum formula.

In fact the Hurwitz numbers are the coefficients, in a specific
basis, of the GW invariants of $\P^1$ relative to a point $V=p$. More
precisely,  each unordered partition $\al=(\al_1,\al_2,\dots)$ of $d$
defines numbers $m_a=\#\{i\;|\; \al_i=a\}$; let $C_m$ be the
corresponding basis (\ref{last.basis2}) (in this case the basis
$\{\gamma_i\}$ of $H^*(V)$ has only one element). Then
\best
N_{d,g}(\al)=GW_{\P^1,d,g}^p(b^r; C_m)
\eest
is the number of degree $d$, genus $g$ covers that have the branching
pattern over $p\in \P^1$ determined by $\al$, and $r=2d-2+2g-2-\sum (a-1)m_a$
other fixed, distinct branch points. (Note that the branching order
is the order of contact to $p=V$). The corresponding generating
function (\ref{TW.gen.fcn})  is
\best
G\ =\ GW_{\P^1}^p\ =\ \sum GW_{\P^1,d,g}^p(b^r; C_m)\prod_{a}
{(z_{a})^{m_{a}}\over m_{a}!}\;\;{u^r\over r!}\;t^d\;\la^{2g-2}.
\eest

Now apply the symplectic sum formula to the decomposition
(\ref{p1.split}), putting $r-1$ branch points on the
first copy of $\P^1$ and one on the second copy.  Since there are no
rim tori and the $S$-matrix vanishes
by Lemma \ref{P1proposgator0} we obtain
\bear
GW_{d,g}^p(b^r;C_m)=\ma\sum {|m'|}\cdot
TW_{d,\chi_1}^p(b^{r-1}; C_{m'}) \cdot
TW_{d,\chi_2}^{p,p}(C_{m'};\; b;\; C_{m})
\label{hurwitzTWformula}
\eear
where the sum is over all $m'=(m'_1, m'_2,\dots)$ and all $\chi_1, \chi_2$
such that $2-2g=\chi_1+\chi_2-2\ell(m')$ and so that the attached
domain is connected.But $TW=\exp GW$ and Lemma \ref{p1invt} implies
that the only possibility for the last factor in
(\ref{hurwitzTWformula}) is a union of trivial spheres together with a
degree $a$ sphere constrained by $C_a$ at one end and $C_{i,j}$ with
$i+j=a$ at the other end (plus the branch point in the
middle). Therefore there are only
two possibilities for the other factor and for  the partition $\al'$
corresponding
to $m'$.
\begin{enumerate}
\item $\al=(i,j, \beta)$ and $\al'=(i+j,\beta)$ for some $i,j,\beta$,
so the covering map has
     genus $g$ and  degree $d$.
\item $\al=(a,\beta)$ and $\al'=(i,j,\beta)$ with $a=i+j$. Then
$\chi_1=2g-4$ so the covering map  is either genus $g-1$ and degree $d$
    or genus $g_1,g_2$ of degrees $d_1$, $d_2$.
\end{enumerate}
The sum formula (\ref{hurwitzTWformula}) can then be written as a
relation for the
generating function, namely
$$
\partial_u G={1\over 2} \sum_{i,j\ge 1}\l(
    ij\la^2 z_{i+j}\l[\partial_{z_i}\partial_{z_j}G+
\partial_{z_i}G\cdot\partial_{z_j}G\r]
+(i+j)z_iz_j\partial_{z_{i+j}}G\r)
$$
This is the `cut-join' operator equation $DG=0$ of \cite{gjv}.  It
clearly determines
the Hurwitz numbers recursively.   The same formula works to give the
`Hurwitz numbers' counting
branched covers of higher genus curves.

\vskip.4in
\subsection{Curves in the rational elliptic surface}
\label{CurvesinE}

We next consider the  invariants of the rational elliptic surface
$E\to\P^1$.  Using the notation of section \ref{rational elliptic
subsection} will focus on the
classes $A=s+df$ where $d$ is an integer.  The numerical
invariants
$GW_{A,g}(p^g)$ then  count the number of connected genus $g$ stable maps in
the class
$s+df$ through $g$ generic points (with no constraints on the
complex structure of the domain). For each $g$ these define  power series
$$
F_g(t)=\ma\sum_{d\ge 0} GW_{d,g}(p^g)\  t^{d}\,t_s
$$
where $t=t_f$.  Recently,  Bryan-Leung \cite{bl} proved that
\bear
F_g(t)=F_0(t) \left[G'(t)\right]^g\label{BLg}
\label{Fgformula}
\eear
with $G$ as in (\ref{defG(t)}) and
\bear
F_0(t)=t_s\, \left(\ma\prod_{d}{1 \over 1-t^d}\right)^{12}.
\label{BL}
\eear 
As mentioned in the introduction, this formula is related to the work 
of Yau-Zaslow \cite{yz} and
to more  general conjectures (such as those stated in \cite{go}) 
about counts of nodal curves in
complex surfaces.

     We will use our symplectic sum theorem to give a
short proof of this formula, beginning with the $g=0$ case.   The
proof is accomplished by
relating $F(t)$ to  the similar series of elliptic
($g=1$)  invariants
\best
H(t)&=&\ma\sum_{d\ge 0}\ GW_{d,1}(\tau_1[f^*]) \  t^{d}\,t_s
\eest
where $f^*\in H^2(E)$ is the Poincar\'{e} dual of the fiber class and where
$\tau_1[f^*] =\ev_1^*(f^*)\cap \psi_1$ is the corresponding `descendent
constraint' described at the end of section \ref{The General Sum Formula}.

We will compute $H$ in two different ways.  The first is based on the
standard method of `splitting the
domain', which yields the following general facts for 4-manifolds.

\begin{lemma}
Let $X$ be a symplectic 4-manifold with canonical class $K$. (a)\ For
$A=0$ and $g=1$ the
GW invariant with a single constraint $B\in H^2(X)$ is
\bear
GW_{0,1}(B)\ =\ \frac{1}{24} K\cdot B,
\label{ghosttoricount}
\eear
     (b)\ For any classes $A,f \in H_2(X)$
satisfying  $A\cdot K=-1$
$$
GW_{A,1}(\tau_1[f^*])\ =\ {(f\cdot A)\over 24}
      (A^2+K\cdot A)\,GW_{A,0}+
\ma\sum_{A_1+A_2=A\atop A_1\neq 0,\,A_2\neq 0}(f\cdot A_2) (A_1\cdot A_2)
GW_{A_1,1}GW_{A_2,0}.
$$
\label{TRRlemma}
\end{lemma}

\pf  (a)\ For $\nu=0$, $\overline{\M}_{1,1}(X,0)$ is the space
$\overline{\M}_{1,1}\ti X$  of `ghost tori' $f:(T^2,j)\ra X$ with  $f(z)=p$
a constant map.  At such $f$, the fiber of the obstruction bundle
is $H^1( T^2, f^*TX)=H^1(T^2, {\cal O}) \otimes TX$.  The dual of
the bundle $H^1(T^2,{\cal O})$ over $\overline{\M}_{1,1}$ is the
      Hodge bundle.  Since the first chern number of the  Hodge bundle is
$-1/24$,  the
euler class the obstruction bundle is
$$
\chi(X)[\overline{\M}_{1,1}]\otimes1+\frac{1}{24} 1\otimes K \in
H_2(\overline{\M}_{1,1}\ti X).
$$
For $\nu\neq 0$, the (virtual) moduli space is the zeros of a generic
section of the obstruction bundle, which consists of (i) maps from a
torus with any complex structure to $\chi(X)$ specified points of $X$
and (ii) maps from a torus of specified complex structure into some
point on the canonical divisor.  Generically, the images of the type
(i) maps will miss the constraint surface representing $A$.  The maps
of type (ii) give the formula (\ref{ghosttoricount}).

\medskip

(b)\ The genus 1 topological recursive relation says
$$
GW_{A,1}(\tau[f^*])\ =\  {1\over 24}  GW_{A,0}(H_\al,H^\al,f)
\ +\  \ma\sum_{A_1+A_2=A}\ma\sum_{\alpha} GW_{A_1,1}(H_\al)
GW_{A_2,0}(H^\al,f)
$$
where $\{H_\alpha\}$ and $\{H^\alpha\}$ are bases of $H^*(X)$
dual by the intersection form.  But for $A\neq 0$ $GW_{A,0}(H^\al,H^\b,f)$
vanishes by
dimension count unless $H^\al$ and $H^\b$ are two-dimensional, and then  each
$A$-curve hits a generic geometric representative of $H_\al$ at
\,$H^\al\cdot A$\, points counted with algebraic multiplicity.  A dimension
count also shows that the moduli spaces with $A_1=A$ and $A_2=0$ are of the
wrong dimension to contribute to
the double sum above.   Hence the expression above becomes
$$
{1\over 24}\sum
(H^\al\cdot A)(H_\al\cdot A)(f\cdot A)\,GW_{A,0}
\ +\  \ma\sum_{A_1+A_2=A\atop A_1\neq 0,\,A_2\neq 0}(H_\al\cdot
A_1)(H^\al\cdot A_2)(f\cdot
A_2) GW_{A_1,1}GW_{A_2,0}
$$
plus the term with  $A_1=0$, which by (\ref{ghosttoricount}) is
\best
\frac{1}{24} \left(K\cdot H_\alpha\right) GW_{A,0}\left(H^\alpha, f\right)\
=\ \frac{1}{24} \left(K\cdot H_\alpha\right)
      \left(A\cdot H^\alpha\right) \left(A\cdot f\right) GW_{A,0}.
\eest
The lemma follows  because $\sum (H_\al\cdot A_1)(H^\al\cdot A_2)=A_1\cdot
A_2$.
\qed

\bigskip

Taking $X$ to be the rational elliptic surface $E$, we can apply Lemma
\ref{TRRlemma}
with $A=s+df$.
Then $K=-f$, $A\cdot f=1$ and $A^2=2d-1$. The only
possible decompositions are $A_1=kf$ and $A_2=s+(d-k)f$ so:
\best
GW_{s+df,1}(\tau[f^*])\ =\ {d-1\over 12}GW_{s+df,0}+
\ma\sum_{k=1}^{d}k\ GW_{kf,1}\ GW_{s+(d-k)f,0}
\eest
But for  the rational elliptic surface the invariant $GW_{kf,1}(s)$ is
$\sigma(k)$ for
$k>0$. (Since in $\P^2$ there is a unique  cubic through 9 generic points.
      As in section 4 of \cite{ip1}, for each $k$ there are
$\sigma(k)$ distinct $k$-fold covers an elliptic curve with marked point, all
with positive
sign).  Because the marked point can go to any of $s\cdot kf=k$ points, this
means that the
unconstrained invariant is
$$
GW_{kf,1}\ =\  \sigma(k)/k \quad \mbox{for } k> 0.
$$
It follows that
\bear
H(t)&=& {1\over 12} \left(t\,F_0' - F_0\right) +F_0\cdot G.
\label{rel1}
\eear

\bigskip

On the other hand, we can calculate $H(t)$ by splitting the target
and using the symplectic sum
theorem.  Let ${\Bbb F}=T^2\ti S^2$, and
let
$F$ denote both a fiber of of the elliptic fibration $E$ and a fixed torus
$T^2\ti \{\mbox{pt}\}$ inside ${\Bbb F}$.  We can apply sum formula by  writing
$E=E\#_{F}{\Bbb F}$ for the class $A=s+df$ with the constraint
on the ${\Bbb F}$ side.  Since $A\cdot F=1$, the connected curves
representing $A$ split into the
union of  connected curves in $E$ and in ${\Bbb F}$; thus the symplectic sum
formula  applies for the $GW$ (as well
as the $TW$ invariants).

If we have a genus 0 curve on the ${\Bbb F}$ side in the class
$s+d_1F$, then by projecting onto the
$T^2$ factor and noting that there are no maps from $S^2$ to $T^2$ of
non-zero degree, we conclude that
$d_1=0$. But the moduli space of genus 0 curves in ${\Bbb F}$ representing
$s$ and passing through $F$ is
isomorphic to $F=T^2$,
and moreover the relative cotangent bundle to them along $F$ is isomorphic to
the normal bundle to $F$. So
\best
GW_{s,0}(\tau_1[f^*])= GW_{s,0}\left((f^*\right)^2)=0.
\eest
Thus there is no contribution from genus 0 curves on the ${\Bbb F}$ side or
in the neck
(which is also a copy of ${\Bbb F}$).  The same argument shows that there
are no rational curves in $F$, so
the $g=0$ absolute and relative invariants are the same.

With these observations, the only possibility is to have a genus 1 curve on
the ${\Bbb
F}$ side, genus 0  on  the $E$ side, and
no contribution from the neck.  The symplectic sum formula thus says
$$
GW_{d,1}(\tau_1[f^*]) \ =\ \sum_{d_1+d_2=d} GW_{s+d_1f,0}(E)  \cdot
GW_{s+d_2f,1}({\Bbb F})(\tau_1[f^*])
$$
This last invariant can be computed by applying the topological recursive
relation
to
$X={\Bbb F}$ just as in  Lemma \ref{TRRlemma}:
\best
GW_{s+df,1}(\tau_1[f^*])\,=\,\frac{d-1}{12}GW_{s+df,0}+
\ma\sum_{d_1+d_2=d\atop d_1\neq 0,\,d_2\neq
0}{d_1\,GW_{d_1f,1}\,GW_{s+d_2f,0}\ +\
d_2\,GW_{s+d_1f,1}\,GW_{d_2f,0}}.
\eest
But the invariants of ${\Bbb F}$ satisfy  $GW_{df,0}=GW_{s+df,0}=0$ for
$d\neq 0$ by the projection argument
      above, while  for $d\ne 0$ Lemma \ref{L.14.4} gives
$d_1GW_{d_1f,1}=GW_{d_1f,1}(s)=2\sigma(d_1)$.  We therefore get
\bear
H=2 F_0\cdot\l(G-{1\over 24}\r).
\label{rel2}
\eear

Combining (\ref{rel1})  with (\ref{rel2}) and noting that
$F_0(0)=GW_{s,0}=1\cdot t_s$  we see that $F_0$ satisfies the ODE
\best
t\, F_0'=12\,G\cdot F_0
\eest
with  $F_0(0)=1\cdot t_s$.  Hence
\best
F_0(t)=t_s\, \exp \left(12 \int G(t)/t\ dt \right).
\eest
Using the Taylor series of $\log(1-t)$ and some elementary combinatorics,
this becomes
\best
F_0(t)=t_s\, \left(\ma\prod_{d}{1 \over 1-t^d}\right)^{12}.
\eest

\bigskip

It remains to show (\ref{Fgformula}) for $g>0$.  This case is different
because
for  genus $g>0$ the relative invariants are no longer equal to the
absolute invariants. We start by fixing a fiber $F$ of  $E$ and introducing
two
generating functions  for the genus $g$ relative invariant: one recording
the number of curves passing through $g$ points in $E\setminus F$, the other
recording the number of curves passing through $g-1$ points in $E\setminus F$
plus a fixed point on $F$:
\best
F^V_g(f)=\sum_d \ov{GW}^F_{s+df,g}(p^{g}; C_1(f))t^d,\\
F^V_g(p)=\sum_d \ov{GW}^F_{s+df,g}(p^{g-1}; C_1(p))t^d.
\eest
Using Lemma \ref{L.avg=abs}, we can
relate the absolute and relative $g=1$ invariants  of $E$.
\begin{lemma} For $X=E$, the absolute and relative $g=1$ invariants in
the  classes $s+df\in H_2(E(1))$   are related by equations
\best
& (a) &  F_g \, =\,  F^V_g(p)+F^V_{g-1}(f)\cdot G' \\
& (b) &  F_g \, =\,  F^V_g(f)\\
& (c) &  0 \, =\,   F^V_g(p) \cdot F_0+ F_{g-1} \cdot F_1^V(p).
\eest
\label{L.split}
\end{lemma}
\pf To prove (a), we again write $E=E\#_F {\Bbb F}$ where ${\Bbb F}=
T^2\ti S^2$, and put $g-1$ points on $E$ and the remaining point on
${\Bbb F}$. If we start with a class $s+df$ the only possible
decompositions are $s+af$ and $s+bf$ where $d=a+b$. Since there are
$g-1$ points on the $E$ side, then the genus $g_1\ge g-1$. There are
two possibilities:
\begin{enumerate}
\item genus $g$ in class $s+df$ on $E$ and genus 0 in class $s+bf$
      on ${\Bbb F}$. But that forces $b=0$ so $a=d$.
\item genus $g-1$ in class $s+df$ on $E$ and genus 1 in class $s+bf$
on ${\Bbb F}$
\end{enumerate}
Putting then together gives (a). Relation (b) is a reformulation of
Lemma \ref{L.avg=abs}.

\medskip

      Relation (c) is seen by applying the symplectic sum formula to the
sum $K3= E\#_F E$ (the elliptic surface $K3=E(2)$ is the fiber sum of
$E=E(1)$ with itself). Because a generic complex structure on $K3$
admits no holomorphic curves, then all relative and absolute
invariants of K3 vanish. In particular, the genus $g$ invariants
through $g-1$ points in the class $[s+df]\in H_2(K3)/ {\cal R}$
vanish, where ${\cal R}$ is the set of rim tori corresponding to the
gluing $K3=E\#_F E$.

So, for any $g\ge 1$, put all the $g-1$ points on $X_1$ and split as above.
A dimension count
shows that the genus of the curve on $X_1$ must be at least $g-1$, so
the only possible decompositions are:
\begin{enumerate}
\item a genus $g$ curve in the class $s+d_1F$ on $X_1$  and a genus 0 curve
on $X_2$ in the class $s+d_2 f$, $d=d_1+d_2$;
\item a genus $g-1$ curve in the class $s+d_1F$ on $X_1$  and a genus 1 curve
on $X_2$ in the class $s+d_2 f$, $d=d_1+d_2$;
\end{enumerate}
The symplectic sum formula then gives $0 \, =\,   F^V_g(p) \cdot F_0+
F_{g-1}^V(f)\cdot F_1^V(p)$, which
simplifies by (b).
\qed

\bigskip

Formula (\ref{Fgformula}) follows quickly from Lemma
\ref{L.split}. Taking $g=1$ in Lemma \ref{L.split}d and factoring out
$F_0\ne 0$ yields $F_1^V(p)=0$.  Putting that in Lemma \ref{L.split}a
and again noting that $F_0\ne 0$ shows that $F_g^V(p)=0$ for all
$g>0$.  Parts (a) and (b) of Lemma \ref{L.split} then reduce to
\best
F_g=F_{g-1}\cdot G'
\eest
which gives (\ref{BLg}) by induction.

\vskip.4in




\setcounter{equation}{0}
\setcounter{section}{5}
\
{
\renewcommand{\theequation}{A.\arabic{equation}}
\renewcommand{\thetheorem}{A.\arabic{theorem}}
\section{Appendix -- Expansions of Relative TW Invariants}
\bigskip

     The Gromov-Witten invariants described in Section 1 are homology
elements --- the pushforward of
the compactified moduli space under (\ref{MHVSV}).  These can be
assembled into a power series
(\ref{defnRelInvt2}) with coefficients in homology.   Often, however, it is
convenient to write the GW and TW invariants as  power series   whose
coefficients are {\em
numbers}, preferably numbers with clear geometric interpretations.
This appendix describes
how that can be done for   the relative TW invariants which   appear in
the symplectic sum formula.

     Such series expansions are easiest when we can ignore the
complications caused by the covering
(\ref{HVcover}), replacing the space ${\cal H}_{X,A,s}^V$ by the more
easily understood space $V_s\cong
V^{\ell(s)}$.   That can be done by pushing the homology class of the
invariant down under the
projection $\ep$ of (\ref{HVcover}), obtaining  a `summed' GW series
\bear
\overline{GW}^V_X\ =\ \ep_*\left( GW^V_X\right)\ =\ \sum_{A\in H_2(X)}
\overline{GW}^V_{X,A}\ t_A
\label{A.1}
\label{def.avgGW}
\eear
whose coefficients are homology classes in $\sqcup_s V_s$.  This is
a less refined invariant, but has
the advantage that its coefficients become numbers after choosing a
basis of $H^*(V)$.

Of course (\ref{A.1}) is the  same as the original GW invariant when
the set ${\cal R}$ of (\ref{def.Rim})
vanishes, that is, there are no rim tori.  That occurs  whenever
$H_1(V)=0$ or more
generally when every rim tori represents zero in $H_2(X\setminus V)$.
We will describe the numerical
expansion under that assumption;  the same discussion applies to (\ref{A.1}).

\bigskip

When there are no rim tori  ${\cal H}_{X,A}^V$ is the union of those
$V_s\cong V^{\ell(s)}$ with
$\deg s=A\cdot V$.  Fix   a basis $\gamma_i$  of $H_*(V;\Q)$.
Then a basis for the tensor algebra on
${\Bbb N} \ti H_*(V)$  is given by  elements of the form
\bear
C_{s,I}=C_{s_1, \gamma_{i_1}}\otimes\dots \otimes  C_{s_\ell, \gamma_{i_\ell}}
\label{last.basis1}
\eear
where  $s_i\ge 1$ are integers. Let $\{C^*_{s,I}\}$ denote the
dual basis.  When  $\kappa\in H^*(\ov\M)$ and $\alpha \in {\Bbb T}(H^*(X))$,
we can expand
\bear
TW_X^V(\kappa,\al)=  \sum_{s,I} {1\over \ell(s)!}\;
TW_{X,A,\chi}^V(\kappa,\alpha; C_{s,I})\ \
C^*_{s,I}\ t_A\ \la^{-\chi}.
\label{2.last}
\label{A.3}
\eear

The coefficients in (\ref{A.3}) have a direct geometric
interpretation.  Choose generic pseudomanifolds
     $K\subset \ov{\M}_{g,n}$, $A_i\subset X$, and $\Gamma_j\subset V$
representing  the
Poincar\'{e} duals of  $\kappa$, $\al$, and the  $\gamma_j$ in their
respective spaces.  Then
$TW_{X,A,\chi}^V(\kappa,\alpha; C_{s,I})$ is the oriented number of genus $g$
$(J,\nu)$-holomorphic, $V$-regular
maps $f:C \to X$ with $C\in K$, $f(x_i)\in A_i$, and having a contact of
order $s_j$ with $V$ along  $\Gamma_j$.  Because of that
interpretation, the $C_{s,I}$ are called ``contact constraints''.

\bigskip

While for the analysis is important to work with {\em ordered} sequences $s$,
in applications it is more convenient to forget the ordering. The symmetries of
the GW invariants allow us to replace the basis (\ref{last.basis1})
with the one having elements of the form
\bear
{\bf C}_{m}=\prod_{a,i} \l(C_{a, \gamma_{i}}\r)^{m_{a,i}}
\label{last.basis2}
\eear
where  $m=(m_{a,i})$ is a finite sequence of nonnegative integers.
Generalizing (\ref{def.ls}), we write
\bear
     |m|=\prod_{a,i} a^{m_{a,i}}\qquad  m!=\prod_{a,i} m_{a,i}!\qquad
\ell(m)=\sum_{a,i} m_{a,i}\qquad \deg m= \sum_{a,i} a\cdot m_{a,i}.
\label{nr.m}
\eear
Let $\{z_{a,i}\}$ denote the dual basis; these generate a (super)
polynomial algebra with the relations
$z_{a,i}\, z_{b,j} =\pm\, z_{b,j}\, z_{a,i}$ where the sign is $+$
when $(\deg \gamma_i)(\deg \gamma_j)$
is even.  Then the generating series
of the relative $TW$ invariant is
\bear
TW_X^V(\kappa,\al)=  \sum_{A,g}\sum_{m}\;
TW_{X,A,\chi}^V(\kappa,\alpha; {\bf C}_{m})\; \; \prod_{a,i}
{(z_{a,i})^{m_{a,i}}\over m_{a,i}!}\;\; t_A\ \la^{-\chi}
\label{TW.gen.fcn}
\eear
where the sum is over all  sequences $m=(m_{a,i})$ as above and where
the coefficients
$TW_{X,A,\chi}^V(\kappa,\alpha; {\bf C}_{m})$ vanish unless
$ \deg m= A\cdot V$.  This  generating series  (\ref{TW.gen.fcn}) is
formally given by
\best
TW_X^V(\kappa,\al)=\sum_{A,g}
TW_{X,A,g}^V\l(\kappa,\alpha; \exp\l( {\sum_{a,i}
C_{a,\gamma_i}z_{a,i}}\r)\r)\;\; t_A\ \la^{-\chi}.
\eest
\label{rem.genTW}
\vskip.4in


\end{document}